\documentclass[12pt]{amsart}

\usepackage{amsfonts,amssymb,amscd}
\newcommand{\Q}{{\mathbb Q}}

\newcommand{\C}{{\mathbb C}}
\newcommand{\LL}{{\mathbb L}}

\newcommand{\ph}{{\varphi}}

\newcommand{\al}{{\alpha}}
\newcommand{\be}{{\beta}}
\newcommand{\la}{{\lambda}}
\newcommand{\ga}{{\gamma}}
\newcommand{\te}{{\theta}}
\newcommand{\ep}{{\varepsilon}}
\newcommand{\de}{{\delta}}
\newcommand{\ze}{{\zeta}}
\newcommand{\De}{{\Delta}}
\newcommand{\Si}{{\Sigma}}

\newcommand{\om}{{\omega}}

\newcommand{\ti}{\tilde}

\newcommand{\st}{\stackrel}
\newcommand{\lra}{\leftrightarrow}
\newcommand{\Lra}{\Leftrightarrow}
\newcommand{\Ra}{\Rightarrow}
\newcommand{\pd}[2]{\dfrac{\partial #1}{\partial #2}}
\newcommand{\ed}{{\hfill $\blacksquare$}}
\newcommand{\ee}{{\hfill $\blacktriangle$}}
\newcommand{\bt}{\blacktriangle}

\newcommand{\ol}[1]{ \overline{ #1 } }
\newcommand{\wh}[1]{ \widehat{ #1 } }

\newtheorem{theorem}{Theorem}[section]
\newtheorem{corollary}[theorem]{Corollary}
\newtheorem{proposition}[theorem]{Proposition}
\newtheorem{lemma}[theorem]{Lemma}
\newtheorem{claim}[theorem]{Claim}
\newtheorem{problem}[theorem]{Problem}

\theoremstyle{definition}
\newtheorem{definition}[theorem]{Definition}
\newtheorem{example}[theorem]{Example}

\title{ Compressed Drinfeld associators}

\author[kurlin]{V.~Kurlin$^*$}
\address{Institut de Math\'ematiques de Bourgogne,
 BP 47870, 21078 Dijon cedex, France}
\email{ kourline@topolog.u-bourgogne.fr, vak26@yandex.ru }

\thanks{$^*$
 The author was supported by a postdoctoral grant of the Regional
  Council of Burgundy (France).}
\subjclass[2000]{}
\keywords{Drinfeld associator, compressed associator, Kontsevich integral,
 zeta function, knot, hexagon equation, pentagon equation,
 Bernoulli numbers, extended Bernoulli numbers,
 Campbell-Baker-Hausdorff formula, Lie algebra, chord diagrams,
 Vassiliev invariants, compressed Vassiliev invariants}
\date{ this version: 18/10/2004\hspace{1cm} first version: 29/08/2004}

\begin{document}
\vspace*{-12mm}
\begin{center}
\it This is a preprint.
Comments and suggestions are welcome.
\end{center}
\vspace{5mm}

\begin{abstract}
Drinfeld associator is a key tool in computing the Kontsevich
 integral of knots.
A Drinfeld associator is a series in two
 non-commuting variables, satisfying highly complicated algebraic
 equations --- hexagon and pentagon.
The logarithm of a Drinfeld associator lives in the Lie algebra
 $L$ generated by the symbols $a,b,c$ modulo $[a,b]=[b,c]=[c,a]$.
The main result is a description of compressed
 associators that satisfy the compressed pentagon and hexagon
 in the quotient $L/\bigl[[L,L],[L,L]\bigr]$.
The key ingredient is an explicit form of
 Campbell-Baker-Hausdorff formula in the case
 when all commutators commute.
\end{abstract}

\maketitle
\vspace{-4mm}


\section{Introduction}


\subsection{Motivation and the previous results}

Let $A$ be a quasi-Hopf algebra \cite{Dr1} with a non-commutative
 non-associative coproduct $\De$.
Roughly speaking, \emph{an associator} is an element
 $\Phi\in A^{\otimes 3}$ controlling
 non-coassociativity of the coproduct $\Delta$.
Another element $R\in A^{\otimes 2}$ measures non-cocommutativity
 of $\Delta$.
For the representations of $A$ to form a tensor category,
 $R$ and $\Phi$ have to obey the so-called "pentagon" and
 "hexagon" equations.
V.~Drinfeld found a "universal" formula $(R_{KZ},\Phi_{KZ})$
 by using analytic methods --- differential equations and
 iterated integrals.
Also Drinfeld proved that there is an iterative algebraic
 procedure for finding a universal formula for an associator over
 the rationals.
Although this procedure is constructive, it does not
 give a close explicit formula.

The main motivation is the construction of the Kontsevich
 integral of knots via associators, investigated by T.~Q.~T.~Le,
 J.~Murakami \cite{LM3}, and D.~Bar-Natan \cite{BN2}.
Another combinatorial constructions of the universal Vassiliev
 invariant are in \cite{Piu, Car1}.
Recall that the Kontsevich integral takes values in the algebra
 $A$ of chord diagrams.
The LM-BN constructon gives an isotopy invariant of
 parenthesized framed tangles \cite{BN3} expressed via a Drinfeld
 associator that is a solution of rather complicated equations
 --- hexagon and pentagon (the same as mentioned above).
Any solution of these equations gives rise to a knot invariant.
Le and Murakami \cite{LM2} have proved that
 the resulting invariant is independent of a particularly choosen
 associator and coincides with the Kontsevich integral from \cite{Kon}
 provided $R=\exp(t^{12})$.
In other words, if one would know all coefficients of at least
 one associator, then one can calculate the whole Kontsevich integral
 for any knot.
Another approach of Bar-Natan, Le, and D.~Thurston has led to
 a formula for the Kontsevich integral of the unknot and all
 torus knots in the space of Jacobi diagrams \cite{BLT}.

One of non-even associators was expressed
 via multiple zeta values \cite{LM1}, i.e.
 via transcendental numbers.
Drinfeld computed the logarithm of the same
 associator in the case when all commutators commute
 by using classical zeta values \cite{Dr2}.
This result was the starting point of the present researches.
Bar-Natan calculated one of even rational Drinfeld associators up to degree
 7 in \cite{BN2}.
J.~Lieberum \cite{Lie} determined explicitly a rational even associator
 in a completion of the universal enveloping algebra of the Lie
 superalgebra $\mathrm{gl}(1|1)^{\otimes 3}$.
Up to now a close formula of a rational associator
 is stil unknown \cite[p.~433, Problem~3.13]{Oht}.

Extreme coefficients of all Drinfeld associators will be
 calculated in Theorem~1.5c below.
It turns out that they are rational and expressed via classical
 Bernoulli numbers $B_n$.


\subsection{Basic definitions}

\begin{definition}[associative algebra $A_n$,
 algebra of chord diagrams $A(X)$]
\noindent
\smallskip

\noindent
(a) For each $n\geq 2$, let the associative algebra $A_n$
 over the field $\C$ be generated by the symbols
 $t^{ij}=t^{ji}$ with $1\leq i\neq j\leq n$ and the relations
$$[t^{ij},t^{kl}]=0 \mbox{ if } i,j,k,l\mbox{ are pairwise disjoint},\quad
  [t^{ij},t^{jk}+t^{ki}]=0 \mbox{ if }i,j,k\mbox{ are pairwise disjoint},$$
 where \emph{the bracket} $[\, ,]:A_n\oplus A_n\to A_n$ is defined
 by $[a,b]:=ab-ba$.
Observe that the relations $[t^{ij},t^{jk}+t^{ki}]=0$ of $A_n$
 are equivalent to
$$[t^{ij},t^{jk}]=[t^{jk},t^{ki}]=[t^{ki},t^{ij}]
 \mbox{ for all pairwise disjoint }i,j,k\in\{1,\ldots,n\}.$$
The associative algebra $A_n$ is graded by \emph{the degree}
 defined by $\deg(t^{ij})=1$.
\smallskip

\noindent
(b) Let us define the same object $A_n$ geometrically.
Let $X$ be a 1-dimensional oriented compact manifold, possibly
 non-connected and with boundary.
\emph{A chord diagram on $X$} is a collection of non-oriented dashed lines
 (\emph{chords}) with endpoints on $X$.
Let $A(X)$ be the linear space generated by all chord diagrams on
 $X$ modulo \emph{the 4T relations}:

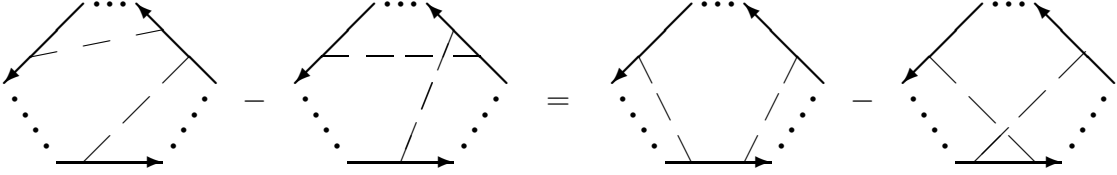
\begin{picture}(400,70)(-30,0)

{\thicklines
\put(0,0){\vector(1,0){40}}
\put(10,60){\vector(-1,-1){30}}
\put(60,30){\vector(-1,1){30}}
\put(70,20){$-$}
\put(110,0){\vector(1,0){40}}
\put(120,60){\vector(-1,-1){30}}
\put(170,30){\vector(-1,1){30}}
\put(185,20){$=$}
\put(230,0){\vector(1,0){40}}
\put(240,60){\vector(-1,-1){30}}
\put(290,30){\vector(-1,1){30}}
\put(300,20){$-$}
\put(340,0){\vector(1,0){40}}
\put(350,60){\vector(-1,-1){30}}
\put(400,30){\vector(-1,1){30}}
}

\multiput(15,60)(5,0){3}{\circle*{2}}
\multiput(-4,6)(-4,6){4}{\circle*{2}}
\multiput(44,6)(4,6){4}{\circle*{2}}

\put(-10,40){\line(5,1){10}}
\put(10,44){\line(5,1){10}}
\put(30,48){\line(5,1){10}}

\put(10,0){\line(1,1){10}}
\put(25,15){\line(1,1){10}}
\put(40,30){\line(1,1){10}}

\multiput(125,60)(5,0){3}{\circle*{2}}
\multiput(106,6)(-4,6){4}{\circle*{2}}
\multiput(154,6)(4,6){4}{\circle*{2}}

\put(100,40){\line(1,0){10}}
\put(117,40){\line(1,0){10}}
\put(132,40){\line(1,0){10}}
\put(150,40){\line(1,0){10}}

\put(130,0){\line(2,5){6}}
\put(138,20){\line(2,5){4}}
\put(144,35){\line(2,5){6}}



\multiput(245,60)(5,0){3}{\circle*{2}}
\multiput(226,6)(-4,6){4}{\circle*{2}}
\multiput(274,6)(4,6){4}{\circle*{2}}

\put(240,0){\line(-1,2){5}}
\put(232,16){\line(-1,2){5}}
\put(225,30){\line(-1,2){5}}

\put(260,0){\line(1,2){5}}
\put(268,16){\line(1,2){5}}
\put(275,30){\line(1,2){5}}

\multiput(355,60)(5,0){3}{\circle*{2}}
\multiput(336,6)(-4,6){4}{\circle*{2}}
\multiput(384,6)(4,6){4}{\circle*{2}}

\put(347,0){\line(1,1){13}}
\put(364,17){\line(1,1){10}}
\put(379,32){\line(1,1){10}}

\put(370,0){\line(-1,1){10}}
\put(355,15){\line(-1,1){10}}
\put(340,30){\line(-1,1){10}}

\end{picture}
\vspace{2mm}

The dotted arcs represent parts of the diagrams that are not
 shown in the figure.
These parts are assumed to be the same in all four diagrams.

If $X=X_n$ is the disjoint union of $n$ oriented segments
 (\emph{strands}), then $A(X_n)$ can be equipped with a natural product.
If in the definition of $A(X_n)$ one allows only
 horizontal chords with endpoints on $n$ vertical strands,
 then the resulting algebra $A^{hor}(X_n)$ is isomorphic
 to the algebra $A_n$.
Indeed, thinking of $t^{ij}$ as a horizontal chord connecting the
 $i$th and $j$th vertical strands, the relations between the $t^{ij}$
 become the 4T relations:
\begin{picture}(440,40)(-20,0)

{\thicklines
\put(0,10){$[t^{12},t^{23}]=$}
\put(55,30){\vector(0,-1){30}}
\put(80,30){\vector(0,-1){30}}
\put(105,30){\vector(0,-1){30}}
\put(115,10){$-$}
\put(130,30){\vector(0,-1){30}}
\put(155,30){\vector(0,-1){30}}
\put(180,30){\vector(0,-1){30}}
\put(190,10){$=[t^{23},t^{13}]=$}
\put(260,30){\vector(0,-1){30}}
\put(285,30){\vector(0,-1){30}}
\put(310,30){\vector(0,-1){30}}
\put(320,10){$-$}
\put(335,30){\vector(0,-1){30}}
\put(360,30){\vector(0,-1){30}}
\put(385,30){\vector(0,-1){30}}
\put(440,0){$\blacksquare$}
}

\put(55,25){\line(1,0){10}}
\put(70,25){\line(1,0){10}}
\put(80,15){\line(1,0){10}}
\put(95,15){\line(1,0){10}}
\put(55,25){\circle*{3}}
\put(80,25){\circle*{3}}
\put(80,15){\circle*{3}}
\put(105,15){\circle*{3}}

\put(130,15){\line(1,0){10}}
\put(145,15){\line(1,0){10}}
\put(155,25){\line(1,0){10}}
\put(170,25){\line(1,0){10}}
\put(130,15){\circle*{3}}
\put(155,15){\circle*{3}}
\put(155,25){\circle*{3}}
\put(180,25){\circle*{3}}

\put(285,25){\line(1,0){10}}
\put(300,25){\line(1,0){10}}
\put(260,15){\line(1,0){8}}
\put(273,15){\line(1,0){8}}
\put(288,15){\line(1,0){8}}
\put(302,15){\line(1,0){8}}
\put(260,15){\circle*{3}}
\put(285,25){\circle*{3}}
\put(310,15){\circle*{3}}
\put(310,25){\circle*{3}}

\put(360,15){\line(1,0){10}}
\put(375,15){\line(1,0){10}}
\put(335,25){\line(1,0){8}}
\put(348,25){\line(1,0){8}}
\put(364,25){\line(1,0){8}}
\put(377,25){\line(1,0){8}}
\put(335,25){\circle*{3}}
\put(360,15){\circle*{3}}
\put(385,15){\circle*{3}}
\put(385,25){\circle*{3}}

\end{picture}
\end{definition}

\begin{definition}[Lie algebra $L_n$, quotient $\bar L_n$,
 long commutators ${[a_1\ldots a_k]}$]
\noindent
\smallskip

\noindent
(a) The Lie algebra $L_n$ is generated by the same generators and
 relations as the associative algebra $A_n$ of Definition~1.1.
The Lie algebra $L_n$ is graded  with respect to $\deg(t^{ij})=1$.
\smallskip

\noindent
(b)
By $[L_n,L_n]$ denote the Lie subalgebra of $L_n$,
 generated by all commutators $[a,b]$ with $a,b\in L_n$.
Introduce \emph{the compressed quotient}
 $\bar L_n=L_n/\bigl[[L_n,L_n],[L_n,L_n]\bigr]$.
Let $\hat A_n$, $\hat L_n$, and $\wh{\bar L}_n$ be the algebras of
 formal series of elements from $A_n$, $L_n$, and $\bar L_n$, respectively.
\smallskip

\noindent
(c) For elements $a_1,\ldots,a_n$ of a Lie algebra $L$, set
 $[a_1a_2\ldots a_k]=[a_1,[a_2,[\ldots,a_k]\ldots]]$.
For example, the algebras $\hat A_2$ and $\hat L_3$
 contain the series $\exp(t^{12})$ and
 $\sum\limits_{k=1}^{\infty}[(t^{12})^kt^{23}]$, resp.
\ed
\end{definition}

\begin{definition}[operators $\De_k$ and $\ep_k$,
 Drinfeld associators and compressed associators]
\noindent
\smallskip

\noindent
(a)
Let $t^{ij}$ be the generators of $A_n$.
Let $\De_k:A_n\to A_{n+1}$ for $0\leq k\leq n+1$
 and $\ep_k:A_n\to A_{n-1}$ for $1\leq k\leq n$
 be the algebra morphisms defined by their action on $t^{ij}$
 (here $i<j$):
$$\De_k(t^{ij})=\left\{ \begin{array}{lll}
t^{ij},& \mbox{if} & i<j<k,\\
t^{i,j+1},& \mbox{if} & i<k<j,\\
t^{i+1,j+1},& \mbox{if} & k<i<j,\\
t^{ij}+t^{i,j+1},& \mbox{if} & i<j=k,\\
t^{i,j+1}+t^{i+1,j+1},& \mbox{if} & i=k<j;
\end{array}\right.\qquad
\ep_k(t^{ij})=\left\{ \begin{array}{lll}
t^{ij},& \mbox{if} & i<j<k,\\
t^{i,j-1},& \mbox{if} & i<k<j,\\
t^{i-1,j-1},& \mbox{if} & k<i<j,\\
0,& \mbox{if} & i<j=k,\\
0,& \mbox{if} & i=k<j.
\end{array}\right.$$
$\De_0$ ($\De_{n+1}$) acts by adding a strand on
 the left (right), $\De_k$ for $1\leq k\leq n$ acts by doubling
 the $i$th strand and summing up all the possible ways of lifting
 the chords that were connected to the $i$th strand to the two
 daughter strands.
The operator
$\ep_k$ acts by deleting the $i$th strand and
 mapping the chord diagram to 0, if any chord in it was connected
 to the $i$th strand.
\smallskip

\noindent
(b) \emph{A horizontal Drinfeld associator}
 (briefly, \emph{a Drinfeld associator}) is
 an element $\Phi\in\hat A_3$ satisfying
 the following equations (here set
 $\Phi^{ijk}:=\Phi(t^{ij},t^{jk})$ and $\Phi:=\Phi^{123}$)
$$\begin{array}{llr}
\mbox{(symmetry)} &
 \Phi\cdot\Phi^{321}=1\mbox{ in }\hat A_3, & (1.3a)\\
\mbox{(hexagon)} &
 \De_1(\exp(t^{12})) = \Phi^{312}\cdot \exp(t^{13})\cdot (\Phi^{-1})^{132}
  \cdot \exp(t^{23})\cdot \Phi^{123}\mbox{ in }\hat A_3, & (1.3b)\\
\mbox{(pentagon)} &
 \De_0(\Phi)\cdot \De_2(\Phi)\cdot\De_4(\Phi)
  =\De_3(\Phi)\cdot \De_1(\Phi)\mbox{ in }\hat A_4, & (1.3c)\\
\mbox{(non-degeneracy)} &
 \ep_1\Phi=\ep_2\Phi=\ep_3\Phi=1\mbox{ in } \hat A_2, & (1.3d)\\
\mbox{(group-like)} &
 \Phi=\exp(\ph)\mbox{ in }\hat A_3\mbox{ for some element }
 \ph\in\hat L_3. & (1.3e)
\end{array}$$
A geometric interpretation of the hexagon and pentagon is shown below:

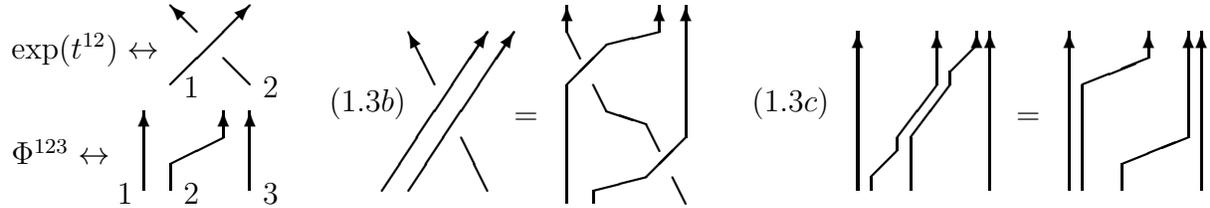
\begin{picture}(450,85)(0,0)

{\thicklines
\put(0,60){$\exp(t^{12})\leftrightarrow$}
\put(60,50){\vector(1,1){30}}
\put(70,70){\vector(-1,1){10}}
\put(90,50){\line(-1,1){10}}
\put(65,45){$1$}
\put(95,45){$2$}

\put(0,20){$\Phi^{123}\leftrightarrow$}
\put(50,10){\vector(0,1){30}}
\put(60,10){\line(0,1){10}}
\put(60,20){\line(2,1){20}}
\put(80,30){\vector(0,1){10}}
\put(90,10){\vector(0,1){30}}
\put(40,5){$1$}
\put(65,5){$2$}
\put(95,5){$3$}

\put(120,40){$(1.3b)$}
\put(140,10){\vector(2,3){40}}
\put(150,10){\vector(2,3){40}}
\put(160,50){\vector(-1,2){10}}
\put(180,10){\line(-1,2){10}}
\put(190,35){$=$}
\put(210,5){\line(0,1){45}}
\put(210,50){\line(1,1){15}}
\put(225,65){\line(4,1){20}}
\put(245,70){\vector(0,1){10}}
\put(220,5){\line(0,1){5}}
\put(220,10){\line(4,1){20}}
\put(240,15){\line(1,1){15}}
\put(255,30){\vector(0,1){50}}
\put(255,5){\line(-1,2){5}}
\put(245,25){\line(-1,2){5}}
\put(240,35){\line(-3,1){15}}
\put(225,40){\line(-1,2){5}}
\put(215,60){\line(-1,2){5}}
\put(210,70){\vector(0,1){10}}

\put(280,40){$(1.3c)$}
\put(320,10){\vector(0,1){60}}
\put(325,10){\line(0,1){5}}
\put(325,15){\line(1,1){10}}
\put(335,25){\line(0,1){5}}
\put(335,30){\line(3,4){15}}
\put(350,50){\vector(0,1){20}}
\put(340,10){\line(0,1){20}}
\put(340,30){\line(3,4){15}}
\put(355,50){\line(0,1){5}}
\put(355,55){\line(1,1){10}}
\put(365,65){\vector(0,1){5}}
\put(370,10){\vector(0,1){60}}
\put(380,35){$=$}
\put(400,10){\vector(0,1){60}}
\put(405,10){\line(0,1){40}}
\put(405,50){\line(5,2){25}}
\put(430,60){\vector(0,1){10}}
\put(420,10){\line(0,1){10}}
\put(420,20){\line(5,2){25}}
\put(445,30){\vector(0,1){40}}
\put(450,10){\vector(0,1){60}}
}
\end{picture}

\noindent
(c) If an associator $\Phi\in\hat A_3$ vanishes in all odd degrees,
 then $\Phi$ is said to be \emph{even}.
Note that the symmetry (1.3a) implies
 $\ph(a,b)=-\ph(b,a)$ in $\hat L_3$.
By taking the logarithms of (1.3b) and (1.3c)
 and projecting them under $\hat L_3\to\wh{\bar L}_3$
 and $\hat L_4\to\wh{\bar L}_4$ one gets
 \emph{the compressed hexagon~$\ol{(1.3b)}$
 and pentagon~$\ol{(1.3c)}$}, respectively.
\emph{A compressed associator} $\bar\ph\in\wh{\bar L}_3$
 is a solution of $\ol{(1.3b)}$ and $\ol{(1.3c)}$,
 satisfying $\bar\ph(a,b)=-\bar\ph(b,a)$ and
 $\bar\ph(a,0)=\bar\ph(0,b)=0$.
\ed
\end{definition}

Definition~1.3b uses a non-classical normalization.
Drinfeld considered the two hexagons \cite{Dr1}:
$\De_1\left(\exp\left(\pm\dfrac{t^{12}}{2}\right)\right)
 = \Phi^{312}\cdot \exp\left(\pm\dfrac{t^{13}}{2}\right)\cdot
 (\Phi^{-1})^{132} \cdot \exp\left(\pm\dfrac{t^{23}}{2}\right)
 \cdot \Phi^{123}.$
To avoid huge denominators in future
 the change of the variables $t^{ij}\mapsto 2t^{ij}$ was made.
Moreover, Bar-Natan has proved that both above hexagons
 are equivalent to the positive hexagon (with the sign "+")
 and the symmetry~(1.3a), see \cite[Proposition~3.7]{BN2}.
The logarithm $\ph=\log(\Phi)$ of any Drinfeld associator
 projects under $\hat L_3\to \hat{\bar L}_3$ onto a compressed
 associator.

\begin{example}
Bar-Natan calculated\footnote{
 By the above normalization one needs to divide
  the denominators at all terms of the degree $n$ by $2^n$.}
 one of even Drinfeld associators up to degree 7
$$\ph^B(a,b)=\left(\frac{[ab]}{12}
 -\frac{8[a^3b]+[abab]}{720}
 +\frac{96[a^5b]+4[a^3bab]+65[a^2b^2ab]+68[aba^3b]+4[(ab)^3]}{90720}\right)-$$
 $-$(interchange of $a\lra b$), where $a=t^{12}$, $b=t^{23}$.
Degree 7 is the maximal achievement of Bar-Natan's computer programme.
Then in $\hat{\bar L}_3$ one gets:
 $\quad\bar\ph^B(a,b)=\dfrac{[ab]}{6}-$
$$-\frac{4[a^3b]+[abab]+4[b^2ab]}{360}
  +\frac{[a^5b]+[b^4ab]}{945}
  +\frac{[a^3bab]+[ab^3ab]}{1260}
  +\frac{23[a^2b^2ab]}{30240}+\cdots\quad(1.4)\eqno{\bt}$$
\end{example}
\smallskip


\subsection{Main results}

For a series $f(\la,\mu)$, let us introduce its \emph{even} and
 \emph{odd} parts:
$$Even\bigl(f(\la,\mu)\bigr)=\dfrac{f(\la,\mu)+f(-\la,-\mu)}{2},\quad
   Odd\bigl(f(\la,\mu)\bigr)=\dfrac{f(\la,\mu)-f(-\la,-\mu)}{2}.$$
A simple reformulation of Theorem~1.5 will be given in
 Corollary~1.6c.

\begin{theorem}
(a) Any compressed Drinfeld associator $\bar\ph\in\wh{\bar L}_3$
 from Definition~1.3c is
$$\bar\ph(a,b)=\sum\limits_{k,l\geq 0} \al_{kl}[a^kb^lab],\;
 \mbox{ where } a=t^{12},\; b=t^{23},\quad
 \al_{kl}\in\C. \leqno{(1.5a)}$$
Moreover, for all $k,l\geq 1$, the coefficients $\al_{kl}$
 are symmetric ($\al_{kl}=\al_{lk}$) and
 could be expressed linearly in terms of $\al_{ij}$ with $i+j<k+l$.
\medskip

\noindent
(b) Let $f(\la,\mu)=\sum\limits_{k,l\geq 0}\al_{kl}\la^k\mu^l$ be
 the generating function of the coefficients $\al_{kl}=\al_{lk}$.
Then the compressed hexagon~$\ol{(1.3b)}$ from Definition~1.3c
 is equivalent to the equation
$$f(\la,\mu)+e^{\mu}f(\mu,-\la-\mu)+e^{-\la}f(\la,-\la-\mu)
=\dfrac{1}{\la+\mu}\left(
 \dfrac{e^{\mu}-1}{\mu}+\dfrac{e^{-\la}-1}{\la}
 \right).\leqno{(1.5b)}$$
Moreover, the compressed pentagon~$\ol{(1.3c)}$
 for $\bar\ph\in\wh{\bar L}_3$
 follows from the symmetry $\al_{kl}=\al_{lk}$.
\medskip

\noindent
(c) The general solution of (1.5b) is
 $f(\la,\mu)=Even(f(\la,\mu))+Odd(f(\la,\mu))$, where
$$\left\{ \begin{array}{l}
1+\la\mu\cdot Even\bigl(f(\la,\mu)\bigr)
=\dfrac{e^{\la+\mu}-e^{-\la-\mu}}{2(\la+\mu)} \left(
 \dfrac{2\om}{e^{\om}-e^{-\om}}
 +\sum\limits_{n=3}^{\infty} h_n(\la,\mu)\right),\\ \\
Odd\bigl(f(\la,\mu)\bigr)
=\dfrac{e^{\la+\mu}-e^{-\la-\mu}}{2} \left(
 \sum\limits_{n=0}^{\infty} \ti\be_{n0} \om^{2n}
 +\sum\limits_{n=3}^{\infty} \ti h_n(\la,\mu)\right),
 \end{array} \right.
 \leqno{(1.5c)}$$
$$h_n(\la,\mu)
 =\sum_{k=1}^{[\frac{n}{3}]}\be_{nk}\la^{2k}\mu^{2k}(\la+\mu)^{2k}
  \om^{2n-6k}\mbox{ for }n\geq 3,\quad
 \be_{nk}\in\C,\quad \om=\sqrt{\la^2+\la\mu+\mu^2}.$$
In particular, any honest Drinfeld associator has the extreme coefficients
 $\al_{2k,0}=\dfrac{2^{2k+1}B_{2k+2}}{(2k+2)!}$ for every $k\geq 0$,
 where $B_n$ are Bernoulli numbers.
The polynomials $\ti h_n(\la,\mu)$ are defined by the same
 formula as $h_n(\la,\mu)$, except the coefficients
 $\ti\be_{nk}\in\C$ are substituted for $\be_{nk}$.
The coefficients $\ti\be_{n0}$ (for $n\geq 0$),
 $\be_{nk}$, and $\ti\be_{nk}$ are free parameters for
 $1\leq k\leq[\frac{n}{3}]$, $n\geq 3$.
\end{theorem}

By Theorem~1.5c the differences of all compressed associators form
 a linear space generated by $\be_{nk},\ti\be_{nk}$.
The projection of any Drinfeld associator is in (1.5c), see
 Problem~6.10a.

\begin{corollary}
(a,b) There are two distinguished even compressed Drinfeld associators\\
 $\bar\ph=\sum\limits_{k,l\geq 0}\al_{kl}[a^kb^lab]\in\wh{\bar L}_3$
 defined by the function
 $f(\la,\mu)=\sum\limits_{k,l\geq 0}\al_{kl}\la^k\mu^l$
 as follows

$$\begin{array}{ll}
(1.6a)\;\mbox{the first series:} &
1+\la\mu f^I(\la,\mu)
 =\dfrac{e^{\la+\mu}-e^{-\la-\mu}}{2(\la+\mu)}
  \cdot \dfrac{2\om}{e^{\om}-e^{-\om}}\,,\;
 \om=\sqrt{\la^2+\la\mu+\mu^2}; \\ \\
(1.6b)\;\mbox{the second series:} &
1+2\la\mu f^{II}(\la,\mu)
 =\dfrac{e^{\la+\mu}-e^{-\la-\mu}}{2(\la+\mu)}
 \left( \dfrac{2\la}{e^{\la}-e^{-\la}}
  +\dfrac{2\mu}{e^{\mu}-e^{-\mu}} -1 \right).
\end{array}$$

\noindent
(c) Any compressed Drinfeld associator $\bar\ph\in \hat{\bar L}_3$
 can be defined by the Drinfeld series
$$1+\la\mu f^D(\la,\mu)=\exp\left(\sum\limits_{n=2}^{\infty}
 \dfrac{\ze(n)}{n}\cdot
 \dfrac{\la^n+\mu^n-(\la+\mu)^n}{(\pi\sqrt{-1})^n} \right),$$
 where after rewriting the above exponent as a series in $\la,\mu$
 one substitutes a free parameter for each monomial
 consisting of odd zeta numbers like $\ze(3)^p\ze(5)^q\ze(7)^r\ldots$
 (see Definition~6.1).

In particular, one gets the third compressed associator
$$(1.6c)\;\mbox{the third series:}\quad
1+\la\mu f^{III}(\la,\mu)=\exp\Biggl(\sum\limits_{n=1}^{\infty}
 \dfrac{2^{2n}B_{2n}}{4n(2n)!}
 \Bigl( (\la+\mu)^{2n}-\la^{2n}-\mu^{2n} \Bigr) \Biggr).$$
\end{corollary}


\subsection{Scheme of proofs}

Key points of proofs are listed below.
\smallskip

\noindent
\emph{First key point:}
 a behavior of the Bernoulli numbers ($B_n=0$ for each odd $n\geq 3$).
\smallskip

\noindent
\emph{Second key point:}
 the Bernoulli numbers $B_n$ can be extended in a natural way,
 this extension gives a compressed variant of CBH formula
 (Definition~2.4 and Proposition~2.8).
\smallskip

\noindent
\emph{Third key point:}
 properties of the extended Bernoulli
 numbers $C_{mn}$ and their generating function $C(\la,\mu)$:
 a non-trivial symmetry $C(\la,\mu)=C(-\mu,-\la)$ (Lemma~2.10)
 and an explicit expression of $C(\la,\mu)$ (Proposition~2.12).
\smallskip

\noindent
\emph{4th key point:}
 the original hexagon equation~(1.3b) can be simplified
 in such a way that it remains to apply CBH formula
 in an essential way exactly once (Lemma~3.1).
\smallskip

\noindent
\emph{5th key point:}
 the quotient $\bar L_3=L_3/\bigl[[L_3,L_3],[L_3,L_3]\bigr]$, where
 a compressed associator lives,
 is isomorphic to a Lie algebra with a small basis of commutators
 (Proposition~3.4).
\smallskip

\noindent
\emph{6th key point:}
 the compressed hexagon equation~$\ol{(1.3b)}$ is equvalent to
 a recursive linear system for the coefficients $\al_{kl}$
 (Proposition~3.9 and Lemma~4.1).
\smallskip

\noindent
\emph{7th key point:}
 the extreme coefficients $\al_{2k,0}$ of the exact logarithm
 of any Drinfeld associator (not only the compressed one)
 are expressed via Bernoulli numbers $B_{2n}$ (Lemma~4.2).
\smallskip

\noindent
\emph{8th key point:}
 for any compressed associator,
 the compressed hexagon~$\ol{(1.3b)}$ can be split into two equations
 for the even and odd parts of this associator (Lemma~4.5).
\smallskip

\noindent
\emph{9th key point:}
 up to certain factor the general solution of
 the compressed hexagon~$\ol{(1.3b)}$ is
 a series $h(\la,\mu)$ with the symmetry
 $h(\la,\mu)=h(\la,-\la-\mu)$ (Lemmas~4.6 and 4.12).
\smallskip

\noindent
\emph{10th key point:}
 non-uniqueness of compressed associators
  is closely related with non-uniqueness of
  associator polynomials (Definition~4.7 and Lemma~4.9).
\smallskip

\noindent
\emph{11th key point:}
 all associator polynomials can be described explicitly:
 in each degree $2n$ the family of all associator
 polynomials depends on $[\frac{n}{3}]$ free parameters
 (Proposition~4.10).
\smallskip

\noindent
\emph{12th key point:}
 for any compressed associator, the compressed pentagon
 equation~$\ol{(1.3c)}$ follows
 from the symmetry condition $\al_{kl}=\al_{lk}$
 (Proposition~5.10).
\medskip

\noindent
All the above results are the positive ones.
But there is also the negative one.
\smallskip

\noindent
\emph{13th key point:} the Drinfeld series
 (a compressed associator expressed via zeta values)
 does not lead to non-trivial polynomial relations between
 odd zeta values (Proposition~6.9).

\noindent
The 13th key point means that odd zeta values are too complicated numbers.
\smallskip

The paper is organized as follows.
In Section~2 the extended Bernoulli numbers $C_{mn}$ are introduced,
 one deduces a compressed variant of CBH formula in the case
 when all commutators commute with each other.
In Sections~3 and 4 the compressed hexagon is written explicitly and solved.
Section~5 is devoted to checking the compressed pentagon.
In Section~6 one explains, why the Drinfeld series does not lead
 to polynomial relations between odd zeta values.
Theorem~1.5a is proved at the end of Subsection~3.3.
The hexagon and pentagon parts of Theorem~1.5b are checked
 in Subsections~4.1 and 5.3, respectively.
Theorem~1.5c, Corollaries~1.6a-b and 1.6c are verified in
 Subsections~4.3 and 6.2, respectively.
In Subsection~6.3 open problems and suggestions for future
 researches are formulated.
Moreover, Appendix contains a lot of explicit formulae
 discussed in the paper in their general forms.

The following diagram shows the scheme for the proof of Theorem~1.5.
The most important steps are called propositions, they are of
 independent interest, especially Propositions~2.8 and 2.12 together.
Lemmas are of less importance.
Claims are technical assertions.

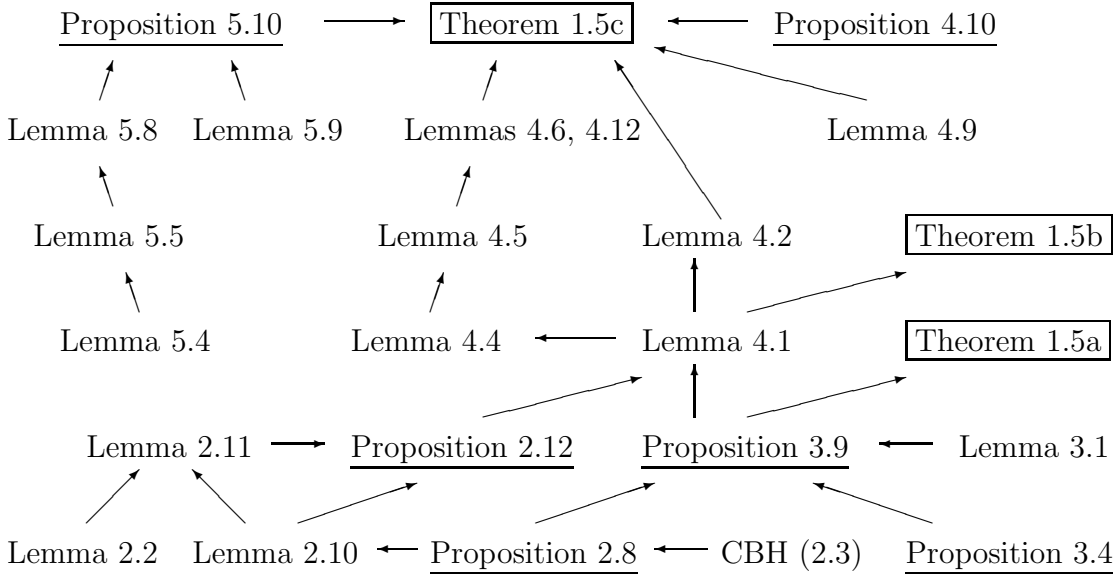
\begin{picture}(400,230)(-10,10)

\put(160,220){\fbox{Theorem~1.5c}}
\put(290,220){\underline{Proposition~4.10}}
\put(20,220){\underline{Proposition~5.10}}

\put(35,195){\vector(1,3){5}}
\put(90,195){\vector(-1,3){5}}
\put(180,195){\vector(1,3){5}}
\put(325,195){\vector(-4,1){80}}
\put(280,225){\vector(-1,0){30}}
\put(120,225){\vector(1,0){30}}

\put(0,180){Lemma~5.8}
\put(70,180){Lemma~5.9}
\put(150,180){Lemmas~4.6,~4.12}
\put(310,180){Lemma~4.9}

\put(40,155){\vector(-1,3){5}}
\put(170,155){\vector(1,3){5}}
\put(270,150){\vector(-2,3){40}}

\put(10,140){Lemma~5.5}
\put(140,140){Lemma~4.5}
\put(240,140){Lemma~4.2}
\put(340,140){\fbox{Theorem~1.5b}}

\put(50,115){\vector(-1,3){5}}
\put(160,115){\vector(1,3){5}}
\put(260,115){\vector(0,1){20}}
\put(280,115){\vector(4,1){60}}

\put(20,100){Lemma~5.4}
\put(130,100){Lemma~4.4}
\put(240,100){Lemma~4.1}
\put(340,100){\fbox{Theorem~1.5a}}

\put(180,75){\vector(4,1){60}}
\put(260,75){\vector(0,1){20}}
\put(280,75){\vector(4,1){60}}
\put(230,105){\vector(-1,0){30}}

\put(30,60){Lemma~2.11}
\put(130,60){\underline{Proposition~2.12}}
\put(240,60){\underline{Proposition~3.9}}
\put(360,60){Lemma~3.1}

\put(100,65){\vector(1,0){20}}
\put(350,65){\vector(-1,0){20}}
\put(30,35){\vector(1,1){20}}
\put(90,35){\vector(-1,1){20}}
\put(110,35){\vector(3,1){45}}
\put(200,35){\vector(3,1){45}}
\put(350,35){\vector(-3,1){45}}

\put(0,20){Lemma~2.2}
\put(70,20){Lemma~2.10}
\put(160,20){\underline{Proposition~2.8}}
\put(270,20){CBH~(2.3)}
\put(340,20){\underline{Proposition~3.4}}

\put(155,25){\vector(-1,0){15}}
\put(265,25){\vector(-1,0){20}}

\end{picture}
\vspace{2mm}

\noindent
{\bf Acknowledgement.}
The author is grateful to O.~Popov for his help in the proof of
 Proposition~3.4 and to D.~Moskovich for his useful suggestions.
Also he thanks D.~Bar-Natan, L.~Funar, C.~Kassel, G.~Masbaum,
 and P.~Vogel for their attention and comments.


\section{Campbell-Baker-Hausdorff formula (CBH)}

This section is devoted to an explicit form of CBH formula in
 the case when all commutators commute with each other,
 see Propositions~2.8 and 2.12.


\subsection{Classical recursive CBH formula}

Recall a classical recursive CBH formula (Theorem~2.3)
 originally proved by J.~Campbell \cite{Cam}, H.~Baker \cite{Bak},
 and F.~Hausdorff \cite{Hau}.

\begin{definition}[Hausdorff series $H$, Bernoulli numbers $B_n$,
 derivative $D=H_1\pd{}{Q}\;$]
\noindent
\smallskip

\noindent
(a) Let $L$ be the free Lie algebra generated by the symbols $P,Q$.
By $\hat L$ denote the algebra of formal series of elements
 from $L$.
\emph{The Hausdorff series} is
 $H=\log\bigl(\exp(P)\cdot\exp(Q)\bigr)\in\hat L$.
\smallskip

\noindent
(b) \emph{The classical Bernoulli numbers} $B_n$ are defined by
 the exponential generating function:
 $\sum\limits_{n=0}^{\infty}\dfrac{B_n}{n!}x^n=\dfrac{x}{e^x-1}$.
One can verify that $B_0=1$, $B_1=-\dfrac{1}{2}$ and that
 $\dfrac{x}{e^x-1}+\dfrac{x}{2}$ is an even function, which shows
 that $B_n=0$ for odd $n\geq 3$.
\emph{The first key point}: the function $\dfrac{x}{e^x-1}$
 vanishes in almost all odd degrees, while $\dfrac{e^x-1}{x}$ does not.
The first few Bernoulli numbers are
 $B_0=1$, $B_1=-\dfrac{1}{2}$, $B_2=\dfrac{1}{6}$, $B_3=0$,
 $B_4=-\dfrac{1}{30}$, $B_5=0$, $B_6=\dfrac{1}{42}$, $B_7=0$.
\medskip

\noindent
(c) \emph{A derivative} of a Lie algebra $L$ is a linear function
 $D:\hat L\to\hat L$ such that $D[x,y]=[Dx,y]+[x,Dy]$ for all
 $x,y\in\hat L$.
For an element $H_1\in\hat L$, denote by $D=H_1\pd{}{Q}$
 the derivative of $L$, which maps $P$ onto $0$ and $Q$ onto $H_1$.
\ed
\end{definition}

\emph{The first key point} implies recursive formulae
 for the Bernoulli numbers $B_n$.

\begin{lemma}
For each $m\geq 1$, the Bernoulli numbers $B_n$ satisfy the relations
$$(a)\; \sum\limits_{n=1}^m\binom{m+1}{n} B_n=-1,\;
  (b)\; \sum\limits_{k=1}^{[\frac{m}{2}]}\binom{m+1}{2k} B_{2k}
   =\dfrac{m-1}{2},\;
  (c)\; \sum\limits_{n=1}^m (-1)^n\binom{m+1}{n} B_n=m.$$
In particular, the first four relations from the item (a) are
$$2B_1=-1,\quad 3B_1+3B_2=-1,\quad 4B_1+6B_2+4B_3=-1,\quad
  5B_1+10B_2+10B_3+5B_4=-1.$$
\end{lemma}
\begin{proof}
(a) One obtains:
$$1=\dfrac{e^x-1}{x}\cdot\dfrac{x}{e^x-1}
   =\left(\sum\limits_{k=1}^{\infty}\dfrac{x^{k-1}}{k!}\right)\cdot
    \left(\sum\limits_{n=0}^{\infty}\dfrac{B_n}{n!}x^n\right)
   =1+\sum_{m=1}^{\infty}\left(\dfrac{1}{(m+1)!}
    +\sum\limits_{\binom{k+n=m+1}{k,n\geq 1}}
    \dfrac{1}{k!}\cdot\dfrac{B_n}{n!}\right)x^m,$$
 i.e.
 $\sum\limits_{n=1}^m\dfrac{1}{(m-n+1)!}\cdot\dfrac{B_n}{n!}
  =-\dfrac{1}{(m+1)!}$ as required.
\smallskip

\noindent
(b) Since $B_{2k+1}=0$ for every $k\geq 1$, then
 $-1=\sum\limits_{n=1}^m\binom{m+1}{n} B_n
 =(m+1)B_1+\sum\limits_{k=1}^{[\frac{m}{2}]}\binom{m+1}{2k} B_{2k}$,
 hence
 $\sum\limits_{k=1}^{[\frac{m}{2}]}\binom{m+1}{2k} B_{2k}
 =-1+\dfrac{1}{2}(m+1)=\dfrac{m-1}{2}$.
\smallskip

\noindent
(c) The required formula is equivalent to the item~(b),
 since $B_{2k+1}=0$ for each $k\geq 1$.
\end{proof}

The following theorem is quoted from
 \cite[Corollaries 3.24--3.25, p.~77--79]{Reu}.

\begin{theorem} \cite{Cam, Bak, Hau} 
The Hausdorff series $H=\log\bigl(\exp(P)\cdot\exp(Q)\bigr)$ is
 $H=\sum\limits_{m=0}^{\infty} H_m$,\\
 $H_0=Q$,
 $H_1=P-\dfrac{1}{2}[Q,P]
  +\sum\limits_{n=1}^{\infty}\dfrac{B_{2n}}{(2n)!}[Q^{2n}P]$,
 $\;H_m=\dfrac{1}{m}\left(H_1\pd{}{Q}\right)(H_{m-1})$ for $m\geq 2$.
\qed
\end{theorem}
\smallskip


\subsection{Extended Bernoulli numbers $C_{mn}$ and
 compressed variant of CBH formula}

\emph{The compressed variant} is the case
 when all commutators commute.
To get an explicit form of CBH formula in this setting,
 one needs extended Bernoulli numbers.

\begin{definition}[extended Bernoulli numbers $C_{mn}$,
 generating function $C(\la,\mu)\;$]
\noindent
\smallskip

\noindent
(a) Introduce \emph{the extended Bernoulli numbers} $C_{mn}$
 in terms of the classical ones:
$$C_{1n}=B_n,\quad
  C_{m+1,n}=\dfrac{n}{n+1} C_{m,n+1}
   -\dfrac{1}{n+1}\sum\limits_{k=1}^n
    \binom{n+1}{k} B_{k} C_{m,n-k+1}
 \mbox{ for }m,n\geq 1. \leqno{(2.4a)}$$
The numbers $C_{mn}$ are calculated in Table A.2 of Appendix
 for $m+n\leq 12$.
\medskip

\noindent
(b) Let us introduce \emph{the generating function}
$$C(\la,\mu)=\sum\limits_{m=1}^{\infty}\sum\limits_{n=1}^{\infty}
  \dfrac{C_{mn}}{m!n!}\la^{n-1}\mu^{m-1}. \eqno{(2.4b)\quad\bt}$$
\end{definition}

The formula~(2.4a) does not look very naturally.
But there is a more natural definition of $C_{mn}$ equivalent
 to~(2.4a), see Proposition~2.8.

\begin{example}
The first few values of the extended Bernoulli numbers are
$$C_{21}=-\frac{1}{6},\quad
  C_{22}=\frac{1}{6},\quad
  C_{23}=-\frac{1}{15},\quad
  C_{31}=0,\quad
  C_{32}=\frac{1}{15},\quad
  C_{41}=\frac{1}{30}.$$
Then the generating function $C(\la,\mu)$ starts with
$$C(\la,\mu)=-\frac{1}{2}+\frac{1}{12}(\la-\mu)+\frac{1}{24}\la\mu
 +\frac{1}{720}(\mu^3+4\la\mu^2-4\la^2\mu-\la^3)+\cdots$$
Up to degree 10 the function $C(\la,\mu)$ is computed in
 Example~A.3 of Appendix.
\ee
\end{example}

\noindent
If $W$ is a word in $P,Q$, then the expression $[W]$ in a long
 commutator is regarded as a formal symbol, i.e.
 $\bigl[PQ[W]\bigr]:=\bigl[P,[Q,[W]]\bigr]$.
But the symbol $W$ in a long commutator is considered
 as the word in $P,Q$, for $W=Q^nP^m$ one gets $[WQP]:=[Q^nP^mQP]$.

\begin{claim}
Let $P,Q$ be two elements of a Lie algebra $L$,
 $W$ be a word in the letters $P,Q$.
\smallskip

\noindent
(a) In the quotient $\bar L=L/\bigl[[L,L],[L,L]\bigr]$,
 for any word $W$ containing at least
 one letter $P$ and at least one letter $Q$,
 one has $[PQ[W]]=[QP[W]]$.
\smallskip

\noindent
(b) In the quotient $\bar L$, for any word $W$ containing
 exactly $m$ letters $P$ and exactly $n$ letters $Q$,
 one has $[WQP]=[Q^nP^mQP]$.
\end{claim}
\begin{proof}
(a) Since the element $[W]$ contains at least one commutator,
 then $\bigl[[P,Q],[W]\bigr]=0$ in the quotient $\bar L$.
The following calculations in $\bar L$ imply~(a):
 $\quad[PQ[W]]-[QP[W]]=$
$$=\bigl[P,Q[W]-[W]Q\bigr]-\bigl[Q,P[W]-[W]P\bigr]
  =\bigl(PQ[W]-P[W]Q-Q[W]P+[W]QP\bigr)-$$
$$\bigl(QP[W]-Q[W]P-P[W]Q+[W]PQ\bigr)
  =\bigl(PQ-QP\bigr)[W]-[W]\bigl(PQ-QP\bigr)
  =\bigl[[P,Q],[W]\bigr]=0.$$
\vspace{1pt}

\noindent
(b) By the item (a) one can permute the letters of $W$, i.e.
 one may assume $W=Q^nP^m$.
\end{proof}

Let $L$ be the Lie algebra freely generated by the symbols $P,Q$.
Recall that the series $H_1\in\hat L$ was introduced in Theorem~2.3.
Put $\bar L=L/\bigl[[L,L],[L,L]\bigr]$.
As usual by $\wh{\bar L}$ denote the algebra of formal series of
 elements from $\bar L$.

\begin{claim}
In the algebra $\wh{\bar L}$, for the derivative $D=H_1\pd{}{Q}$
 and all $m,n\geq 1$, one has
$$\begin{array}{ll}
(2.7a) &
 [H_1,P]=-\sum\limits_{k=1}^{\infty}\dfrac{B_{k}}{k!} [PQ^{k}P];\\ \\
(2.7b) &
 D[Q^nP]=(n-1)[Q^{n-2}PQP]
  -\sum\limits_{k=1}^{\infty}\dfrac{B_{k}}{k!}[Q^{n-1}PQ^{k}P];\\ \\
(2.7c) &
 D[Q^{n-1}P^{m-1}QP]=(n-1)[Q^{n-2}P^{m}QP]
  -\sum\limits_{k=1}^{\infty}\dfrac{B_{k}}{k!}[Q^{k+n-2}P^{m}QP].
\end{array}$$
\end{claim}
\begin{proof}
(a) It suffices to rewrite the formula of $H_1$ from Theorem~2.3
 as follows:
$$H_1\st{(2.3)}{=}P+\sum\limits_{k=1}^{\infty}\dfrac{B_k}{k!}[Q^kP]\quad
 \Ra\quad
 [H_1,P]=-\left[P,\sum_{k=1}^{\infty}\frac{B_{k}}{k!} [Q^{k}P]\right]
 =-\sum_{k=1}^{\infty}\frac{B_{k}}{k!} [PQ^{k}P].$$
Observe that here \emph{the first key point}
 ($B_n=0$ for odd $n\geq 3$) was used.
\smallskip

\noindent
(b) Induction on $n$.
The base $n=1$ follows from~(a):
 $D[QP]\st{(2.1c)}{=}[DQ,P]+[Q,DP]=$
 $=[H_1,P]\st{(2.7a)}{=}
  -\sum\limits_{k=1}^{\infty}\dfrac{B_{k}}{k!} [PQ^{k}P]$.
Induction step (from $n$ to $n+1$): $D[Q^{n+1}P]\st{(2.1c)}{=}$
$$\st{(2.1c)}{=}\bigl[H_1,[Q^nP]\bigr]+\bigl[Q,D[Q^nP]\bigr]
 =[PQ^nP]+(n-1)[Q^{n-1}PQP]
  -\left[Q,\sum_{k=1}^{\infty}\frac{B_{k}}{k!}[Q^{n-1}PQ^{k}P]\right].$$
It remains to apply Claim~2.6a:
 $[PQ^nP]+(n-1)[Q^{n-1}PQP]=n[Q^{n-1}PQP]$.
The item (b) and
 $D[Q^{n-1}P^{m-1}QP]\st{(2.6a)}{=}D[P^{m-1}Q^nP]
 \st{(2.1c)}{=}\underbrace{\bigl[P,[P,[\cdots}_{m-1},D[Q^nP]\cdots]\bigr]$
 imply (c).
\end{proof}

The following result gives a natural definition of the extended
 Bernoulli numbers $C_{mn}$: they give rise to an explicit
 compressed CBH formula (\emph{the second key point}).

\begin{proposition}[compressed variant of CBH]
Let $L$ be the Lie algebra freely generated by the symbols $P,Q$.
Under the natural projection
 $\hat L\to \wh{\bar L}$, where $\bar L=L/\bigl[[L,L],[L,L]\bigr]$,
 the Hausdorff series $H=\log\bigl(\exp(P)\cdot\exp(Q)\bigr)$ maps
 onto the series
$$\bar H=P+Q+\sum_{m=1}^{\infty}\sum_{n=1}^{\infty}
 \frac{C_{mn}}{m!n!} [Q^{n-1}P^{m-1}QP].\leqno{(2.8)}$$
\end{proposition}
\noindent
\emph{Proof.}
By Theorem~2.3 the series $H$ maps onto the series
 $\bar H=\sum\limits_{m=0}^{\infty}\bar H_m$, where
$$\bar H_0=Q,\quad
  \bar H_1=H_1=P+\sum_{n=1}^{\infty}\frac{B_n}{n!}[Q^nP],\quad
  \bar H_{m+1}=\frac{1}{m!}D^m(H_1)
 \mbox{ for }m\geq 1 \mbox{ and }D=H_1\pd{}{Q}.$$
It remains to prove the following formula:
$$D^m(H_1)=\sum\limits_{n=1}^{\infty}
  \dfrac{C_{m+1,n}}{n!}[Q^{n-1}P^{m}QP]\mbox{ for each }m\geq 1.
  \leqno{(2.8'_m)}$$
The base $m=1$ is completely analogous to the inductive step
 (from $m-1$ to $m$):
$$D^m(H_1)\st{(2.8'_{m-1})}{=}
  D\left( \sum\limits_{n=1}^{\infty}
  \dfrac{C_{mn}}{n!}[Q^{n-1}P^{m-1}QP] \right)
  =\sum\limits_{n=1}^{\infty}
  \dfrac{C_{mn}}{n!} D[Q^{n-1}P^{m-1}QP] \st{(2.7c)}{=}$$
$$ \st{(2.7c)}{=}\sum\limits_{n=1}^{\infty}
  \dfrac{C_{mn}}{n!} \left(
 (n-1)[Q^{n-2}P^{m}QP]
  -\sum_{k=1}^{\infty}\frac{B_{k}}{k!}[Q^{k+n-2}P^{m}QP] \right)=$$
$$=\sum\limits_{n=0}^{\infty}
  \dfrac{C_{m,n+1}}{(n+1)!} n [Q^{n-1}P^{m}QP]
 -\sum\limits_{n=1}^{\infty}\sum\limits_{k=1}^{\infty}
  \dfrac{C_{mn}}{n!} \frac{B_k}{k!} [Q^{k+n-2}P^{m}QP]=$$
$$=\sum\limits_{n=1}^{\infty} \left(
  \dfrac{C_{m,n+1}}{(n+1)!}n-\sum_{k=1}^n
  \frac{C_{m,n-k+1}}{(n-k+1)!} \frac{B_{k}}{k!} \right)
   [Q^{n-1}P^{m}QP]   \st{(2.4a)}{=}
  \sum\limits_{n=1}^{\infty}
  \dfrac{C_{m+1,n}}{n!}[Q^{n-1}P^{m}QP].\eqno{\square}$$
\medskip

The series $\bar H$ will be calculated up to degree 10
 in Proposition~A.4 of Appendix.


\subsection{Properties of the extended Bernoulli numbers
 and $C(\la,\mu)$}

\begin{claim}
Under the natural projection
 $\hat L\to \wh{\bar L}$, where $\bar L=L/\bigl[[L,L],[L,L]\bigr]$,
 the Hausdorff series $H=\log\bigl(\exp(P)\cdot\exp(Q)\bigr)$ maps
 onto the series\footnote{
  Note that $B_1=-\dfrac{1}{2}$, but one sets $C'_{11}=\dfrac{1}{2}$. }
$$\bar H=P+Q+\sum_{m=1}^{\infty}\sum_{n=1}^{\infty}
  \frac{C'_{mn}}{m!n!} [P^{n-1}Q^{m-1}PQ],\mbox{ where }
 C'_{11}=\frac{1}{2},\quad
 C'_{1n}=B_n, \leqno{(2.9a)}$$
$$C'_{m+1,n}=\dfrac{n}{n+1} C'_{m,n+1}
   -\dfrac{1}{n+1}\sum\limits_{k=1}^n
    \binom{n+1}{k} C'_{1k} C'_{m,n-k+1}
 \mbox{ for }m\geq 1,\; n\geq 2. \leqno{(2.9b)}$$
\end{claim}
\noindent
\emph{Proof} is completely analogous to the proof of
 Proposition~2.8.
One can use the following analog of Theorem~2.3
 \cite[the remark after Corollary~3.25 on p.~80]{Reu}:\\
 the Hausdorff series $H=\log\bigl(\exp(P)\cdot\exp(Q)\bigr)$
 is equal to
 $H=\sum\limits_{m=0}^{\infty} H'_m$, where
$$H'_0=P,\quad
  H'_1=Q+\dfrac{1}{2}[P,Q]
  +\sum\limits_{n=1}^{\infty}\dfrac{B_{2n}}{(2n)!}[P^{2n}Q],\quad
  H'_m=\dfrac{1}{m}\left(H'_1\pd{}{P}\right)(H'_{m-1})\mbox{ for }m\geq 2.$$
Here the derivative $D'=H'_1\pd{}{P}$ maps $P$ onto $H'_1$ and
 $Q$ onto 0.
To get the formula (2.9a) it suffices
 to apply the equations similar to Claim~2.7, where
 $P,Q$ are interchanged.
$$\begin{array}{ll}
(2.7a') & [H'_1,Q]
 =-\dfrac{1}{2}[QPQ]-\sum\limits_{k=2}^{\infty}\dfrac{B_{k}}{k!} [QP^{k}Q]
 =-\sum\limits_{k=1}^{\infty}\dfrac{C'_{1k}}{k!} [QP^{k}Q];\\ \\
(2.7b') &
 D'[P^nQ]=(n-1)[P^{n-2}QPQ]
  -\sum\limits_{k=1}^{\infty}\dfrac{C'_{1k}}{k!}[P^{n-1}QP^{k}Q];\\ \\
(2.7c') &
 D'[P^{n-1}Q^{m-1}PQ]=(n-1)[P^{n-2}Q^{m}PQ]
  -\sum\limits_{k=1}^{\infty}\dfrac{C'_{1k}}{k!}[P^{k+n-2}Q^{m}PQ].
\end{array}$$
Actually, it remains to deduce by $(2.7c')$ the formula analogous to
 $(2.8'_m)$:
$$(D')^m(H'_1)=\sum\limits_{n=1}^{\infty}
  \dfrac{C'_{m+1,n}}{n!}[P^{n-1}Q^{m}PQ]\mbox{ for every }m\geq 1.
  \eqno{\square}$$

\begin{lemma}
The extended Bernoulli numbers are symmetric
 in the following sense:
$$C_{mn}=(-1)^{m+n}C_{nm} \mbox{ for all }m,n\geq 1.\leqno{(2.10)}$$
Hence the generating function $C(\la,\mu)$ obeys the symmetry
 $C(\la,\mu)=C(-\mu,-\la)$.
\end{lemma}
\begin{proof}
Let us rewrite the recursive formula (2.9b) in a more explicit form:
$$C'_{m+1,n}=\dfrac{n}{n+1} C'_{m,n+1}-\dfrac{1}{2}C'_{mn}
   -\dfrac{1}{n+1}\sum\limits_{k=1}^{[\frac{n}{2}]}
    \binom{n+1}{2k} B_{2k} C'_{m,n-2k+1}.$$
In the same form the formula (2.4a) looks like
$$C_{m+1,n}=\dfrac{n}{n+1} C_{m,n+1}+\dfrac{1}{2}C_{mn}
   -\dfrac{1}{n+1}\sum\limits_{k=1}^{[\frac{n}{2}]}
    \binom{n+1}{2k} B_{2k} C_{m,n-2k+1}.$$
If the latter equation is multiplied by $(-1)^{m+n}$,
 then one obtains
$$(-1)^{m+n}C_{m+1,n}=\dfrac{n}{n+1}(-1)^{m+n} C_{m,n+1}
   -\dfrac{1}{2}(-1)^{m+n-1}C_{mn}-$$
$$ -\dfrac{1}{n+1}\sum\limits_{k=1}^{[\frac{n}{2}]}
    \binom{n+1}{2k} B_{2k} (-1)^{m+n-2k} C_{m,n-2k+1}.$$
Hence the numbers $C'_{mn}$ and $(-1)^{m+n-1}C_{mn}$ obey
 the same recursive relation.
Since $C'_{1n}=(-1)^nB_n=(-1)^nC_{1n}$ (\emph{the first key point}) for each
 $n\geq 1$, then $C'_{mn}=(-1)^{m+n-1}C_{mn}$ for all $m,n\geq 1$.
Then the formula (2.9a) converts to
$$\bar H-P-Q=\sum_{m,n\geq 1}
  (-1)^{m+n-1}\frac{C_{mn}}{m!n!} [P^{n-1}Q^{m-1}PQ]
 \st{(2.6)}{=}\sum_{m,n\geq 1}
  (-1)^{m+n}\frac{C_{mn}}{m!n!} [Q^{m-1}P^{n-1}QP].$$
By comparing the above formula with (2.8), one gets
 $C_{mn}=(-1)^{m+n}C_{nm}$ as required.
\end{proof}

The following two assertions show that the extended Bernoulli
 numbers $C_{mn}$ are not too complicated,
 contrary they can be expressed via binoms and classical
 Bernoulli numbers.

\begin{lemma}
The extended Bernoulli numbers can be expressed via
 the classical ones:
$$C_{mn}=\sum\limits_{k=0}^{m-1}\binom{m}{k}B_{n+k}
 \mbox{ for each }m\geq 1.\qquad
 \mbox{ In particular, one gets}\leqno{(2.11)}$$
$$C_{1n}=B_n,\;
  C_{2n}=B_n+2B_{n+1},\;
  C_{3n}=B_n+3B_{n+1}+3B_{n+2},\;
  C_{4n}=B_n+4B_{n+1}+6B_{n+2}+4B_{n+3}.$$
\end{lemma}
\noindent
\emph{Proof.}
Multiplying (2.4a) by $\frac{n+1}{n}$, one has
 $\frac{n+1}{n}C_{m+1,n}=C_{m,n+1}
   -\frac{1}{n}\sum\limits_{k=1}^n
    \binom{n+1}{k} B_k C_{m,n-k+1}$,~i.e.
$$C_{m,n+1}=\dfrac{n+1}{n}C_{m+1,n}
   +\dfrac{1}{n}\sum\limits_{k=1}^n
    \binom{n+1}{k} B_k C_{m,n-k+1}
    \mbox{ for all } m,n\geq 1.\leqno{(2.4a')}$$
The equation (2.11) will be checked by induction on $m$.
The base $m=1$ follows from Definition~2.4a.
Suppose that the formula (2.11) holds for $m$,
 let us prove it for $m+1$.
$$C_{m+1,n}\st{(2.10)}{=}(-1)^{m+n+1}C_{n,m+1}
\st{(2.4a')}{=}\dfrac{(-1)^{m+n+1}}{m}\left(
 (m+1)C_{n+1,m} +\sum\limits_{k=1}^m
  \binom{m+1}{k} B_k C_{n,m-k+1}\right)$$
$$\st{(2.10)}{=}(-1)^{m+n+1}\left(
 \dfrac{m+1}{m} (-1)^{m+n+1} C_{m,n+1} +\dfrac{1}{m}\sum\limits_{k=1}^m
  \binom{m+1}{k} B_k (-1)^{m+n-k+1} C_{m-k+1,n}\right)=$$
$$\mbox{(by hypothesis)}=
 \dfrac{m+1}{m} \sum_{k=0}^{m-1}\binom{m}{k} B_{n+k+1}
 +\dfrac{1}{m}\sum\limits_{k=1}^m
  (-1)^k \binom{m+1}{k} B_k \sum_{l=0}^{m-k} \binom{m-k+1}{l}B_{n+l}$$
$$=\dfrac{m+1}{m} \sum_{l=1}^{m}\binom{m}{l-1} B_{n+l}
 +\dfrac{1}{m}\sum\limits_{k=1}^m (-1)^k \binom{m+1}{k} B_k B_n+$$
$$+\dfrac{1}{m}\sum\limits_{k=1}^m (-1)^k \binom{m+1}{k} B_k
  \sum_{l=1}^{m-k} \binom{m-k+1}{l}B_{n+l}=$$
$$\frac{B_n}{m}\sum_{k=1}^m(-1)^k\binom{m+1}{k}B_k
 +\sum_{l=1}^{m} \frac{B_{n+l}}{m}\biggl( (m+1)\binom{m}{l-1}
 +\sum_{k=1}^{m-l}(-1)^k\binom{m+1}{k}B_k\binom{m-k+1}{l}\biggr).$$
Since by Lemma~2.2c the first term is equal to $B_n$, then
 it remains to check the formula
$$\frac{m+1}{m}\binom{m}{l-1}+\frac{1}{m}\sum_{k=1}^{m-l}
 (-1)^k\binom{m+1}{k}B_k\binom{m-k+1}{l}=\binom{m+1}{l}
 \mbox{ for all }m\geq l\geq 1.$$
The left hand side of the above equation is
$$\frac{(m+1)\cdot(m-1)!}{(l-1)!(m-l+1)!}
 +\sum_{k=1}^{m-l}\left(
  (-1)^kB_k\cdot\frac{(m+1)\cdot(m-1)!}{k!(m-k+1)!}\cdot
  \frac{(m-k+1)!}{l!(m-k-l+1)!} \right)=$$
$$=\frac{(m+1)\cdot(m-1)!}{(l-1)!(m-l+1)!}
 +\frac{(m+1)\cdot(m-1)!}{l!(m-l+1)!}\cdot\sum_{k=1}^{m-l}
 (-1)^kB_k\binom{m-l+1}{k}\st{(2.2c)}{=}$$
$$\st{(2.2c)}{=}
 \frac{(m+1)\cdot(m-1)!}{(l-1)!(m-l+1)!}\left(1+\frac{m-l}{l}\right)
 =\frac{(m+1)!}{l!(m-l+1)!}=\binom{m+1}{l}
\mbox{ as required}.\eqno{\qed}$$
\smallskip

Now one can get an explicit expression of $C(\la,\mu)$,
 which is not followed immediately from Definition~2.4 or Proposition~2.8
 (\emph{the third key point}).

\begin{proposition}
The generating function $C(\la,\mu)$ from Definition~2.4b is equal to
$$C(\la,\mu)=\dfrac{e^{\mu}-1}{\la\mu}\cdot\left(
 \dfrac{\la+\mu}{e^{\la+\mu}-1}-\dfrac{\mu}{e^{\mu}-1}
 \right).\leqno{(2.12)}$$
\end{proposition}
\noindent
\emph{Proof.}
The next trick is applied: to get a non-trivial relation,
 one destroys a symmetry.
$$\la\mu C(\la,\mu)
=\sum_{m,n\geq 1}\frac{C_{mn}}{m!n!}\la^n\mu^m
\st{(2.11)}{=}\sum_{n=1}^{\infty}\sum_{m=1}^{\infty}\left(
 \sum_{k=0}^{m-1}\binom{m}{k}B_{n+k}\right)
 \cdot \frac{\la^n}{n!} \cdot \frac{\mu^m}{m!}=$$
$$\sum_{k=0}^{\infty} \sum_{n=1}^{\infty} \sum_{m=k+1}^{\infty}
 \frac{B_{n+k}\la^n\mu^m}{n!(m-k)!k!}
=\sum_{k=0}^{\infty} \sum_{n=1}^{\infty} \sum_{m=1}^{\infty}
 B_{n+k} \frac{\la^n\mu^m\mu^k}{n!m!k!}
=\left( \sum_{m=1}^{\infty}\frac{\mu^m}{m!}\right)
 \left(\sum_{k=0}^{\infty} \sum_{n=1}^{\infty} B_{n+k}
 \frac{\la^n\mu^k}{n!k!}\right)$$
$$=(e^{\mu}-1)\cdot \left( \sum_{n=1}^{\infty}B_n\frac{\la^n}{n!}
 +\sum_{k,n\geq 1} B_{n+k} \frac{\la^n\mu^k}{n!k!} \right)
\st{(2.1b)}{=}(e^{\mu}-1)\cdot \left(
 \frac{\la}{e^{\la}-1}-1+\ti C(\la,\mu)\right),\mbox{ where}$$
$\ti C(\la,\mu)=\sum\limits_{k,n\geq 1}B_{n+k}\dfrac{\la^n\mu^k}{n!k!}$.
Since $(-1)^{n+k}B_{n+k}=B_{n+k}$ for $n,k\geq 1$
 (\emph{the first key point}),
$$\mbox{then }\ti C(-\mu,-\la)
 =\sum\limits_{k,n\geq 1} (-1)^kB_{n+k}\dfrac{\la^k\mu^n}{k!n!}
 =\sum\limits_{k,n\geq 1} B_{n+k}\dfrac{\la^k\mu^n}{k!n!}
 =\ti C(\la,\mu).\qquad \mbox{ One obtains}$$
$$(e^{\mu}-1)\cdot\left(
 \frac{\la}{e^{\la}-1}-1+\ti C(\la,\mu)\right)=\la\mu
 C(\la,\mu)\st{(2.10)}{=}\la\mu C(-\mu,-\la)=$$
$$=(e^{-\la}-1)\cdot\left(
 \frac{-\mu}{e^{-\mu}-1}-1+\ti C(-\mu,-\la)\right)
=(e^{-\la}-1)\cdot\left(
 \frac{-\mu}{e^{-\mu}-1}-1+\ti C(\la,\mu)\right),\mbox{ hence}$$
$$(e^{-\la}-e^{\mu})\ti C(\la,\mu)
=(e^{\mu}-1)\Bigl(\frac{\la}{e^{\la}-1}-1\Bigr)-
 (e^{-\la}-1)\Bigl(\frac{-\mu}{e^{-\mu}-1}-1\Bigr)
=\la\dfrac{e^{\mu}-1}{e^{\la}-1}+\mu\dfrac{e^{-\la}-1}{e^{-\mu}-1}
 +e^{-\la}-e^{\mu}.$$
By substituting the above formula into the expression of
 $\la\mu C(\la,\mu)$ via $\ti C(\la,\mu)$, one gets
$$C(\la,\mu)=\dfrac{e^{\mu}-1}{\la\mu}\cdot\left(
 \dfrac{\la}{e^{\la}-1} +\left(
 \la\dfrac{e^{\mu}-1}{e^{\la}-1}+\mu\dfrac{e^{-\la}-1}{e^{-\mu}-1}\right)
 :(e^{-\la}-e^{\mu}) \right).$$
Note that after dividing by $\la+\mu$ both series
 $\la\dfrac{e^{\mu}-1}{e^{\la}-1}+\mu\dfrac{e^{-\la}-1}{e^{-\mu}-1}$
 and $e^{-\la}-e^{\mu}$ begin with 1.
It remains to prove the following equation:
$$\la\dfrac{e^{\mu}-1}{e^{\la}-1}
 +\mu\dfrac{e^{-\la}-1}{e^{-\mu}-1}
 =(e^{-\la}-e^{\mu})\cdot \left( \dfrac{\la+\mu}{e^{\la+\mu}-1}
 -\dfrac{\la}{e^{\la}-1}-\dfrac{\mu}{e^{\mu}-1} \right),\mbox{ or}$$
$$\la e^{\la}\dfrac{e^{\mu}-1}{e^{\la}-1}
 +\mu e^{\mu}\dfrac{1-e^{\la}}{1-e^{\mu}}
 =(1-e^{\la+\mu})\cdot \left( \dfrac{\la+\mu}{e^{\la+\mu}-1}
 -\dfrac{\la}{e^{\la}-1}-\dfrac{\mu}{e^{\mu}-1} \right),\mbox{ or}$$
$$(e^{\mu}-1)\dfrac{\la e^{\la}}{e^{\la}-1}
 +(e^{\la}-1)\dfrac{\mu e^{\mu}}{e^{\mu}-1}
 =-\la-\mu-(1-e^{\la+\mu})\cdot \left( \dfrac{\la}{e^{\la}-1}
 +\dfrac{\mu}{e^{\mu}-1} \right),\mbox{ or}$$
$$(e^{\mu}-1)\cdot\left( \dfrac{\la}{e^{\la}-1}+\la \right)
 +(e^{\la}-1)\cdot\left( \dfrac{\mu}{e^{\mu}-1}+\mu \right)
 =-\la-\mu+(e^{\la+\mu}-1)\cdot \left( \dfrac{\la}{e^{\la}-1}
 +\dfrac{\mu}{e^{\mu}-1} \right),$$
$$\mbox{or  }\qquad\la e^{\mu}+\mu e^{\la}
 =(e^{\la+\mu}-e^{\mu})\dfrac{\la}{e^{\la}-1}
 +(e^{\la+\mu}-e^{\la})\dfrac{\mu}{e^{\mu}-1},\mbox{ that is clear}.
 \eqno{\square}$$


\section{Compressed hexagon equation}

In this section one shall write explicitly the compressed hexagon
 and prove Theorem~1.5a.


\subsection{First simplification of the hexagon}

Lemma~3.1 gives a first simplification of the hexagon~(1.3b).
Due to this \emph{4th key point} it remains to apply CBH formula exactly once.

\begin{lemma}
In $\hat A_3$ the hexagon equation~(1.3b) is equivalent
 to the following equation:
$$\exp(a+b+c)=\exp(\psi(c,a))\cdot\exp(\psi(b,c))\cdot\exp(\psi(a,b)),
 \leqno{(3.1)}$$
 where $\psi(a,b):=\log\bigl(\exp(\ph(a,b))\cdot\exp(a)\bigr)\in\hat L_3$,
 $a:=t^{12}$, $b:=t^{23}$, $c:=t^{13}$.
\end{lemma}
\begin{proof}
Let us rewrite the original hexagon (1.3b) in a more explicit form:
$$\exp(t^{13}+t^{23})=
  \exp\bigl(\ph(t^{13},t^{12})\bigr)\cdot \exp(t^{13})\cdot
  \exp\bigl(-\ph(t^{13},t^{23})\bigr)\cdot \exp(t^{23})\cdot
  \exp\bigl(\ph(t^{12},t^{23})\bigr).$$
Introduce the symbols $a=t^{12}$, $b=t^{23}$, $c=t^{13}$
 and apply the symmetry $\ph(b,c)=-\ph(c,b)$.
$$\exp(b+c)=\exp\bigl(\ph(c,a)\bigr)\cdot \exp(c)\cdot
 \exp\bigl(\ph(b,c)\bigr)\cdot \exp(b)\cdot \exp\bigl(\ph(a,b)\bigr).$$
It remains to multiply the above equation by $\exp(a)$ from the right.
Note the element $a$ commutes with $b+c$ in $L_3$ by definition.
Moreover, $a+b+c$ is a central element of $L_3$.
\end{proof}

Recall $\wh{\bar L}_3$ is the algebra of formal series of
 elements from $\bar L_3=L_3/\bigl[[L_3,L_3],[L_3,L_3]\bigr]$.

\begin{claim}
In $\wh{\bar L}_3$, for the compressed series
 $\bar\psi(a,b)=\log\bigl(\exp(\bar\ph(a,b))\cdot\exp(a)\bigr)$, one has
$$\bar\psi(a,b)=a+\sum_{n=0}^{\infty}\frac{B_n}{n!}[a^n\bar\ph(a,b)]
 =a+\bar\ph(a,b)-\frac{1}{2}[a,\bar\ph(a,b)]
 +\sum_{n=1}^{\infty}\frac{B_{2n}}{(2n)!}[a^{2n}\bar\ph(a,b)].
 \leqno{(3.2)}$$
\end{claim}
\begin{proof}
Apply Proposition~2.8 for the Hausdorff series
 $\bar\psi(a,b)=\log\bigl(\exp(\bar\ph(a,b))\cdot\exp(a)\bigr)$
 and $Q=a$, $P=\bar\ph(a,b)$.
Since all commutators including in $\bar\ph(a,b)$ commute with
 each other in $\bar L_3$, then $[Q^nP]=[a^n\bar\ph(a,b)]$ and
 $[Q^{n-1}P^{m-1}QP]=0$ for all $m\geq 2$, $n\geq 1$.
\end{proof}

In particular, by Claim~3.2 the series $\bar\psi(a,b)$
 starts with
$$\bar\psi(a,b)=a+\bar\ph(a,b)-\frac{1}{2}[a,\bar\ph(a,b)]
 +\frac{1}{12}[a,[a,\bar\ph(a,b)]]
 -\frac{1}{720}[a,[a,[a,[a,\ph(a,b)]]]]
 +\cdots$$


\subsection{The Lie algebra $L_3(\la,\mu)$}

Recall that $\bar L_3=L_3/\bigl[[L_3,L_3],[L_3,L_3]\bigr]$ by definition.

\begin{definition}[Lie algebra $L_3(\la,\mu)$]
Let us introduce the linear space $L_3(\la,\mu)$ generated by
 $a,b$, and $\la^k\mu^l[ab]$ for $k,l\geq 0$,
 where $\la,\mu\in\C$ are formal parameters.
At this moment the expression $[ab]$ means a formal symbol,
 which denotes an element of $L_3(\la,\mu)$.
Define in $L_3(\la,\mu)$ \emph{the bracket} by
 $[a,b]:=[ab]$, $[a,\la^k\mu^l[ab]]:=\la^{k+1}\mu^l[ab]$,
 $[b,\la^k\mu^l[ab]]:=\la^{k}\mu^{l+1}[ab]$,
 the other brackets are zero.
One can check that this bracket satisfies the Jacobi identity.
Actually, only the following identities
 (up to permutation $a\lra b$) contain non-zero terms:
$$\bigl[ a, [ a, \la^k\mu^l[ab] ] \bigr]
 +\bigl[ a, [ \la^k\mu^l[ab], a ] \bigr] =0,\;
  \bigl[ a, [ b, \la^k\mu^l[ab] ] \bigr]
 +\bigl[ b, [ \la^k\mu^l[ab], a ] \bigr]
 +\bigl[ \la^k\mu^l[ab],   [ ab ] \bigr] =0.$$
Hence the space $L_3(\la,\mu)$ becomes a Lie algebra.
Now the formal symbol $[ab]$ is
 the true commutator $[a,b]\in L_3(\la,\mu)$.
Note that in $L_3(\la,\mu)$ one has
 $[a^{k}b^{l}ab]=\la^{k}\mu^{l}[a,b]$.

By $\hat L_3(\la,\mu)$ denote the algebra of formal series
 of elements from $L_3(\la,\mu)$.
Observe that in $\hat L_3(\la,\mu)$ any series of commutators
 $[a^kb^lab]$ is (a series in $\la,\mu$)$\times[a,b]$.
\ed
\end{definition}

\noindent
Remind the notation
 $[a_1a_2\ldots a_k]:=[a_1,[a_2,[\ldots,a_k]\ldots]]$
 from Definition~1.2c.

\begin{proposition}
(a) Let $L(a+b+c)\subset L_3$ be the 1-dimensional Lie subalgebra
 generated by the element $a+b+c$,
 $L(a,b)\subset L_3$ be the Lie subalgebra freely generated by $a,b$.
Then the Lie algebra $L_3$ is isomorphic to the direct sum
 $L(a+b+c)\oplus L(a,b)$.
\smallskip

\noindent
(b) The quotient
 $\bar L(a,b):=L(a,b)/\bigl[[L(a,b),L(a,b)],[L(a,b),L(a,b)]\bigr]$
 is isomorphic to the Lie algebra $L_3(\la,\mu)$
 introduced in Definition~3.3.
\smallskip

\noindent
(c) The quotient $\bar L_3=L_3/\bigl[[L_3,L_3],[L_3,L_3,]\bigr]$
 is isomorphic to $L(a+b+c)\oplus L_3(\la,\mu)$.
\end{proposition}
\begin{proof}
(a) In the setting $a:=t^{12}$, $b:=t^{23}$, $c=t^{13}$,
 the element $a+b+c$ is a central element of $L_3$.
Moreover, the defining relations of the Lie algebra $L_3$ become
 $[a,a+b+c]=[b,a+b+c]=0$.
Now the isomorphism $L_3\cong L(a+b+c)\oplus L(a,b)$ is obvious.
\smallskip

\noindent
(b) By definition the Lie algebra $L(a,b)$ is linearly generated by
 $a,b$ and all commutators $[wab]$, where $w$ is a word in $a,b$.
By Claim~2.6b the quotient $\bar L(a,b)$ is linearly generated
 by $a,b$ and commutators $[a^kb^lab]$ with $k,l\geq 0$.
The only non-zero brackets of these elements in $\bar L(a,b)$ are
 $[a,b]$, $[a,a^kb^lab]=[a^{k+1}b^lab]$, $[b,a^kb^lab]=[a^{k}b^{l+1}ab]$.
One has
 $\bigl[[L_3(\la,\mu),L_3(\la,\mu)],[L_3(\la,\mu),L_3(\la,\mu)]\bigr]=0$.
The required isomorphism $\bar L(a,b)\to L_3(\la,\mu)$
 is defined by $[a^kb^lab]\mapsto\la^k\mu^l[a,b]$.
The item (c) follows from (a) and (b).
\end{proof}

Due to Proposition~3.4 the image of the logarithm of a Drinfeld associator
 is $[a,b]\times$ (a series in the commuting parameters $\la,\mu$).
Hence the compressed hexagon lives in $\hat L_3(\la,\mu)$.
This simplification (\emph{the 5th key point}) allows
 to solve the compressed hexagon completely.

\begin{definition}[generating series $f(\la,\mu)$ and $g(\la,\mu)$;
 symmetric coefficients $\al_{kl}$]
\noindent
\smallskip

\noindent
Let $\C[[\la,\mu]]$ be the set of formal power series
 with complex coefficients in the commuting arguments $\la,\mu$.
Introduce series $f,g\in \C[[\la,\mu]]$ by the formulae
 in the algebra $\hat L_3(\la,\mu)$
 $$\bar\ph(a,b)=\sum\limits_{k,l\geq 0}\al_{kl}[a^kb^lab]
   =f(\la,\mu)\cdot [a,b],\quad
   \bar\psi(a,b)=a+g(\la,\mu)\cdot [a,b],\mbox{ respectively}.$$
The formula for $\bar\ph(a,b)$ is a general form
 of a series in $\hat L_3(\la,\mu)$.
A series $g(\la,\mu)$ exists due to Claim~3.2.
Since $\bar\ph(a,b)=-\bar\ph(b,a)$, then
 the series $f(\la,\mu)=\sum\limits_{k,l\geq 0}\al_{kl}\la^k\mu^l$
 (but not $g(\la,\mu)$) is symmetric in the arguments $\la,\mu$,
 i.e. $\al_{kl}=\al_{lk}$ are symmetric coefficients.
\ed
\end{definition}

\begin{claim}
(a) The series $f(\la,\mu)$ and $g(\la,\mu)$ are related as
 follows
$$g(\la,\mu)=\sum_{k=0}^{\infty}\frac{B_k}{k!}\la^k\cdot f(\la,\mu)
=\left(1-\frac{\la}{2}+\frac{\la^2}{12}-\frac{\la^4}{720}+\cdots\right)
 f(\la,\mu)\st{(2.1b)}{=}\frac{\la f(\la,\mu)}{e^{\la}-1}. \leqno{(3.6a)}$$

\noindent
(b) In the algebra $\hat L_3(\la,\mu)$ one has
$$\bar\psi(b,c)=b+g(\mu,\rho)\cdot[a,b],\quad
  \bar\psi(c,a)=c+g(\rho,\la)\cdot[a,b],\mbox{ where }
  \rho=-\la-\mu.\leqno{(3.6b)}$$
\end{claim}
\noindent
\emph{Proof}.
The following calculations imply the item~(a):
 $g(\la,\mu)[ab]
 \st{(3.5)}{=} \bar\psi(a,b)-a\st{(3.2)}{=}$
$$\st{(3.2)}{=}
   \sum\limits_{k=0}\dfrac{B_k}{k!}[a^k\bar\ph(a,b)]
  \st{(3.5)}{=}
   \sum\limits_{k=0}\dfrac{B_k}{k!}\bigl[a^k f(\la,\mu)[ab]\bigr]
  \st{(3.3)}{=}
   \sum\limits_{k=0}\dfrac{B_k}{k!}\la^k f(\la,\mu) \cdot [ab].$$

\noindent
The item (b) follows from Definition~3.5 and
 $\bar\ph(b,c)=f(\mu,\rho)[ab]$, $\bar\ph(c,a)=f(\rho,\la)[ab]$.
The latter formula is proved similarly to the former one
 by using (1.2a) $[ab]=[bc]=[ca]$:
 $\bar\ph(b,c)\st{(3.5)}{=}$
$$\sum\limits_{k,l\geq 0}\al_{kl}[b^kc^lbc]
 \st{(1.2a)}{=}\sum\limits_{k,l\geq 0}\al_{kl}[b^k(-a-b)^lab]
 \st{(3.3)}{=} \sum\limits_{k,l\geq 0}\al_{kl}\mu^k(-\la-\mu)^l[ab]
 \st{(3.5)}{=} f(\mu,\rho)[ab].\eqno{\qed}$$

\begin{example}
For the series $\bar\ph^B(a,b)=f^B(\la,\mu)\cdot[a,b]$
 from Example~(1.4), one has
 $$\al_{00}=\dfrac{1}{6};\quad
   \al_{10}=\al_{01}=0;\quad
   \al_{20}=\al_{02}=-\dfrac{1}{90},\;
   \al_{11}=-\dfrac{1}{360};\quad
   \al_{30}=\al_{21}=\al_{12}=\al_{03}=0.$$
\noindent
Then up to degree 3 the corresponding series
 $f^B(\la,\mu)$ and $g^B(\la,\mu)$ are
$$f^B(\la,\mu)=\dfrac{1}{6}-\dfrac{4\la^2+\la\mu+4\mu^2}{360},\;
  g^B(\la,\mu)=\dfrac{1}{6}-\dfrac{\la}{12}
   +\dfrac{\la^2-\la\mu-4\mu^2}{360}
   +\dfrac{4\la^3+\la^2\mu+4\la\mu^2}{720}.$$
Moreover, one can compute the series from Claim~3.6b
 (they will be needed later):
$$g^B(\mu,\rho)=\dfrac{1}{6}-\dfrac{\mu}{12}
   -\dfrac{4\la^2+7\la\mu+2\mu^2}{360}
   +\dfrac{4\la^2\mu+7\la\mu^2+7\mu^3}{720}+\cdots,$$
$$g^B(\rho,\la)=\dfrac{1}{6}+\dfrac{\la+\mu}{12}
   -\dfrac{2\la^2+3\la\mu+\mu^2}{360}
   -\dfrac{7\la^3+14\la^2\mu+11\la\mu^2+4\mu^3}{720}+\cdots,\mbox{ hence}$$
$$G^B(\la,\mu):=g^B(\la,\mu)+g^B(\mu,\rho)+g^B(\rho,\la)
  =\dfrac{1}{2}-\frac{\la^2+\la\mu+\mu^2}{72}
  +\frac{\mu^3-3\la^2\mu-\la^3}{240}+\cdots,$$
$$T^B(\la,\mu):=1+\la g^B(\mu,\rho)-\mu g^B(\la,\mu)
  =1+\frac{\la-\mu}{6}+\frac{4\mu^3-\la\mu^2-8\la^2\mu-4\la^3}{360}+\cdots
  \eqno{\bt}$$
\end{example}


\subsection{Explicit form of the compressed hexagon~$\ol{(1.3b)}$}

\noindent
\begin{claim}
For the generators $a=t^{12}$, $b=t^{23}$, $c=t^{13}$ of the Lie
 algebra $L_3$, set $P:=\bar\psi(b,c)$, $Q:=\bar\psi(a,b)$.
Then in the algebra $\hat L_3(\la,\mu)$ one has\footnote{
 The series $g(\la,\mu)$ was introduced in Definition~3.5.}
$$[Q,P]=T(\la,\mu)\cdot[a,b],\mbox{ where }
 T(\la,\mu)=1+\la g(\mu,-\la-\mu)-\mu g(\la,\mu).\leqno{(3.8)}$$
\end{claim}
\begin{proof}
By Definition~3.4a and Claim~3.6 one gets
$$[Q,P]=\bigl[a+g(\la,\mu)[ab],b+g(\mu,-\la-\mu)[ab]\bigr]
 =[ab]+g(\mu,-\la-\mu)[aab]-g(\la,\mu)[bab].$$
It remains to use the relations $[aab]=\la[ab]$ and
 $[bab]=\mu[ab]$ of $L_3(\la,\mu)$.
\end{proof}

\begin{proposition}
Let $\bar\ph=\sum\limits_{k,l\geq 0}\al_{kl}[a^kb^lab]$
 be a compressed Drinfeld associator, $\al_{kl}\in\C$,\\
 $f(\la,\mu)=\sum\limits_{k,l\geq 0}\al_{kl}\la^k\mu^l$ be the generating
 function of the coefficients $\al_{kl}$.
Then the compressed hexagon~$\ol{(1.3b)}$ is equivalent to
 the following equation in the algebra $\C[[\la,\mu]]$
$$G(\la,\mu)+C(\la,\mu)\cdot T(\la,\mu)=0,
 \mbox{ where}\leqno{(3.9)}$$
$$G(\la,\mu):=g(\la,\mu)+g(\mu,\rho)+g(\rho,\la),\quad
  T(\la,\mu):=1+\la g(\mu,\rho)-\mu g(\la,\mu),\quad
  \rho:=-\la-\mu.$$
\end{proposition}
\noindent
\emph{Proof}.
Let us apply Proposition~2.8 for the Hausdorff series
 $\bar H=\log\bigl(\exp(P)\cdot\exp(Q)\bigr)$, where
 $P=\bar\psi(b,c)=b+g(\mu,\rho)[a,b]$,
 $\;Q=\bar\psi(a,b)=a+g(\la,\mu)[a,b]$.
One has
$$\bar H=P+Q+\sum_{m=1}^{\infty}\sum_{n=1}^{\infty}
 \frac{C_{mn}}{m!n!} [Q^{n-1}P^{m-1}QP].$$
Since in the quotient $\bar L_3$ all commutators commute, then
 the formula~(3.8) implies
$$[Q^{n-1}P^{m-1}QP]=T(\la,\mu)\cdot[a^{n-1}b^{m-1}ab]
 =\la^{n-1}\mu^{m-1}T(\la,\mu)\cdot[a,b],\mbox{ hence}$$
$$\bar H=a+g(\la,\mu)[a,b]+b+g(\mu,\rho)[a,b]+
 +\sum_{m=1}^{\infty}\sum_{n=1}^{\infty}
 \frac{C_{mn}}{m!n!} \la^{n-1}\mu^{m-1} T(\la,\mu)\cdot[a,b]=$$
$$=\mbox{ (by Definition~2.4b) }=
 a+b+(g(\la,\mu)+g(\mu,\rho))\cdot[a,b]+C(\la,\mu)T(\la,\mu)\cdot[a,b].$$
On the other hand, by Lemma~3.1 one has
 $\bar H=\log\bigl(\exp(-\bar\psi(c,a))\cdot\exp(a+b+c)\bigr)$.
But $a+b+c$ is a central element of $L_3$, hence
 $\bar H=a+b+c-\bar\psi(c,a)$.

Let us take together both above expressions for $\bar H$
 and use Claim~3.6b for $\bar\psi(c,a)$
$$a+b+\bigl(g(\la,\mu)+g(\mu,\rho)\bigr)[ab]+C(\la,\mu)T(\la,\mu)[ab]
 =a+b+c-\bigl(c+g(\rho,\la)[ab]\bigr)\Lra(3.9).\quad\square$$

\begin{example}
By $G^B_k,C_k,T^B_k$ denote the degree $k$ parts of the
  functions $G^B(\la,\mu)$, $C(\la,\mu)$, $T^B(\la,\mu)$, respectively.
Due to Examples 2.5 and 3.7 one can calculate
$$G^B_0(\la,\mu)=\frac{1}{2},\quad
  G^B_1(\la,\mu)=0,\quad
  G^B_2(\la,\mu)=-\frac{\la^2+\la\mu+\mu^2}{72},\quad
  G^B_3(\la,\mu)=\frac{\mu^3-3\la^2\mu-\la^3}{240};$$
$$C_0(\la,\mu)=-\frac{1}{2},\quad
  C_1(\la,\mu)=\frac{\la-\mu}{12},\quad
  C_2(\la,\mu)=\frac{\la\mu}{24},\quad
  C_3(\la,\mu)=\frac{\mu^3+4\la\mu^2-4\la^2\mu-\la^3}{720};$$
$$T^B_0(\la,\mu)=1,\quad
  T^B_1(\la,\mu)=\frac{\la-\mu}{6},\quad
  T^B_2(\la,\mu)=0,\quad
  T^B_3(\la,\mu)=\frac{4\mu^3-\la\mu^2-8\la^2\mu-4\la^3}{360}.$$
Now one can check by hands the first four compressed hexagons:
$$G^B_0+C_0=0,\;
  G^B_1+C_0T^B_1+C_1=0,\;
  G^B_2+C_1T^B_1+C_2=0,\;
  G^B_3+C_0T^B_3+C_2T^B_1+C_3=0.\eqno{\bt}$$
\end{example}
\smallskip

\noindent
{\bf Proof of Theorem~1.5a.}
The expression $\bar\ph(a,b)=\sum\limits_{k,l\geq 0}\al_{kl}[a^kb^lab]$
 follows from Definition~3.5.
The condition~$\bar\ph(a,b)=-\bar\ph(b,a)$ is the symmetry
 $\al_{kl}=\al_{lk}$.
The non-degeneracy~$\bar\ph(a,0)=\bar\ph(0,b)=0$ holds trivially.
Let us rewrite (3.9) as follows:
$$g(\la,\mu)\cdot\bigl(1-\mu C(\la,\mu)\bigr)
 +g(\mu,\rho)\cdot\bigl(1+\la C(\la,\mu)\bigr)
 +g(\rho,\la)+C(\la,\mu)=0.\leqno{(3.9')}$$
Then, for each $n\geq 1$, the degree $n$ parts of
 the functions $g(\la,\mu)$, $g(\mu,\rho)$, and $g(\rho,\la)$
 can be expressed linearly via the degree $k<n$ parts of the same functions.
Due to Claim~3.6a the degree $n$ parts of these functions contain
 the coefficients $\al_{kl}$ only for $k+l=n$.
Hence the coefficients $\al_{kl}$ could be expressed linearly
 in terms of $\al_{ij}$ with $i+j<k+l$.
\qed


\section{Solving the compressed hexagon equation}

Here one solves completely the compressed hexagon~$\ol{(1.3b)}=(3.9)$
 from Proposition~3.9.


\subsection{Further simplifications of the compressed hexagon}

\begin{lemma}
Let $f(\la,\mu)=\sum\limits_{k,l\geq 0}\al_{kl}\la^k\mu^l$ be
 the generating function of the coefficients $\al_{kl}$ of
 a compressed associator $\bar\ph\in\wh{\bar L}_3$.
Then the compressed hexagon~$\ol{(1.3b)}=(3.9)$ is equivalent
 to the equation~(4.1a) and to the equation~(4.1b)
 in the algebra $\C[[\la,\mu]]$
$$\dfrac{\la f(\la,\mu)}{e^{\la}-1}\bigl(1-\mu C(\la,\mu)\bigr)
 +\dfrac{\mu f(\mu,\rho)}{e^{\mu}-1}\bigl(1+\la C(\la,\mu)\bigr)
 +\dfrac{\rho f(\rho,\la)}{e^{\rho}-1}
 +C(\la,\mu)=0, \mbox{ where}\leqno{(4.1a)}$$
$$\rho=-\la-\mu,\quad
 C(\la,\mu)=\dfrac{e^{\mu}-1}{\la\mu}\cdot\left(
 \dfrac{\la+\mu}{e^{\la+\mu}-1}-\dfrac{\mu}{e^{\mu}-1} \right);$$
$$f(\la,\mu)+e^{\mu}f(\mu,-\la-\mu)+e^{-\la}f(\la,-\la-\mu)
=\dfrac{1}{\la+\mu}\left(
 \dfrac{e^{\mu}-1}{\mu}+\dfrac{e^{-\la}-1}{\la}
 \right).\leqno{(4.1b)}$$
\end{lemma}
\begin{proof}
(a) The equation (4.1a) follows from the formula $(3.9')$ and Claim~3.6a.
\smallskip

\noindent
(b) By using the formula (2.12) in the form
 $C(\la,\mu)=\dfrac{e^{\mu}-1}{\la\mu}\cdot
  \dfrac{\la+\mu}{e^{\la+\mu}-1}-\dfrac{1}{\la}$, one gets
$$\dfrac{\mu}{e^{\mu}-1}(1+\la C(\la,\mu))
 =\dfrac{\mu}{e^{\mu}-1}\left(
  \dfrac{e^{\mu}-1}{\mu}\cdot \dfrac{\la+\mu}{e^{\la+\mu}-1}\right)
 =\dfrac{\la+\mu}{e^{\la+\mu}-1}.$$
Now apply the symmetry (\emph{the third key point})
$$C(\la,\mu)=C(-\mu,-\la)=\dfrac{e^{-\la}-1}{\la\mu}\cdot
  \dfrac{-\la-\mu}{e^{-\la-\mu}-1}+\dfrac{1}{\mu},\mbox{ hence}$$
$$\dfrac{\la}{e^{\la}-1}(1-\mu C(\la,\mu))
 =\dfrac{\la}{e^{\la}-1}\left(
  \dfrac{e^{-\la}-1}{\la}\cdot \dfrac{\la+\mu}{e^{-\la-\mu}-1}\right)
 =\dfrac{1-e^{\la}}{e^{\la}-1}\cdot e^{\mu} \dfrac{\la+\mu}{1-e^{\la+\mu}}
 =e^{\mu} \dfrac{\la+\mu}{e^{\la+\mu}-1}.$$
Then the equation~(4.1a) converts to
$$e^{\mu} \dfrac{\la+\mu}{e^{\la+\mu}-1}f(\la,\mu)
 +\dfrac{\la+\mu}{e^{\la+\mu}-1}f(\mu,\rho)
 -\dfrac{\la+\mu}{e^{-\la-\mu}-1}f(\rho,\la)
 +\dfrac{e^{\mu}-1}{\la\mu}\cdot
  \dfrac{\la+\mu}{e^{\la+\mu}-1}-\dfrac{1}{\la}=0,\mbox{ or}$$
$$e^{\mu} \dfrac{\la+\mu}{e^{\la+\mu}-1}f(\la,\mu)
 +\dfrac{\la+\mu}{e^{\la+\mu}-1}f(\mu,\rho)
 +e^{\la+\mu}\dfrac{\la+\mu}{e^{\la+\mu}-1}f(\rho,\la)
 =\dfrac{1}{\la}\left(1-\dfrac{e^{\mu}-1}{\mu}\cdot
  \dfrac{\la+\mu}{e^{\la+\mu}-1} \right).$$
One can multiply both sides by $\dfrac{e^{\la+\mu}-1}{\la+\mu}$ since
 this series begins with 1, i.e.
$$e^{\mu}f(\la,\mu)+f(\mu,\rho)+e^{\la+\mu}f(\rho,\la)
 =\dfrac{1}{\la}\left(\dfrac{e^{\la+\mu}-1}{\la+\mu}
  -\dfrac{e^{\mu}-1}{\mu} \right).$$
Let us swap $\la$ and $\rho=-\la-\mu$ (in other words one substitutes
 $(-\la-\mu)$ for $\la$):
$$e^{\mu}f(-\la-\mu,\mu)+f(\mu,\la)+e^{-\la}f(\la,-\la-\mu)
 =\dfrac{1}{\la+\mu}\left(\dfrac{e^{-\la}-1}{\la}
  +\dfrac{e^{\mu}-1}{\mu} \right).$$
To get the equation (4.1b) it remains
 to use the symmetry $f(\la,\mu)=f(\mu,\la)$.
\end{proof}
\smallskip

\noindent
{\bf The hexagon part of Theorem~1.5b} follows from Lemma~4.1b:
 the compressed hexagon equation $\ol{(1.3b)}$
 is equivalent to $(4.1b)=(1.5b)$.
\qed
\medskip

An explicit form of the compressed hexagon~$\ol{(1.3b)}$
 is \emph{the 6th key point}.

\begin{lemma}
For the generating function $f(\la,\mu)$ of any compressed associator,
 one has
$$Even\bigl(f(\la,0)\bigr)=\dfrac{1}{2\la^2}
 \left(\dfrac{2\la}{e^{2\la}-1}+\la-1\right),\quad
 f(\la,-\la)=\dfrac{1}{\la^2}-\dfrac{2}{\la(e^{\la}-e^{-\la})}.
 \leqno{(4.2)}$$
In particular, one obtains the extreme coefficients
 $\al_{2k,0}=\dfrac{2^{2k+1}B_{2k+2}}{(2k+2)!}$ for every $k\geq 0$.
\end{lemma}
\noindent
\emph{Proof.}
Let us solve the equation (4.1b) explicitly for $\mu=0\;\Bigl($then
 set $\left.\dfrac{e^{\mu}-1}{\mu}\right|_{\mu=0}=1$ $\Bigr)$.
$$\mbox{One has }\quad f(\la,0)+f(0,-\la)+e^{-\la}f(\la,-\la)
 =\dfrac{1}{\la}\left( 1+\dfrac{e^{-\la}-1}{\la} \right)
 =\dfrac{e^{-\la}+\la-1}{\la^2}.$$
Since $Odd\bigl(f(\la,-\la)\bigr)=0$,
 $Odd\bigl(f(\la,0)+f(0,-\la)\bigr)=0$, and
$$Even\bigl(f(0,-\la)\bigr)=Even\bigl(f(\la,0)\bigr)
  =\sum\limits_{k=0}^{\infty}\al_{2k,0}\la^{2k},
  \mbox{ then one obtains}$$
$$2\sum\limits_{k=0}^{\infty}\al_{2k,0}\la^{2k}
 +e^{-\la}Even\bigl(f(\la,-\la)\bigr)
 =\dfrac{e^{-\la}+\la-1}{\la^2}.$$
By substituting $(-\la)$ for $\la$ and using
 $Even\bigl(f(-\la,\la)\bigr)=Even\bigl(f(\la,-\la)\bigr)$, one gets
$$2\sum\limits_{k=0}^{\infty}\al_{2k,0}\la^{2k}
 +e^{\la}Even\bigl(f(\la,-\la)\bigr)
 =\dfrac{e^{\la}-\la-1}{\la^2}.$$
From the two above equations one deduces
$$f(\la,-\la)=Even\bigl(f(\la,-\la)\bigr)
 =\dfrac{1}{\la^2}-\dfrac{2}{\la(e^{\la}-e^{-\la})}\quad\mbox{ and}$$
$$Even\bigl(f(\la,0)\bigr)=\sum\limits_{k=0}^{\infty}\al_{2k,0}\la^{2k}
 =\dfrac{1}{2\la^2}\left(\dfrac{2\la}{e^{2\la}-1}+\la-1\right)
 =\dfrac{1}{2\la^2}\sum_{n=1}^{\infty}\dfrac{B_{2n}}{(2n)!}(2\la)^{2n}.
 \eqno{\square}$$

\begin{example}
Note that $\al_{2k,0}=\dfrac{2^{2k+1}B_{2k+2}}{(2k+2)!}$
 are correct coefficients of the logarithm
 of any honest Drinfeld associator
 $\Phi(a,b)=\exp\bigl(\ph(a,b)\bigr)\in\hat L_3$,
 not only the compressed one.

This form of the extreme coefficients $\al_{2k,0}$
 is \emph{the 7th key point}:
$$\al_{00}=\dfrac{2B_2}{2!}=\dfrac{1}{6},\quad
  \al_{20}=\dfrac{2^3B_4}{4!}=-\dfrac{1}{90},\quad
  \al_{40}=\dfrac{2^5B_6}{6!}=\dfrac{1}{945},\quad
  \al_{60}=\dfrac{2^7B_8}{8!}=-\dfrac{1}{9450}.$$
For the coefficients $\ga_k$ of the series
 $\dfrac{2\la}{e^{\la}-e^{-\la}}=\sum\limits_{k=0}^{\infty}\ga_k\la^{2k}$,
 similarly to Lemma~2.2 one can find the recursive formula
 $\sum\limits_{k=0}^n\dfrac{\ga_{n-k}}{(2k+1)!}=0$ for each
 $n>1$, $\ga_0=0$.
Hence
$$\dfrac{2\la}{e^{\la}-e^{-\la}}=
 1-\dfrac{1}{6}\la^2+\dfrac{7}{360}\la^4-\dfrac{31}{3\cdot 7!}\la^6
 +\dfrac{127}{15\cdot 8!}\la^8
 +\cdots\leqno{(4.3)}$$
For the generating function of any compressed Drinfeld associator
 $\bar\ph\in\hat{\bar L}_3$, one gets
$$f(\la,-\la)=\dfrac{1}{6}
 -\dfrac{7}{360}\la^2+\dfrac{31}{15120}\la^4
 -\dfrac{127}{604800}\la^6
 +\cdots\eqno{\bt}$$
\end{example}

\begin{lemma}
Under the conditions of Lemma~4.1,
 set $\ti f(\la,\mu):=1+\la\mu f(\la,\mu)$.
\smallskip

\noindent
(a) The new function $\ti f(\la,\mu)$ obeys the same symmetry:
 $\ti f(\la,\mu)=\ti f(\mu,\la)$.
Moreover, the function $f(\la,\mu)$ is even (i.e.
 $f(\la,\mu)=f(-\la,-\mu)$) if and only if $\ti f(\la,\mu)$ is even.
\smallskip

\noindent
(b) The compressed hexagon $\ol{(1.3b)}=(4.1b)$ is equivalent
 to the following equation in $\C[[\la,\mu]]$
$$(\la+\mu)\ti f(\la,\mu)
 =\la e^{\mu}\ti f(\mu,-\la-\mu)
 +\mu e^{-\la}\ti f(\la,-\la-\mu).\leqno{(4.4)}$$
\end{lemma}
\begin{proof}
The item (a) follows from the definition of $\ti f(\la,\mu)$.
\smallskip

\noindent
(b) The equation~(4.1b) can be rewritten as follows:
$$\left( f(\la,\mu)+\dfrac{1}{\la\mu} \right)
 +e^{\mu}\left( f(\mu,-\la-\mu)-\dfrac{1}{\mu(\la+\mu)} \right)
 +e^{-\la}\left( f(\la,-\la-\mu)-\dfrac{1}{\la(\la+\mu)} \right)=0.$$
To get (4.4) it remains to multiply by $\la\mu(\la+\mu)$.
\end{proof}

Recall that the even and odd parts of a series were introduced
 before Theorem~1.5.

\begin{lemma}
In the notations of Lemma~4.4,
 the compressed hexagon $\ol{(1.3b)}=(4.4)$ can be split into
 the two following equations in the algebra $\C[[\la,\mu]]$
$$(\la+\mu)Even\bigl(\ti f(\la,\mu)\bigr)
 =\la e^{\mu}Even\bigl(\ti f(\mu,-\la-\mu)\bigr)
 +\mu e^{-\la}Even\bigl(\ti f(\la,-\la-\mu)\bigr);\leqno{(4.5a)}$$
$$Odd\bigl(f(\la,\mu)\bigr)+e^{\mu}Odd\bigl(f(\mu,-\la-\mu)\bigr)
 +e^{-\la}Odd\bigl(f(\la,-\la-\mu)\bigr)=0.\leqno{(4.5b)}$$
\end{lemma}
\begin{proof}
Let us substitute the pair $(-\la,-\mu)$ for $(\la,\mu)$ into
 the equation~(4.4)
$$-(\la+\mu)\ti f(-\la,-\mu)
 =-\la e^{-\mu}\ti f(-\mu,\la+\mu)
 -\mu e^{\la}\ti f(-\la,\la+\mu).$$
Now swap the arguments $\la$, $\mu$ and use the symmetry
 $\ti f(\la,\mu)=\ti f(\mu,\la)$
$$-(\la+\mu)\ti f(-\la,-\mu)
 =-\mu e^{-\la}\ti f(-\la,\la+\mu)
 -\la e^{\mu}\ti f(-\mu,\la+\mu).$$
If one subtracts the latter equation from (4.4), then
 after dividing by 2 one obtains the equation (4.5a).
If one adds the latter equation to (4.4), then one has
$$(\la+\mu)Odd\bigl(\ti f(\la,\mu)\bigr)
 =\la e^{\mu}Odd\bigl(\ti f(\mu,-\la-\mu)\bigr)
 +\mu e^{-\la}Odd\bigl(\ti f(\la,-\la-\mu)\bigr).$$
Since $Odd\bigl(\ti f(\la,\mu)\bigr)=\la\mu Odd\bigl(f(\la,\mu)\bigr)$,
 after dividing by $\la\mu(\la+\mu)$ one gets~(4.5b).
\end{proof}

The splitting of Lemma~4.5 is \emph{the 8th key point}.


\subsection{Explicit description of all even compressed associators}

\begin{lemma}
The general solution of the equation~(4.5a) is
$$Even\bigl(\ti f(\la,\mu)\bigr)=1+\la\mu\cdot Even\bigl(f(\la,\mu)\bigr)
 =\dfrac{e^{\la+\mu}-e^{-\la-\mu}}{2(\la+\mu)}h(\la,\mu),\leqno{(4.6)}$$
 where $h(\la,\mu)$ is a function satisfying the boundary
 condition $h(\la,0)=\dfrac{2\la}{e^{\la}-e^{-\la}}$ and
 the symmetry relations
 $h(\la,\mu)=h(\mu,\la)=h(-\la,-\mu)=h(\la,-\la-\mu)$.
\end{lemma}
\begin{proof}
Let us swap the arguments $\la,\mu$ in~(4.5a).
Due to $\ti f(\la,\mu)=\ti f(\mu,\la)$ one gets
$$(\la+\mu)Even\bigl(\ti f(\la,\mu)\bigr)
 =\mu e^{\la}Even\bigl(\ti f(\la,-\la-\mu)\bigr)
 +\la e^{-\mu}Even\bigl(\ti f(\mu,-\la-\mu)\bigr).$$
By subtracting the above equation from (4.5a) one obtains
$$\la(e^{\mu}-e^{-\mu})\cdot Even\bigl(\ti f(\mu,-\la-\mu)\bigr)
 +\mu(e^{-\la}-e^{\la})\cdot Even\bigl(\ti f(\la,-\la-\mu)\bigr)=0,\mbox{ or}$$
$$Even\bigl(\ti f(\mu,-\la-\mu)\bigr):\left(\dfrac{e^{\la}-e^{-\la}}{2\la}\right)
 =Even\bigl(\ti f(\la,-\la-\mu)\bigr):\left(\dfrac{e^{\mu}-e^{-\mu}}{2\mu}\right).$$
Introduce the function
$$h(\la,\mu):=
 Even\bigl(\ti f(\mu,-\la-\mu)\bigr):\left(\dfrac{e^{\la}-e^{-\la}}{2\la}\right)
=Even\bigl(\ti f(\la,-\la-\mu)\bigr):\left(\dfrac{e^{\mu}-e^{-\mu}}{2\mu}\right).
 \leqno{(4.6')}$$
Then the new function is even and symmetric:
 $h(\la,\mu)=h(\mu,\la)=h(-\la,-\mu)$.
Let us substitute the expressions
$$Even\bigl(\ti f(\mu,-\la-\mu)\bigr)
 =\left(\dfrac{e^{\la}-e^{-\la}}{2\la}\right)\cdot h(\la,\mu),\quad
  Even\bigl(\ti f(\la,-\la-\mu)\bigr)
 =\left(\dfrac{e^{\mu}-e^{-\mu}}{2\mu}\right)\cdot h(\la,\mu)$$
 into the equation (4.5a).
One has
$$(\la+\mu)Even\bigl(\ti f(\la,\mu)\bigr)
 =\left(\la e^{\mu}\dfrac{e^{\la}-e^{-\la}}{2\la}
 +\mu e^{-\la}\dfrac{e^{\mu}-e^{-\mu}}{2\mu}\right)h(\la,\mu)
 =\dfrac{e^{\la+\mu}-e^{-\la-\mu}}{2}h(\la,\mu).$$
So, (4.6) is proved.
Let us substitute $(-\la-\mu)$ for $\mu$ into the above equation:
$$Even\bigl(\ti f(\la,-\la-\mu)\bigr)
 =\dfrac{e^{-\mu}-e^{\mu}}{2(-\mu)}h(\la,-\la-\mu)
 =\dfrac{e^{\mu}-e^{-\mu}}{2\mu}h(\la,-\la-\mu).$$
By comparing the latter formula with the equation~$(4.6')$,
 one gets $h(\la,\mu)=h(\la,-\la,-\mu)$.
To obtain the condition
 $h(\la,0)=\dfrac{2\la}{e^{\la}-e^{-\la}}$
 it remains to substite $\mu=0$ into~(4.6).
\end{proof}

Observe that the conditions of Lemma~4.6 imply also
 $h(\la,-\la)=h(\la,0)=\dfrac{2\la}{e^{\la}-e^{-\la}}$.
This agrees with the formula~(4.2):
 $h(\la,-\la)\st{(4.6)}{=}1-\la^2f(\la,-\la)
 \st{(4.2)}{=}\dfrac{2\la}{e^{\la}-e^{-\la}}$.
So, the compressed hexagon~$\ol{(1.3b)}$ was reduced to
 $h(\la,\mu)=h(\la,-\la-\mu)$, this is \emph{the 9th key point}.
To describe all series $h(\la,\mu)$ with this symmetry one needs
 associator polynomials.

\begin{definition}[associator polynomials]
For each $n\geq 0$, a homogeneous polynomial
 $F_n(\la,\mu)=\sum\limits_{k=0}^{n}\de_k \la^k\mu^{n-k}$
 with $\de_k\in\C$ is called \emph{an associator polynomial},
 if $F_n(\la,\mu)=F_n(\mu,\la)=F_n(\la,-\la-\mu)$
 hold for all $\la,\mu\in\C$.
The Drinfeld series $s(\la,\mu)$ from Definition~6.1b
 contains associator polynomials
 $F^D_{2n+1}(\la,\mu)=(\la+\mu)^{2n+1}-\la^{2n+1}-\mu^{2n+1}$.
\ed
\end{definition}

\begin{example}
In degrees $0,2,4$ an associator polynomial is unique up to factor:
$$F_0(\la,\mu)=1,\quad
  F_2(\la,\mu)=\la^2+\la\mu+\mu^2,\quad
  F_4(\la,\mu)=\la^4+2\la^3\mu+3\la^2\mu^2+2\la\mu^3+\mu^4.$$
In each odd degree $n=1,3,5$
 there is also exactly one associator polynomial up to factor:
$$F_1(\la,\mu)=0,\quad
  F_3(\la,\mu)=\la^2\mu+\la\mu^2,\quad
  F_5(\la,\mu)=\la^4\mu+2\la^3\mu^2+2\la^2\mu^3+\la\mu^5.$$
Rather surprisingly, that an associator polynomial of the degree $n$
 is not unique in general, for instance, in degree 6.
One can check by hands that in degree 6 there is a 1-parametric
 family of associator polynomials:
$$F_6(\la,\mu)=\la^6+3\la^5\mu+\de\la^4\mu^2+
 (2\de-5)\la^3\mu^3+\de\la^2\mu^4+3\la\mu^5+\mu^6,\quad \de\in\C.$$
However, in degree 7 there is a unique associator
 polynomial up to factor:
$$F_7(\la,\mu)=
 \la^6\mu+3\la^5\mu^2+5\la^4\mu^3+5\la^3\mu^4+3\la^2\mu^5+\la\mu^6.
 \eqno{\bt}$$
\end{example}

\begin{lemma}
Any function $\bar h(\la,\mu)$ satisfying the boundary conditions
 $\bar h(\la,0)=\sum\limits_{n=0}^{\infty}\ga_n\la^{2n}$
 and the symmetry relations
 $\bar h(\la,\mu)=\bar h(\mu,\la)=\bar h(\la,-\la-\mu)$
 is $\bar h(\la,\mu)=\sum\limits_{n=0}^{\infty} \ga_n F_n(\la,\mu)$,
 where $F_n(\la,\mu)$ is an associator polynomial
 with the extreme coefficient 1 (respectively, 0)
 for each even (respectively, odd) $n\geq 0$.
\end{lemma}

\noindent
\emph{Proof} follows from Definition~4.7.
The relation $F_{2n+1}(\la,\mu)=F_{2n+1}(\la,-\la-\mu)$
 implies that the extreme coefficient of any associator polynomial
 $F_{2n+1}(\la,\mu)$ is always 0.
\qed
\smallskip

Non-uniqueness of even compressed associators will follow
 from non-uniqueness of associator polynomials
 (\emph{the 10th key point}), see the proof of Theorem~1.5c.

\begin{proposition}
(a) For each $n\geq 0$, any associator polynomial of the degree $2n$ is
$$F_{2n}(\la,\mu)=\sum_{k=0}^{[\frac{n}{3}]}
 \be_{nk} \la^{2k}\mu^{2k}(\la+\mu)^{2k}(\la^2+\la\mu+\mu^2)^{n-3k},
 \leqno{(4.10a)}$$
 where $\be_{nk}\in\C$ are free parameters for $0\leq k\leq [\frac{n}{3}]$.
\medskip

\noindent
(b) For every $n\geq 1$, any associator polynomial $F_{2n+1}(\la,\mu)$
 has the following form:
$$F_{2n+1}(\la,\mu)=\sum_{k=0}^{[\frac{n-1}{3}]}
 \ti\be_{nk} \la^{2k+1}\mu^{2k+1}(\la+\mu)^{2k+1}
  (\la^2+\la\mu+\mu^2)^{n-3k-1},\leqno{(4.10b)}$$
 where $\ti\be_{nk}\in\C$ are free parameters for
 $0\leq k\leq [\frac{n-1}{3}]$.
\end{proposition}
\begin{proof}
(a) Induction on $n$.
The bases $n=0,1,2$ are in Example~4.8.
\smallskip

Induction step goes down from $n$ to $n-3$.

Suppose that $F_{2n}(\la,\mu)$ is an associator polynomial of the degree $2n$
 with the extreme coefficient $\be_{n0}$.
Then the polynomial
 $\ti F_{2n}(\la,\mu)=F_{2n}(\la,\mu)-\be_{n0}(\la^2+\la\mu+\mu^2)^n$
 satisfies Definition~4.7.
The relations $F_{2n}(\la,\mu)=F_{2n}(\mu,\la)=F_{2n}(\la,-\la-\mu)$ imply
 $$F_{2n}(\la,\mu)=\be_{n0}\la^{2n}+n\be_{n0}\la^{2n-1}\mu+\cdots
 +n\be_{n0}\la\mu^{2n-1}+\be_{n0}\mu^{2n}.$$
The polynomial $\be_{n0}(\la^2+\la\mu+\mu^2)^n$
 has the same form.
Hence the first two (and the last two) coefficients of
 $\ti F_{2n}(\la,\mu)$ are always zero.

One gets $\ti F_{2n}(\la,\mu)=\la^2\mu^2\ol{F}_{2n-4}(\la,\mu)$ for
 a polynomial $\ol{F}_{2n-4}(\la,\mu)$ of the degree $2n-4$
 such that $\ol{F}_{2n-4}(\la,\mu)=\ol{F}_{2n-4}(\mu,\la)$ and
 $\la^2\mu^2\ol{F}_{2n-4}(\la,\mu)
 =\la^2(-\la-\mu)^2\ol{F}_{2n-4}(\la,-\la-\mu)$.
\smallskip

The equation
 $\mu^2\ol{F}_{2n-4}(\la,\mu)=(\la+\mu)^2\ol{F}_{2n-4}(\la,-\la-\mu)$
 implies that there is a polynomial $\ol{\ol{F}}_{2n-6}(\la,\mu)$
 of the degree $2n-6$ such that
 $\ol{F}_{2n-4}(\la,\mu)=(\la+\mu)^2\ol{\ol{F}}_{2n-6}(\la,\mu)$.
\smallskip

Moreover, the new polynomial $\ol{\ol{F}}_{2n-6}(\la,\mu)$ satisfies
 all the conditions of Definition~4.7, i.e.
 $\ol{\ol{F}}_{2n-6}(\la,\mu)$
 is equal to an associator polynomial $F_{2n-6}(\la,\mu)$.
Hence
$$F_{2n}(\la,\mu)-\be_{n0}(\la^2+\la\mu+\la\mu)^n
 =\la^2\mu^2(\la+\mu)^2 F_{2n-6}(\la,\mu).$$
The induction hypothesis
$$F_{2n-6}(\la,\mu)=\sum_{k=0}^{[\frac{n}{3}]-1}
 \be_{n-3,k}\la^{2k}\mu^{2k}(\la+\mu)^{2k}(\la^2+\la\mu+\mu^2)^{n-3k-3}
 \qquad \mbox{ implies}$$
$$F_{2n}(\la,\mu)=\be_{n0}(\la^2+\la\mu+\mu^2)^n
 +\la^2\mu^2(\la+\mu)^2 \sum_{k=0}^{[\frac{n}{3}]-1}
 \be_{n-3,k} \la^{2k}\mu^{2k}(\la+\mu)^{2k}(\la^2+\la\mu+\mu^2)^{n-3-3k}.$$
To get (4.10a) it remains to set
 $\be_{n,k+1}:=\be_{n-3,k}$ for $0\leq k\leq[\frac{n}{3}]-1$.
\medskip

\noindent
(b) The proof is analogous to the item~(a).
If $F_{2n+1}(\la,\mu)$ is an associator polynomial, then
 its extreme coefficient is zero, i.e.
 $F_{2n+1}(\la,\mu)=\la\mu\ol{F_{2n-1}}(\la,\mu)$
 for some polynomial $\ol{F}_{2n-1}(\la,\mu)$ of the degree $2n-1$
 with the properties
$$\ol{F}_{2n-1}(\la,\mu)=\ol{F}_{2n-1}(\mu,\la)\quad \mbox{ and }\quad
  \la\mu\ol{F}_{2n-1}(\la,\mu)=\la(-\la-\mu)\ol{F}_{2n-1}(\la,-\la-\mu).$$
The equation
 $\mu\ol{F}_{2n-1}(\la,\mu)=-(\la+\mu)\ol{F}_{2n-1}(\la,-\la-\mu)$
 implies that there is a polynomial $\ol{\ol{F}}_{2n-2}(\la,\mu)$
 of the degree $2n-2$ such that
 $\ol{F}_{2n-1}(\la,\mu)=(\la+\mu)\ol{\ol{F}}_{2n-2}(\la,\mu)$.
Moreover, the new polynomial $\ol{\ol{F}}_{2n-2}(\la,\mu)$ satisfies
 all the conditions of Definition~4.7, i.e.
 $\ol{\ol{F}}_{2n-2}(\la,\mu)$
 is equal to an associator polynomial $F_{2n-2}(\la,\mu)$.

To get (4.10b) it remains to apply the formula (4.10a) for this
 polynomial $F_{2n-2}(\la,\mu)$ and to set
 $\ti\be_{nk}:=\be_{n-1,k}$ for $0\leq k\leq[\frac{n-1}{3}]$.
\end{proof}
\smallskip

\emph{The 11th key point}:
 the symmetry $F_n(\la,\mu)=F_n(\la,-\la-\mu)$ has led to~(4.10a) and (4.10b).

\begin{example}
By Propositon 4.10 the number of free parameters, on which the family
 $F_{n}(\la,\mu)$ depends, increases by 1 when $n$ increases by 3.
The first six (starting with 0) associator polynomials
 are unique up to factor, but not the seventh $F_6(\la,\mu)$.
One gets
$$\begin{array}{l}
F_6(\la,\mu)=\be_{30}(\la^2+\la\mu+\mu^2)^3
 +\be_{31}\la^{2}\mu^2(\la+\mu)^{2},\\
F_8(\la,\mu)=\be_{40}(\la^2+\la\mu+\mu^2)^4
 +\be_{41}\la^{2}\mu^2(\la+\mu)^{2}(\la^2+\la\mu+\mu^2),\\
F_{10}(\la,\mu)=\be_{50}(\la^2+\la\mu+\mu^2)^5
 +\be_{51}\la^{2}\mu^2(\la+\mu)^{2}(\la^2+\la\mu+\mu^2)^2,\\
F_{12}(\la,\mu)=\be_{60}(\la^2+\la\mu+\mu^2)^6
 +\be_{61}\la^{2}\mu^2(\la+\mu)^{2}(\la^2+\la\mu+\mu^2)^3
 +\be_{62}\la^{4}\mu^4(\la+\mu)^{4};
\end{array}$$
$$\begin{array}{l}
F_7(\la,\mu)=\ti\be_{30}\la\mu(\la+\mu)(\la^2+\la\mu+\mu^2)^2, \\
F_9(\la,\mu)=\ti\be_{40}\la\mu(\la+\mu)(\la^2+\la\mu+\mu^2)^3
 +\ti\be_{41}\la^{3}\mu^3(\la+\mu)^{3},\\
F_{11}(\la,\mu)=\ti\be_{50}\la\mu(\la+\mu)(\la^2+\la\mu+\mu^2)^4
 +\ti\be_{51}\la^{3}\mu^3(\la+\mu)^{3}(\la^2+\la\mu+\mu^2),\\
F_{13}(\la,\mu)=\ti\be_{60}\la\mu(\la+\mu)(\la^2+\la\mu+\mu^2)^5
 +\ti\be_{61}\la^{3}\mu^3(\la+\mu)^{3}(\la^2+\la\mu+\mu^2)^2.
\end{array}$$
One can check that the parameter $\be_{31}$
 is related with $\de$ from Example~4.8 as $\de=\be_{31}+6$.
Up to degree~10 the function $h(\la,\mu)$ from Lemma~4.6 is
$$h(\la,\mu)=1-\dfrac{1}{6}F_2(\la,\mu)+\dfrac{7}{360}F_4(\la,\mu)
 -\dfrac{31}{3\cdot 7!}F_6(\la,\mu)
 +\dfrac{127}{15\cdot 8!}F_8(\la,\mu)
 +\cdots\eqno{\bt}$$
\end{example}
\smallskip


\subsection{Description of all odd compressed associators}

\begin{lemma}
The general solution of the equation (4.5b) is
$$Odd\bigl(f(\la,\mu)\bigr)
 =\dfrac{e^{\la+\mu}-e^{-\la-\mu}}{2}\ti h(\la,\mu),
 \mbox{ where}\leqno{(4.12)}$$
 $\ti h(\la,\mu)$ is a function satisfying
 the relations
 $\ti h(\la,\mu)=\ti h(\mu,\la)=\ti h(-\la,-\mu)=\ti h(\la,-\la-\mu)$.
\end{lemma}
\noindent
\emph{Proof} is analogous to the proof of Lemma~4.6.
Let us swap the arguments $\la,\mu$ in the equation (4.5b).
Due to the symmetry $f(\la,\mu)=f(\mu,\la)$ one gets
$$Odd\bigl(f(\la,\mu)\bigr)
 +e^{\la}Odd\bigl(f(\la,-\la-\mu)\bigr)
 +e^{-\mu}Odd\bigl(f(\mu,-\la-\mu)\bigr)=0.$$
By subtracting the above equation from (4.5b) one obtains
$$(e^{\mu}-e^{-\mu})\cdot Odd\bigl(f(\mu,-\la-\mu)\bigr)
 +(e^{-\la}-e^{\la})\cdot Odd\bigl(f(\la,-\la-\mu)\bigr)=0,\mbox{ or}$$
$$Odd\bigl(f(\mu,-\la-\mu)\bigr):\left(\dfrac{e^{\la}-e^{-\la}}{2}\right)
 =Odd\bigl(f(\la,-\la-\mu)\bigr):\left(\dfrac{e^{\mu}-e^{-\mu}}{2}\right).$$
Introduce the function
$$\ti h(\la,\mu):=
 -Odd\bigl(f(\mu,-\la-\mu)\bigr):\left(\dfrac{e^{\la}-e^{-\la}}{2}\right)
=-Odd\bigl(f(\la,-\la-\mu)\bigr):\left(\dfrac{e^{\mu}-e^{-\mu}}{2}\right).
 \leqno{(4.12')}$$
Then the new function is even and symmetric:
 $\ti h(\la,\mu)=\ti h(\mu,\la)=\ti h(-\la,-\mu)$.
Let us substitute the expressions
$$Odd\bigl(f(\mu,-\la-\mu)\bigr)
 =-\left(\dfrac{e^{\la}-e^{-\la}}{2}\right)\cdot \ti h(\la,\mu),\quad
  Odd\bigl(f(\la,-\la-\mu)\bigr)
 =-\left(\dfrac{e^{\mu}-e^{-\mu}}{2}\right)\cdot \ti h(\la,\mu)$$
 into the equation (4.5b).
One has
$$Odd\bigl(f(\la,\mu)\bigr)
 =\left(e^{\mu}\dfrac{e^{\la}-e^{-\la}}{2}
 +e^{-\la}\dfrac{e^{\mu}-e^{-\mu}}{2}\right) \ti h(\la,\mu)
 =\dfrac{e^{\la+\mu}-e^{-\la-\mu}}{2} \ti h(\la,\mu).$$
So, (4.12) is proved.
Let us substitute $(-\la-\mu)$ for $\mu$ into the above equation:
$$Odd\bigl(f(\la,-\la-\mu)\bigr)
 =-\dfrac{e^{\mu}-e^{-\mu}}{2}\ti h(\la,-\la-\mu).$$
By comparing the latter formula with~$(4.12')$,
 one gets $\ti h(\la,\mu)=\ti h(\la,-\la,-\mu)$.
\qed
\bigskip

\noindent
{\bf Proof of Theorem~1.5c.}
By Lemmas~4.6, 4.9, 4.12, and Proposition~4.10a,
 the general solution of the equation~(1.5b)=(4.1b) is
 $f(\la,\mu)=Even\bigl(f(\la,\mu)\bigr)+Odd\bigl(f(\la,\mu)\bigr)$, where
$$\begin{array}{l}
1+\la\mu\cdot Even\bigl(f(\la,\mu)\bigr)
=\dfrac{e^{\la+\mu}-e^{-\la-\mu}}{2(\la+\mu)} \left(
 \sum\limits_{n=0}^{\infty}\sum\limits_{k=0}^{[\frac{n}{3}]}
 \be_{nk} \la^{2k}\mu^{2k}(\la+\mu)^{2k}\om^{2n-6k}
  \right),\\ \\
Odd\bigl(f(\la,\mu)\bigr)
=\dfrac{e^{\la+\mu}-e^{-\la-\mu}}{2} \left(
 \sum\limits_{n=0}^{\infty}\sum\limits_{k=0}^{[\frac{n}{3}]}
 \ti\be_{nk} \la^{2k}\mu^{2k}(\la+\mu)^{2k}\om^{2n-6k} \right),\;
 \om=\sqrt{\la^2+\la\mu+\mu^2}. \end{array}$$
Let us substitute $\mu=-\la$ (then $\om=\la$) into the former equation
 and apply Lemma~4.2:
$$\sum\limits_{n=0}^{\infty} \be_{n0}\la^{2n}
 =1-\la^2 Even\bigl(f(\la,-\la)\bigr)
 \st{(4.2)}{=}\dfrac{2\la}{e^{\la}-e^{-\la}},\mbox{ hence }
 \sum\limits_{n=0}^{\infty}\be_{n0}\om^{2n}=\dfrac{2\om}{e^{\om}-e^{-\om}}$$
 as required.
Let us relate $\al_{2n+1,0}$ with $\ti\be_{n0}$.
One has
 $Odd\bigl(f(\la,\mu)\bigr)=\sum\limits_{n=0}^{\infty}\sum\limits_{k=0}^{2n+1}
  \al_{2n+1,k}\la^{2n-k+1}\mu^k$.
Hence one obtains
 $Odd(f(\la,0))=\sum\limits_{n=0}^{\infty}\al_{2n+1,0}\la^{2n+1}$.
On the other hand, for $\mu=0$, one gets
$$\dfrac{e^{\la+\mu}-e^{-\la-\mu}}{2}
 \sum\limits_{n=0}^{\infty}\sum\limits_{k=0}^{[\frac{n}{3}]}
 \ti\be_{nk}\la^{2k}\mu^{2k}(\la+\mu)^{2k}\om^{2n-6k}
 =\dfrac{e^{\la}-e^{-\la}}{2}
 \sum_{n=0}^{\infty}\ti\be_{n0}\la^{2n}
 =\sum\limits_{n=0}^{\infty}\al_{2n+1,0}\la^{2n+1}.\eqno{\qed}$$
\smallskip

\noindent
{\bf Proof of Corollaries~1.6a-b}.
(a) The first series $f^I(\la,\mu)$ is obtained from the general
 formula~(1.5c) by taking $\be_{nk}=\ti\be_{nk}=0$.
This solution is related with the associator polynomials
 $F^I_{2n}(\la,\mu)=(\la^2+\la\mu+\mu^2)^n$.
(b) The second series $f^{II}(\la,\mu)$ appears due to
 the associator polynomials
 $F^{II}_{2n}(\la,\mu)=\dfrac{(\la+\mu)^{2n}+\la^{2n}+\mu^{2n}}{2}$,
 $n>1$.
Actually, let $\ga_k$ be the coefficients of the series
 $\dfrac{2\la}{e^{\la}-e^{-\la}}
 =\sum\limits_{k=0}^{\infty}\ga_k\la^{2k}$, see~(4.3).
Then the series $1+\sum\limits_{k=1}^{\infty}\ga_kF^{II}_{2k}(\la,\mu)$
 plays the role of the function $h(\la,\mu)$ from Lemma~4.6.
It remains to compute
$$1+\la\mu\cdot f^{II}(\la,\mu)
 =\dfrac{e^{\la+\mu}-e^{-\la-\mu}}{2(\la+\mu)} \left(
  1+\sum\limits_{k=1}^{\infty}\ga_k
  \dfrac{(\la+\mu)^{2n}+\la^{2n}+\mu^{2n}}{2} \right)=$$
$$=\dfrac{e^{\la+\mu}-e^{-\la-\mu}}{2(\la+\mu)}\cdot\dfrac{1}{2}
 \left( \dfrac{2(\la+\mu)}{e^{\la+\mu}-e^{-\la-\mu}}
 +\dfrac{2\la}{e^{\la}-e^{-\la}}+\dfrac{2\mu}{e^{\mu}-e^{-\mu}} -1 \right)=$$
$$=\dfrac{1}{2}+\dfrac{1}{2}\cdot\dfrac{e^{\la+\mu}-e^{-\la-\mu}}{2(\la+\mu)}
 \left( \dfrac{2\la}{e^{\la}-e^{-\la}}+\dfrac{2\mu}{e^{\mu}-e^{-\mu}}
 -1 \right)\quad\Lra\quad(1.6b).\eqno{\qed}$$
\smallskip


\section{Compressed pentagon equation}

This section is devoted to the proof of the pentagon part
 of Theorem~1.5b.
It turns out that if all commutators commute, then
 the compressed pentagon~$\ol{(1.3c)}$ follows from the symmetry
 $\al_{kl}=\al_{lk}$ (\emph{the 12th key point}), see Proposition~5.10.


\subsection{Generators and relations of the quotient $\bar L_4$}
\noindent

\noindent
Here one studies the quotient $\bar L_4$, where
 the compressed pentagon~$\ol{(1.3c)}=(5.9)$ lives.

\begin{definition}[alphabet $\LL$, algebras $L_4$ and $\bar L_4$,
 simple and non-simple commutators]
\noindent
\smallskip

\noindent
(a) Let the Lie algebra $L_4$ be generated by the letters of the alphabet
$$\LL=\{\;a:=t^{12},\; b:=t^{23},\; c:=t^{13},\;
          d:=t^{24},\; e:=t^{34},\; v:=t^{14}\; \}
 \mbox{ and the relations}$$
$$[a,e]=[b,v]=[c,d]=0\mbox{ and }
\left\{ \begin{array}{ll}
 x:=[a,b]=[b,c]=[c,a],&
 y:=[a,d]=[d,v]=[v,a],\\
 z:=[b,e]=[e,d]=[d,b],&
 u:=[c,e]=[e,v]=[v,c].
 \end{array}\right.$$
The Lie algebra $L_4$ is graded by $\deg(s)=1$ for any letter $s\in\LL$.
Put $\bar L_4=L_4/\bigl[[L_4,L_4],[L_4,L_4]\bigr]$.
By $\wh{\bar L}_4$ denote the algebra of formal series of
 elements from $\bar L_4$.
\medskip

\noindent
(b) Let $w$ be a word in the alphabet $\LL=\{t^{ij},\; 1\leq i<j\leq 4\}$.
Let $I_w$ be the set forming by the upper indices of the letters
 $t^{ij}$ including in $w$.
If the set $I_w$ contains at most three different indices,
 then the commutator $[w]\in\bar L_4$ is called \emph{simple},
 otherwise $[w]$ is \emph{non-simple}.
For example, the commutators $[aab]$ and $[vad]$ are simple, but
 $[dab]$ is not simple.
\ed
\end{definition}

\begin{claim}
The following relations hold in the quotient $\bar L_4$
$$\begin{array}{l}
(a)\quad
 [a+b+c,x]=0,\quad
 {[a+d+v,y]}=0,\quad
 {[b+e+d,z]}=0,\quad
 {[c+e+v,u]}=0;\\
(b)\quad
 [da]+[ea]+[va]=0,\hspace{13mm}
 {[db]}+[eb]+[vb]=0,\hspace{14mm}
 {[dx]}+[ex]+[vx]=0;\\
(c)\quad
 [dx]=-[cy]=-[cz]=[du],\;
 {[ex]}=-[ey]=-[az]=[au],\;
 {[vx]}=-[by]=-[vz]=[bu].
\end{array}$$
\end{claim}
\begin{proof}
(a) By Definition~5.1a one gets
 $[a+b+c,a]=[c,a]-[a,b]\st{(5.1a)}{=}0$ and
 similarly $[a+b+c,b]\st{(5.1a)}{=}0$, hence $[a+b+c,x]=0$.
The other relations are proved analogously.
\medskip

\noindent
(b) By Definition~5.1a one has
 $[d+e+v,a]=\bigl(-[a,d]+[v,a]\bigr)-[a,e]\st{(5.1a)}{=}0$ and
 similarly $[d+e+v,b]\st{(5.1a)}{=}0$, hence $[d+e+v,x]=0$.
\medskip

\noindent
(c) Definition~5.1a, the Jacobi identity $(J)$, and $[cd]=0$ are
 used below.
$$\begin{array}{ll}
{[dx]}\st{(5.1a)}{=}[dbc]\st{(J)}{=}-[bcd]-[cdb]\st{(5.1a)}{=}-[cz],&
{[dx]}\st{(5.1a)}{=}[dab]\st{(J)}{=}-[abd]-[bda]\st{(5.1a)}{=}[az]+[by],\\
{[dx]}\st{(5.1a)}{=}[dca]\st{(J)}{=}-[adc]-[cad]\st{(5.1a)}{=}-[cy],&
{[cy]}\st{(5.1a)}{=}[cva]\st{(J)}{=}-[vac]-[acv]\st{(5.1a)}{=}[vx]+[au],\\
{[cy]}\st{(5.1a)}{=}[cdv]\st{(J)}{=}-[dvc]-[vcd]\st{(5.1a)}{=}-[du],&
{[cz]}\st{(5.1a)}{=}[cbe]\st{(J)}{=}-[bec]-[ecb]\st{(5.1a)}{=}[bu]+[ex].
\end{array}$$

\noindent
One gets $[dx]=-[cy]=-[cz]=[du]=[az]+[by]=-[bu]-[ex]=-[vx]-[au]$, i.e.
 $[bu]=-[dx]-[ex]\st{(5.2b)}{=}[vx]$ and
 $[au]=-[dx]-[vx]\st{(5.2b)}{=}[ex]$.
Similarly, by $[ae]=0$ one has
$$[ex]\st{(5.1a)}{=}[eab]\st{(J)}{=}-[abe]-[bea]\st{(5.1a)}{=}-[az],\quad
  [az]\st{(5.1a)}{=}[aed]\st{(J)}{=}-[eda]-[dae]\st{(5.1a)}{=}[ey].$$

\noindent
Then $[ex]=-[az]=-[ey]$ and $[by]=[dx]-[az]=[dx]+[ex]=-[vx]$.
Finally, one has
 $[du]\st{(5.1a)}{=}[dev]\st{(J)}{=}-[vde]-[evd]\st{(5.1a)}{=}[vz]+[ey]$
 $\Ra$ $[vz]=[du]-[ey]=[dx]+[ex]=-[vx]$.
\end{proof}
\smallskip

\begin{claim}
(a) Every simple (resp., non-simple) commutator $[w]\in\bar L_4$
 of degree 3 can be expressed linearly via
 $[ax]$, $[bx]$, $[ay]$, $[dy]$,
 $[bz]$, $[ez]$, $[cu]$, $[eu]$
 (resp., via $[dx]$ and $[ex]$).
\smallskip

\noindent
(b) The degree 3 part of $\bar L_4$ is linearly generated by
 8 simple commutators $[ax]$, $[bx]$, $[ay]$, $[dy]$,
 $[bz]$, $[ez]$, $[cu]$, $[eu]$ and
 2 non-simple ones $[dx]$, $[ex]$.
\end{claim}
\begin{proof}
(a) The degree 3 part of $\bar L_4$ contains exactly
$$\begin{array}{ll}
\mbox{12 simple commutators} &
 [ax],\; [bx],\; [cx],\; [ay],\; [dy],\; [vy],\;
 [bz],\; [ez],\; [dz],\; [cu],\; [eu],\; [vu]\mbox{ and}\\
\mbox{12 non-simple ones} &
 [dx],\; [ex],\; [vx],\; [by],\; [cy],\; [ey],\;
 [az],\; [cz],\; [vz],\; [au],\; [bu],\; [du].
\end{array}$$
Due to the relations~(5.2a) one can throw out
 the four simple commutators $[cx]$, $[vy]$, $[dz]$, and $[vu]$.
By the relations~(5.2c) any non-simple commutator reduces
 to one of $[dx]$, $[ex]$, and $[vx]$.
To eliminate $[vx]$ it remains to apply
 the relation (5.2b): $[vx]=-[dx]-[ex]$.
\medskip

\noindent
(b) Any element of $\bar L_4$ is a sum of simple and non-simple
 commutators, i.e. $(a)\Ra(b)$.
\end{proof}

\begin{lemma}
(a) For every word $w$ in the alphabet $\LL=\{a,b,c,d,e,v\}$,
 containing at least two letters, and for any letters
 $s,s'\in\LL$, $[ss'w]=[s'sw]$ holds in the quotient $\bar L_4$.
\smallskip

\noindent
(b) Let $w$ be any word in the alphabet $\LL$, containing at
 least one letter $d$ or $e$.
Then in the quotient $\bar L_4$ the following relations hold:
$$\begin{array}{lll}
 [adx]=[edx],&
 [aex]=[e^2x],&
 [awx]=[ewx];\\
 {[bdx]}=[(-d-e)dx],&
 {[bex]}=[(-d-e)ex],&
 [bwx]=[(-d-e)wx];\\
 {[cdx]}=[d^2x],&
 [cex]=[dex],&
 [cwx]=[dwx].
\end{array}$$
\end{lemma}
\begin{proof}
(a) This is an analog of Claim~2.6a.
\smallskip

\noindent
(b) Apply the item~(a) and Claim~5.2 as follows
$$[adx]\st{(5.2c)}{=}
  [adu]\st{(5.4a)}{=}
  [dau]\st{(5.2c)}{=}
  [dex]\st{(5.4a)}{=}[edx],\quad
  [aex]\st{(5.4a)}{=}
  [eax]\st{(5.2a)}{=}
  -[ebx]-[ecx]=[e^2x];$$
$$[bdx]\st{(5.2c)}{=}
  [bdu]\st{(5.4a)}{=}
  [dbu]\st{(5.2c)}{=}
  [dvx]\st{(5.2b)}{=}
  [d(-d-e)x]\st{(5.4a)}{=}
  [(-d-e)dx];$$
$$[bex]\st{(5.2c)}{=}
  -[bey]\st{(5.4a)}{=}
  -[eby]\st{(5.2c)}{=}
  [evx]\st{(5.2b)}{=}
  [e(-d-e)x]\st{(5.4a)}{=}
  [(-d-e)ex];$$
$$[cdx]\st{(5.4a)}{=}
  [dcx]\st{(5.2a)}{=}
  -[dax]-[dbx]\st{(5.4b)}
  [d^2x],\quad
  [cex]\st{(5.2c)}{=}
  -[cey]\st{(5.4a)}{=}
  -[ecy]\st{(5.2c)}{=}[edx].$$
If a word $w$ in $\LL$ contains at least one letter $d$ or $e$,
 then by the item~(a) there is a word $w'$ in $\LL$ such that
 $[wx]=[w'dx]$ (without loss of generality) or $[wx]=[w'ex]$.
Then
$$[awx]=[aw'dx]\st{(5.4a)}{=}
  [w'adx]\st{(5.4b)}{=}
  [w'edx]\st{(5.4a)}{=}
  [ew'dx]=[ewx].$$
The proof in the other cases is similar.
\end{proof}

\begin{lemma}
(a) In each degree $n\geq 2$ every simple commutator $[w]\in\bar L_4$
 can be expressed linearly via $4(n-1)$ simple ones
 $[a^kb^lx]$, $[a^kd^ly]$, $[b^le^lz]$, $[c^ke^lu]$,
 $k+l=n-2$, $k,l\geq 0$.
\smallskip

\noindent
(b) In each degree $n\geq 3$ any non-simple commutator $[w]\in\bar L_4$
 can be expressed linearly via $n-1$ non-simple commutators
 $[d^ke^lx]$, where $k+l=n-2$, $k,l\geq 0$.
\smallskip

\noindent
(c) In any degree $n\geq 3$ the quotient $\bar L_4$ is linearly generated
 by $4(n-1)$ simple commutators $[a^kb^lab]$, $[a^kd^lad]$,
 $[b^ke^lbe]$, $[c^ke^lce]$ and $n-1$ non-simple ones $[d^ke^lab]$,
 $k+l=n-2$, $k,l\geq 0$.
\end{lemma}
\begin{proof}
(a) Let us consider a simple commutator $[w]$ not containing upper index 4.
Hence the commutator $[w]$ contains only the letters $a,b,c$.
By the relations~(5.2a) $[cx]=-[ax]-[bx]$,
 one can express $[w]$ via commutators $[a^kb^lx]$, $k,l\geq 0$.
The proof is analogous for another three types of simple
 commutators.
\smallskip

\noindent
(b) Let $w=w''s_1s_2s_3$ be the given word, where
 $s_1,s_2,s_3$ are the three last letters of $w$.
Since $[w]$ is non-simple, then by Lemma~5.4a one can permute
 the letters of $w$ in such a way that one may assume the commutator
 $[s_1s_2s_3]$ is non-simple.
By Claim~5.3b the commutator $[s_1s_2s_3]$ is expressed via $[dx]$
 and $[ex]$.
Then by Lemma~5.4a $[w]=[w''s_1s_2s_3]$ can be written in terms of
 $[dw''x]$ and $[ew''x]$.
Hence one can apply the induction on the length of $w$.
\smallskip

\noindent
The item (c) follows from (a) and (b).
\end{proof}


\subsection{Calculations in the quotient $\bar L_4$}

\begin{claim}
For all $k\geq 0$, $l\geq 1$, in the quotient $\bar L_4$
 the following relations hold:
$$\begin{array}{lll}
(a) &
[b^kd^l x]=[(-d-e)^kd^l x], &
[d^kb^l y]=-[d^k(-d-e)^l x]; \\
(b) &
[b^kc^l z]=[(-d-e)^kd^l x], &
[c^kb^l u]=[d^k(-d-e)^l x]; \\
(c) &
[d^ke^l y]=-[d^ke^l x], &
[e^kd^l u]=[e^kd^l x]; \\
(d) &
[a^kc^l y]=-[e^kd^l x], &
[c^ka^l u]=[d^ke^l x].
\end{array}$$
\end{claim}
\noindent
\emph{Proof.}
(a) By Lemma~5.4b one has $[bd^lx]=[(-d-e)d^lx]$, i.e.
 $[b^2d^lx]=[(-d-e)bd^lx]=[(-d-e)^2d^lx]$
 and so on, i.e.
 $[b^kd^lx]=[(-d-e)^kd^lx]$ holds for all $k\geq 0$, $l\geq 1$.
Similarly,
$$[by]\st{(5.2c)}{=}
  -[vx]\st{(5.2b)}{=}
  -[(-d-e)x]\Ra
  [b^2y]=-[b(-d-e)x]\st{(5.4b)}{=}-[(-d-e)^2x]\Ra$$
$$\Ra [b^ly]=-[(-d-e)^l x]\Ra
  [d^kb^ly]=-[d^k(-d-e)^l x]
  \mbox{ for all } k\geq 0,\;l\geq 1.$$

\noindent
The items (b), (c), and (d) are proved analogously to (a).
Apply the following formulae:
$$[cz]\st{(5.2c)}{=}-[dx]\Ra
  [c^2z]=-[cdx]\st{(5.4b)}{=}-[d^2x]\Ra
  [c^lz]=-[d^lx]\Ra$$
$$\Ra [b^kc^lz]=-[b^kd^lx]\st{(5.6a)}{=}-[(-d-e)^kd^lx]
  \mbox{ for all } k\geq 0,\;l\geq 1.$$

$$[bu]\st{(5.2c)}{=}[vx]\st{(5.2b)}{=}[(-d-e)x]\Ra
  [b^2u]=[b(-d-e)x]\st{(5.4b)}{=}[(-d-e)^2x]\Ra
  [b^lu]=[(-d-e)^lx]\Ra$$
$$\Ra [cb^lu]=[c(-d-e)^lx]\st{(5.4b)}{=}[d(-d-e)^lx]\Ra
  [c^kb^lu]=[d^k(-d-e)^lx]\mbox{ for all }k\geq 0,\; l\geq  1.$$

$$[ey]\st{(5.2c)}{=}-[ex]\Ra
  [e^2y]=-[e^2x]\Ra
  [e^ly]=-[e^lx]\Ra
  [d^ke^ly]=-[d^ke^lx],\; k\geq 0,\; l\geq 0.$$

$$[du]\st{(5.2c)}{=}[dx]\Ra
  [d^lu]=[d^lx]\Ra
  [e^kd^lu]=[e^kd^lx]\mbox{ for all }k\geq 0,\; l\geq  1.$$

$$[cy]\st{(5.2c)}{=}-[dx]\Ra
  [c^2y]=-[cdx]\st{(5.4b)}{=}-[d^2x]\Ra
  [c^ly]=-[d^lx]\Ra$$
$$\Ra [ac^ly]=-[ad^lx]\st{(5.4b)}{=}-[ed^lx]\Ra
  [a^kc^ly]=-[e^kd^lx]
  \mbox{ for all } k\geq 0,\;l\geq 1.$$

$$\mbox{Finally, one has}\quad
  [au]\st{(5.2c)}{=}[ex]\Ra
  [a^2u]=[aex]\st{(5.4b)}{=}[e^2x]\Ra
  [a^lu]=[e^lx]\Ra$$
$$\Ra [ca^lu]=[ce^lx]\st{(5.4b)}{=} [de^lx]\Ra
  [c^ka^lu]=[d^ke^lx].\eqno{\qed}$$

\begin{claim}
For any $k\geq 0$, in the quotient $\bar L_4$ one has
$$\begin{array}{lll}
(a) &
[(b+d)^k x]=[b^k x]-[(-d-e)^k x]+[(-e)^k x], &
[(b+d)^k y]=[d^k y]+[d^k x]-[(-e)^k x]; \\
(b) &
[(b+c)^k z]=[b^k z]+[(-d-e)^k x]-[(-e)^k x], &
[(b+c)^k u]=[c^k u]-[d^k x]+[(-e)^k x];\\
(c) &
[(d+e)^k y]=[d^k y]+[d^k x]-[(d+e)^k x], &
[(d+e)^k u]=[e^k u]-[e^k x]+[(d+e)^k x]; \\
(d) &
[(a+c)^k y]=[a^k y]+[e^k x]-[(d+e)^k x], &
[(a+c)^k u]=[c^k u]-[d^k x]+[(d+e)^k x].
\end{array}$$
\end{claim}
\noindent
\emph{Proof.}
(a) For each $k\geq 0$, one obtains
$$[(b+d)^k x]=
 [b^kx]+\left[\sum_{j=1}^k\binom{k}{j}b^{k-j}d^jx\right]
 \st{(5.6a)}{=}
 [b^kx]+\left[\sum_{j=1}^k\binom{k}{j}(-d-e)^{k-j}d^jx\right]=$$
$$=[b^kx]-[(-d-e)^kx]+\left[\sum_{j=0}^k\binom{k}{j}(-d-e)^{k-j}d^jx\right]
  =[b^kx]-[(-d-e)^kx]+[(d-d-e)^kx].$$

$$[(b+d)^k y]=
 [d^ky]+\left[\sum_{j=0}^k\binom{k-1}{j}b^{k-j}d^j y\right]
 \st{(5.6a)}{=}
 [d^ky]-\left[\sum_{j=0}^{k-1}\binom{k}{j}(-d-e)^{k-j}d^jx\right]=$$
$$=[d^ky]+[d^kx]-\left[\sum_{j=0}^k\binom{k}{j}(-d-e)^{k-j}d^jx\right]
  =[d^ky]+[d^kx]-[(-d-e+d)^kx].$$

The items (b),(c),(d) are proved analogously to (a).
The following formulae are used:
$$[(b+c)^k z]-[b^kz]=
 \sum_{j=1}^{k}\binom{k}{j}[b^{k-j}c^jz]
 \st{(5.6b)}{=}
 -\sum_{j=1}^k\binom{k}{j}[(-d-e)^{k-j}d^jx]=
 [(-d-e)^k x]-[(-e)^k x],$$

$$[(b+c)^k u]-[c^ku]=
 \sum_{j=0}^{k-1}\binom{k}{j}[b^{k-j}c^ju]
 \st{(5.6b)}{=}
 \sum_{j=0}^{k-1}\binom{k}{j}[(-d-e)^{k-j}d^jx]=
 -[d^k x]+[(-e)^k x],$$

$$[(d+e)^k y]-[d^ky]=
 \sum_{j=1}^{k}\binom{k}{j}[d^{k-j}e^j y]
 \st{(5.6c)}{=}
 -\sum_{j=1}^{k}\binom{k}{j}[d^{k-j}e^j x]=
 [d^k x]-[(d+e)^k x],$$

$$[(d+e)^k u]-[e^k u]=
 \sum_{j=1}^{k}\binom{k}{j}[e^{k-j}d^j u]
 \st{(5.6c)}{=}
 \sum_{j=1}^{k}\binom{k}{j}[e^{k-j}d^j x]=
 -[e^k x]+[(d+e)^k x],$$

$$[(a+c)^k y]-[a^ky]=
 \sum_{j=1}^{k}\binom{k}{j}[a^{k-j}c^j y]
 \st{(5.6d)}{=}
 -\sum_{j=1}^{k}\binom{k}{j}[e^{k-j}d^jx]=
 [e^k x]-[(d+e)^k x],$$

$$[(a+c)^k u]-[c^k u]=
 \sum_{j=1}^{k}\binom{k}{j}[c^{k-j}a^j y]
 \st{(5.6d)}{=}
 \sum_{j=1}^{k}\binom{k}{j}[d^{k-j}e^j x]=
 -[d^k x]+[(d+e)^k x].\eqno{\qed}$$

\begin{lemma}
(a) For all $k,l\geq 0$, in the quotient $\bar L_4$ one gets
$$\left\{\begin{array}{l}
[ a^k(b+d)^l a(b+d) ]
 = \bigl([a^kb^l x]+[a^kd^l y]\bigr) + [e^kd^l x]-[e^k(-d-e)^l x], \\
{[ (b+c)^k e^l(b+c)e ]}
 = \bigl([b^ke^l z]+[c^ke^l u]\bigr) - [d^ke^l x]+[(-d-e)^ke^l x].
\end{array} \right. \leqno{(5.8a)}$$

\noindent
(b) For all $k,l\geq 0$, in the quotient $\bar L_4$ one has
$$[(a+c)^k(d+e)^l(a+c)(d+e)]
 =\bigl([a^kd^l y]+[c^ke^l u]\bigr)+[e^kd^lx]-[d^ke^lx].\leqno{(5.8b)}$$
\end{lemma}
\begin{proof}
(a) One has $[a,b+d]\st{(5.1a)}{=}x+y$, $[b+c,e]\st{(5.1a)}{=}z+u$.
For all $k\geq 0$, $l\geq 1$, one gets
$$[ a^k(b+d)^l a(b+d) ]=[a^k(b+d)^l x]+[a^k(b+d)^l y]
 \st{(5.7a)}{=}
 \bigl([a^kb^l x]-[a^k(-d-e)^l x]+[a^k(-e)^l x]\bigr)+$$
$$+\bigl([a^kd^l y]+[a^kd^l x]-[a^k(-e)^l x]\bigr)
 \st{(5.4b)}{=} \bigl([a^kb^l x]+[a^kd^l y]\bigl)
 + [e^kd^l x]-[e^k(-d-e)^l x].$$
Observe that the above equations hold for $l=0$ also.
Similarly, one gets
$$[ (b+c)^k e^l(b+c)e ]\st{(5.4a)}{=}[e^l(b+c)^kz]+[e^l(b+c)^ku]
 \st{(5.7b)}{=}
 \bigl([e^lb^k z]+[e^l(-d-e)^k x]-[e^l(-e)^k x]\bigr)+$$
$$+\bigl([e^lc^k u]-[e^ld^k x]+[e^l(-e)^k x]\bigr)
 = \bigl([e^lb^k z]+[e^lc^k u]\bigr) - [e^ld^k x]+[e^l(-d-e)^k x].$$
\smallskip

\noindent
(b) Analogously to the item~(a), for all $k,l\geq 0$, one obtains
$$[(a+c)^k(d+e)^l(a+c)(d+e)]\st{(5.1a)}{=}
  [(a+c)^k(d+e)^ly]+[(a+c)^k(d+e)^lu]\st{(5.7c)}{=}$$
$$\st{(5.7c)}{=}
  [(a+c)^kd^ly]+[(a+c)^kd^lx]-[(a+c)^k(d+e)^lx]+
  [(a+c)^ke^lu]-[(a+c)^ke^lx]+[(a+c)^k(d+e)^lx]$$
$$\st{(5.2a),(5.4a)}{=}
  [d^l(a+c)^ky]+[d^l(-b)^kx]+
  [e^l(a+c)^ku]-[e^l(-b)^kx]\st{(5.7d),(5.4a)}{=}$$
$$\bigl([d^la^ky]+[d^le^kx]-[d^l(d+e)^kx]\bigr)+
  \bigl([e^lc^ku]-[e^ld^kx]+[e^l(d+e)^kx]\bigr)+[(-b)^k(d^l-e^l)x]$$
$$\st{(5.4a),(5.4b)}{=}
 \bigl([a^kd^l y]+[c^ke^l u]\bigr)+[e^kd^lx]-[d^ke^lx]+[(e^l-d^l)(d+e)^kx]
 +[(d+e)^k(d^l-e^l)x]$$
$$\Bigl|\; [(-b)^k(d^l-e^l)x]\st{(5.4b)}{=}[(d+e)^k(d^l-e^l)x]
 \mbox{ was used }\Bigr|
 =\bigl([a^kd^l y]+[c^ke^l u]\bigr)+[e^kd^lx]-[d^ke^lx].$$
Note that the relation
 $[(-b)^k(d^l-e^l)x]\st{(5.4b)}{=}[(d+e)^k(d^l-e^l)x]$
 holds for any $l\geq 0$.
\end{proof}


\subsection{Checking the compressed pentagon~$\ol{(1.3c)}$}

\begin{lemma}
For any compressed associator $\bar\ph\in\wh{\bar L}_3$,
 the compressed pentagon~$\ol{(1.3c)}$ is equivalent to
 the following equation in the algebra $\wh{\bar L}_4$:
$$\bar\ph(b,e)+\bar\ph(a+c,d+e)+\bar\ph(a,b)
 =\bar\ph(a,b+d)+\bar\ph(b+c,e).\leqno{(5.9)}$$
\end{lemma}
\begin{proof}
Let us rewrite explicitly the pentagon~(1.3c)
 for a compressed associator $\bar\ph\in\hat{\bar L}_3$
$$\exp\bigl(\bar\ph(b,e)\bigr)\cdot
  \exp\bigl(\bar\ph(a+c,d+e)\bigr)\cdot
  \exp\bigl(\bar\ph(a,b)\bigr)
 =\exp\bigl(\bar\ph(a,b+d)\bigr)\cdot
  \exp\bigl(\bar\ph(b+c,e)\bigr).$$
Since in the quotient $\bar L_4$ all commutators commute, then
 the series $\bar\ph(b,e)$, $\bar\ph(a+c,d+e)$, $\bar\ph(a,b)$,
 $\bar\ph(a,b+d)$, and $\bar\ph(b+c,e)$ commute with each other
 in $\wh{\bar L}_4$.
Hence, taking the logarithm of both sides of the above pentagon,
 one needs to apply the simplest case of CBH formula (2.3):
 $\log\bigl(\exp(P)\cdot\exp(Q)\bigr)=P+Q$ provided that $P,Q$ commute.
\end{proof}

\begin{proposition}
Let $f(\la,\mu)=\sum\limits_{k,l\geq 0}\al_{kl}\la^k\mu^l$ be
 the generating function of the coefficients $\al_{kl}=\al_{lk}$ of
 a compressed Drinfeld associator $\bar\ph\in\wh{\bar L}_3$.
Then in the algebra $\wh{\bar L}_4$ the compressed
 pentagon equation~$\ol{(1.3c)}=(5.9)$
 follows from the symmetry $\al_{kl}=\al_{lk}$.
\end{proposition}
\emph{Proof.}
By Lemma~5.8b the left hand side of (5.9) is
$$\sum_{k,l\geq 0}\al_{kl}\Bigl(
 [b^ke^lbe]+[(a+c)^k(d+e)^l(a+c)(d+e)]+[a^kb^lab]
 \Bigr)\st{(5.8b)}{=}$$
$$\st{(5.8b)}{=}\sum_{k,l\geq 0}\al_{kl}\Bigl(
 [b^ke^lz]
 +\bigl([a^kd^l y]+[c^ke^l u]+[e^kd^lx]-[d^ke^lx]\bigr)
 +[a^kb^lx] \Bigr).$$
Similarly by Lemma~5.8a the right hand side of (5.9) is
$$\sum_{k,l\geq 0}\al_{kl}\Bigl(
 [a^k(b+d)^ka(b+d)]+[(b+c)^ke(b+c)e] \Bigr)\st{(5.8a)}{=}$$
$$\sum_{k,l\geq 0}\al_{kl}\Bigl(
 \bigl([a^kb^l x]+[a^kd^l y] + [e^kd^l x]-[e^k(-d-e)^l x] \bigr)
+\bigl([b^ke^l z]+[c^ke^l u] - [d^ke^l x]+[(-d-e)^ke^l x] \bigr)\Bigr).$$
$$\mbox{The difference is }
\sum_{k,l\geq 0}\al_{kl}\Bigl([e^k(-d-e)^l x]-[(-d-e)^ke^l x]\Bigr)
 \st{(5.4a)}{=}0 \mbox{ if }\al_{kl}=\al_{lk}.\eqno{\qed}$$


\section{Drinfeld series, zeta values, and problems}

In this section one shall check that the Drinfeld series from
 Definition~6.1b (a compressed associator expressed via zeta
 values) is contained in the general family~(1.5c).


\subsection{Riemann zeta-function of even integers}

\begin{definition}[Riemann zeta function $\ze(n)$,
 the Drinfeld series $s(\la,\mu)\;$]
\noindent
\smallskip

\noindent
(a) Let $\ze(n)=\sum\limits_{k=1}^{\infty}\dfrac{1}{k^n}$ be
 \emph{the classical Riemann zeta function}.
Put $\te_n:=\dfrac{\ze(n)}{n(\pi\sqrt{-1})^n}$.
\smallskip

\noindent
(b) Introduce
 $S(\la)=\sum\limits_{n=2}^{\infty}\te_n\la^n
 =\sum\limits_{n=2}^{\infty}
  \dfrac{\ze(n)\la^n}{n(\pi\sqrt{-1})^n}$.
\emph{The Drinfeld series} is\\
$s(\la,\mu)=S(\la)+S(\mu)-S(\la+\mu)
 =\sum\limits_{n=2}^{\infty}
 \dfrac{\ze(n)}{n}\cdot\dfrac{\la^n+\mu^n-(\la+\mu)^n}{(\pi\sqrt{-1})^n}$.
\ed
\end{definition}
\smallskip

Due to the change $t^{ij}\mapsto 2t^{ij}$ in the hexagon~(1.3b),
 one has missed $2^n$ in the denominators of $s(\la,\mu)$.
The following theorem is quoted from
 \cite[Chapter XIX, remark 6.6(b), p.~468]{Kas}.
\smallskip

\begin{theorem} \cite{Dr2} 
There is a compressed Drinfeld associator defined by
$$\ti f^D(\la,\mu)=1+\la\mu f^D(\la,\mu)
 =\exp\bigl(s(\la,\mu)\bigr).\leqno{(6.2)}$$
\end{theorem}
\smallskip

The Drinfeld series $s(\la,\mu)$ leads to the well-known formula
 for even zeta values.
\smallskip

\begin{lemma}
For each $n\geq 1$, one has
 $2n\te_{2n}=-\dfrac{2^{2n}B_{2n}}{2(2n)!}$ and
 $\ze(2n)=(-1)^{n-1}\dfrac{(2\pi)^{2n}}{2(2n)!}B_{2n}$.
\end{lemma}
\begin{proof}
Theorem~1.5c will be used for $\mu=-\la$.
Lemma~4.2 says that
 $f(\la,-\la)=\dfrac{1}{\la^2}-\dfrac{2}{\la(e^{\la}-e^{-\la})}$.
By Definition~6.1b one has
 $s(\la,-\la)=2\sum\limits_{n=1}^{\infty}\te_{2n}\la^{2n}$.
Let us substitute $\mu=-\la$ in the formula (6.2).
Then by using the formula~(1.5c) one gets
$$\log\bigl(\ti f(\la,-\la)\bigr)=\log(1-\la^2f(\la,-\la))
 =\log\left(\frac{2\la}{e^{\la}-e^{-\la}}\right)
 =2\sum\limits_{n=1}^{\infty}\te_{2n}\la^{2n}.$$
Taking the first derivative of the above equation by $\la$,
 one obtains
$$\left(
 \dfrac{e^{\la}-e^{-\la}-\la(e^{\la}+e^{-\la})}{(e^{\la}-e^{-\la})^2}
 \right):\left(\frac{\la}{e^{\la}-e^{-\la}}\right)
 =\dfrac{1}{\la}-1-\dfrac{2}{e^{2\la}-1}
 =2\sum\limits_{n=1}^{\infty}(2n\te_{2n})\la^{2n-1}.$$
Multiply the resulting equation by $\la$ and use
 the definition of the Bernoulli numbers:
$$1-\la-\dfrac{2\la}{e^{2\la}-1}
 =-\sum_{n=1}^{\infty}\dfrac{B_{2n}}{(2n)!}(2\la)^{2n}
 =2\sum\limits_{n=1}^{\infty}(2n\te_{2n})\la^{2n},\mbox{ hence }
 2n\te_{2n}=-\dfrac{B_{2n}}{2(2n)!}2^{2n}.$$
It remains to use the formula from Definition~6.1a:
 $\te_{2n}=\dfrac{\ze(2n)}{2n(\pi\sqrt{-1})^{2n}}
  =(-1)^n\dfrac{\ze(2n)}{2n\pi^{2n}}$.
\end{proof}
\smallskip

\begin{example}
By Lemma~6.3 and Table~A.2 of Appendix one can easily calculate
$$\begin{array}{lllll}
\te_2=-\dfrac{1}{12},\quad &
\te_4=\dfrac{1}{360},\quad &
\te_6=-\dfrac{1}{5670},\quad &
\te_8=\dfrac{1}{75600},\quad &
\te_{10}=-\dfrac{1}{935550};\\ \\
\ze(2)=\dfrac{\pi^2}{6},\quad &
\ze(4)=\dfrac{\pi^4}{90},\quad &
\ze(6)=\dfrac{\pi^6}{945},\quad &
\ze(8)=\dfrac{\pi^8}{9450},\quad &
\ze(10)=\dfrac{\pi^{10}}{93555}.\quad\quad\quad\bt
\end{array}$$
\end{example}
\smallskip


\subsection{Substitution of the Drinfeld series $s(\la,\mu)$}

\begin{claim}
The Drinfeld series $s(\la,\mu)$ satisfies the following equations:
$$\dfrac{e^{\la+\mu}-e^{-\la-\mu}}{2(\la+\mu)} \left(
 \dfrac{2\om}{e^{\om}-e^{-\om}}
 +\sum_{n=3}^{\infty} h_n(\la,\mu)\right)
 =Even\Bigl(\exp\bigl(s(\la,\mu)\bigr)\Bigr),\quad\mbox{ and}\leqno{(6.5a)}$$
$$\la\mu\cdot\dfrac{e^{\la+\mu}-e^{-\la-\mu}}{2} \left(
 \sum_{n=0}^{\infty}\ti\be_{n0}\om^{2n}
 +\sum_{n=3}^{\infty} \ti h_n(\la,\mu)\right)
 =Odd\Bigl(\exp\bigl(s(\la,\mu)\bigr)\Bigr),\leqno{(6.5b)}$$
$$\mbox{where }h_n(\la,\mu)
 =\sum_{k=1}^{[\frac{n}{3}]}\be_{nk}\la^{2k}\mu^{2k}(\la+\mu)^{2k}
  \om^{2n-6k}\mbox{ for }n\geq 3,\quad
 \om=\sqrt{\la^2+\la\mu+\mu^2}.$$
The polynomials $\ti h_n(\la,\mu)$ are defined by the same
 formula as $h_n(\la,\mu)$, except the coefficients
 $\ti\be_{nk}\in\C$ are substituted for $\be_{nk}\in\C$.
\end{claim}
\noindent
\emph{Proof} follows from Theorems~1.5c, 6.2, and
 $Odd\bigl(1+\la\mu f^D(\la,\mu)\bigr)=
  \la\mu Odd\bigl(f^D(\la,\mu)\bigr)$.
\qed
\smallskip

\begin{claim}
Any series $s(\la,\mu)$ satisfies the following relations:
$$Even\Bigl(\exp\bigl(s(\la,\mu)\bigr)\Bigr)
 =\exp\Bigl(Even\bigl(s(\la,\mu)\bigr)\Bigr)\cdot
 Even\Bigl(\exp\bigl( Odd(s(\la,\mu)) \bigr)\Bigr),\leqno{(6.6a)}$$
$$Odd\Bigl(\exp\bigl(s(\la,\mu)\bigr)\Bigr)
 =\exp\Bigl(Even\bigl(s(\la,\mu)\bigr)\Bigr)\cdot
 Odd\Bigl(\exp\bigl( Odd(s(\la,\mu)) \bigr)\Bigr).\leqno{(6.6b)}$$
\end{claim}
\begin{proof}
Actually, one has
 $Even\Bigl(\exp\bigl(s(\la,\mu)\bigr)\Bigr)
 =\dfrac{ \exp\bigl(s(\la,\mu)\bigr)+\exp\bigl(s(-\la,-\mu)\bigr) }{2}=$
$$=\dfrac{ \exp\Bigl( Even\bigl(s(\la,\mu))+Odd\bigl(s(\la,\mu)\bigr) \Bigr)
 +\exp\Bigl( Even\bigl(s(\la,\mu)\bigr)-Odd\bigl(s(\la,\mu)\bigr) \Bigr)}{2}=$$
$$=\exp\Bigl( Even\bigl(s(\la,\mu)\bigr) \Bigr)\cdot
 \dfrac{ \exp\Bigl(Odd\bigl(s(\la,\mu)\Bigr)
  +\exp\Bigl( -Odd\bigl(s(\la,\mu)\bigr) \Bigr)}{2} \Lra (6.6a).$$
The formula (6.6b) is proved absolutely analogously.
\end{proof}

\begin{claim}
(a) For the series $S(\rho)$, one has
 $\exp\Bigl(-2Even\bigl(S(\rho)\bigr)\Bigr)
 =\dfrac{e^{\rho}-e^{-\rho}}{2\rho}$.
\smallskip

\noindent
(b) For the Drinfeld series $s(\la,\mu)$, one has
$$\exp\Bigl( Even\bigl(s(\la,\mu)\bigr)\Bigr)
 =\sqrt{ \dfrac{e^{\la+\mu}-e^{-\la-\mu}}{2(\la+\mu)}\cdot
  \dfrac{2\la}{e^{\la}-e^{-\la}}\cdot
  \dfrac{2\mu}{e^{\mu}-e^{-\mu}} }.\leqno{(6.7)}$$
\end{claim}
\noindent
\emph{Proof.}
(a) By Definition~6.1b, one obtains
 $Even\bigl(S(\rho)\bigr)=\sum\limits_{n=1}^{\infty}\te_{2n}\rho^{2n}$.
So, one needs to get
 $-2\sum\limits_{n=1}^{\infty}\te_{2n}\rho^{2n}
 =\log\left( \dfrac{e^{\rho}-e^{-\rho}}{2\rho} \right)$.
This formula holds for $\rho=0$.
Then it suffices to check that the first derivatives of both
 hand sides are equal.
For the left hand side, apply the formula for $2n\te_{2n}$
 from Lemma~6.3:
$$-2\sum\limits_{n=1}^{\infty}2n\te_{2n}\rho^{2n-1}
 =\dfrac{1}{\rho}\sum\limits_{n=1}^{\infty}\dfrac{B_{2n}}{(2n)!}(2\rho)^{2n}
 \st{(2.1b)}{=}
 \dfrac{1}{\rho}\left(\dfrac{2\rho}{e^{2\rho}-1}-1+\rho\right)
 =1-\dfrac{1}{\rho}+\dfrac{2}{e^{2\rho}-1}.$$
For the right hand side, one gets
$$\left( \dfrac{e^{\rho}-e^{-\rho}}{2\rho} \right)'
 :\left( \dfrac{e^{\rho}-e^{-\rho}}{2\rho} \right)
 =\dfrac{(e^{\rho}+e^{-\rho})\rho-(e^{\rho}-e^{-\rho})}{\rho^2}
  \cdot \dfrac{\rho}{e^{\rho}-e^{-\rho}}
 =1-\dfrac{1}{\rho}+\dfrac{2}{e^{2\rho}-1}\mbox{ as required}.$$
\smallskip

\noindent
(b) Apply the item~(a) and Definition~6.1b:
$$\exp\Bigl(Even\bigl(s(\la,\mu)\bigr)\Bigr)
 \st{(6.1b)}{=}
  \exp\Bigl(Even\bigl(S(\la)\bigr)\Bigr)\cdot
  \exp\Bigl(Even\bigl(S(\mu)\bigr)\Bigr)\cdot
  \exp\Bigl(-Even\bigl(S(\la+\mu)\bigr)\Bigr)
 \st{(6.7a)}{=}$$
$$\st{(6.7a)}{=}
  \sqrt{ \dfrac{2\la}{e^{\la}-e^{-\la}} }
  \sqrt{ \dfrac{2\mu}{e^{\mu}-e^{-\mu}} }
  \sqrt{ \dfrac{e^{\la+\mu}-e^{-\la-\mu}}{2(\la+\mu)} }
 =\sqrt{ \dfrac{e^{\la+\mu}-e^{-\la-\mu}}{2(\la+\mu)}\cdot
  \dfrac{2\la}{e^{\la}-e^{-\la}}\cdot
  \dfrac{2\mu}{e^{\mu}-e^{-\mu}} }\;.\eqno{\qed}$$
\smallskip

\begin{claim}
The equations (6.5a) and (6.5b) are equivalent to
 (6.8a) and (6.8b), respectively
$$\sum_{k=0}^{\infty}\dfrac{\bigl(\te(\la,\mu)\bigr)^{2k}}{(2k)!}
 =h(\la,\mu)\sqrt{ \dfrac{e^{\la+\mu}-e^{-\la-\mu}}{2(\la+\mu)}\cdot
  \dfrac{e^{\la}-e^{-\la}}{2\la}\cdot
  \dfrac{e^{\mu}-e^{-\mu}}{2\mu} }\;,\leqno{(6.8a)}$$
$$\sum_{k=0}^{\infty}\dfrac{\bigl(\te(\la,\mu)\bigr)^{2k+1}}{(2k+1)!}
 =\ti h(\la,\mu)\la\mu(\la+\mu)
  \sqrt{ \dfrac{e^{\la+\mu}-e^{-\la-\mu}}{2(\la+\mu)}\cdot
  \dfrac{e^{\la}-e^{-\la}}{2\la}\cdot
  \dfrac{e^{\mu}-e^{-\mu}}{2\mu} }\;,\leqno{(6.8b)}$$
$$\mbox{where }
\te(\la,\mu)=-\sum_{n=1}^{\infty}\te_{2n+1}F^D_{2n+1}(\la,\mu)
 =\sum_{n=1}^{\infty} \dfrac{(-1)^n\ze(2n+1)}{(2n+1)\pi^{2n+1}\sqrt{-1}}
 \Bigl(\la^{2n+1}+\mu^{2n+1}-(\la+\mu)^{2n+1}\Bigr),$$
$$h(\la,\mu)=\dfrac{2\om}{e^{\om}-e^{-\om}}
 +\sum_{n=3}^{\infty} h_n(\la,\mu),\quad
\ti h(\la,\mu)=\sum_{n=1}^{\infty}\ti\be_{n0}\om^{2n}
 +\sum_{n=3}^{\infty} \ti h_n(\la,\mu).$$
The polynomials $h_n(\la,\mu)$ and $\ti h_n(\la,\mu)$ are the
 same as in Claim~6.5.
\end{claim}
\noindent
\emph{Proof} follows from Claims 6.6--6.7, the formula
 $Odd\bigl(s(\la,\mu)\bigr)=\te(\la,\mu)$, and
$$Even\Bigl(\exp\bigl( Odd(s(\la,\mu))\bigr)\Bigr)
 =\sum_{k=0}^{\infty}\dfrac{\te^{2k}(\la,\mu)}{(2k)!},\;
  Odd\Bigl(\exp\bigl( Odd(s(\la,\mu))\bigr)\Bigr)
 =\sum_{k=0}^{\infty}\dfrac{\te^{2k+1}(\la,\mu)}{(2k+1)!}.
 \eqno{\qed}$$

\begin{proposition}
There exist parameters $\be_{nk},\ti\be_{nk}\in\C$, expressed
 linearly via monomials of odd zeta values, such that
 the equations (6.8a) and (6.8b) hold indentically.
Hence the Drinfeld series $s(\la,\mu)$ does not
 lead to polynomial relations between odd zeta values.
\end{proposition}
\begin{proof}
The left hand side of (6.8a) is unchanged under the
 transformation $\mu\lra(-\la-\mu)$.
By Lemma~4.9 such a series is a sum of associator polynomials
 of even degrees.
In other words, there exist parameters $\be'_{nk}\in\C$,
 expressed in terms of $\te_{2n+1}$, such that
$$\mbox{the left hand side of (6.8a)}
 =1+\sum_{n=1}^{\infty}\sum_{k=1}^{[\frac{n}{3}]}
 \be'_{nk}\la^{2k}\mu^{2k}(\la+\mu)^{2k}(\la^2+\la\mu+\mu^2)^{n-3k}.$$
Actually, the left hand side of (6.8a) minus 1 is divided by
 $\la^2\mu^2(\la+\mu)^2$, hence $\be'_{n0}=0$ for $n\geq 1$.
The right hand side of (6.8a) has the same form:
 it equals to 1 when $\la+\mu=0$.
Then one can find the parameters $\be_{nk}$ via $\be'_{nk}$
 (hence via odd zeta values)
 in such a way that the equation (6.8a) holds indentically.
The proof for the equation (6.8b) is similar.
\end{proof}

So, Proposition~6.9 supports the long standing conjecture in the
 number theory: odd zeta values are algebraically independent
 over the rationals \cite{Car2}.

Proposition~6.9 will be reproved explicitly up to degree~7 in
 Claim~A.6 of Appendix.
\medskip

\noindent
{\bf Proof of Corollary~1.6c} follows from Proposition~6.9
 since all monomials consisting of odd zeta values
 can be considered as free parameters.
The third associator~(1.6c) is obtained from the Drinfeld series
 $f^D(\la,\mu)$ by taking the free parameters
 $\zeta(2k+1)=0$ for each $k\geq 1$.
\qed


\subsection{Conjectures and open problems}

Theorem~1.5 describes only compressed associators.
This is a first step in the general problem
 to find a complete rational associator.
Theorem~1.5c gives a hope to describe Drinfeld associators
 up to triple commutators.

\begin{problem}
(a) Is it true that any compressed Drinfeld associator
 is the projection under
 $\hat L_3\to \bigl[[\hat L_3,\hat L_3],[\hat L_3,\hat L_3]\bigr]$
 of the logarithm $\ph(a,b)$ of an honest one from Definition~1.3b?
\smallskip

\noindent
(b) Describe all Drinfeld associators up to triple commutators.
In other words, solve the hexagon and pentagon in the quotient
 $L_3/[L_3',[L_3',L_3']]$, where $L_3'=[L_3,L_3]$.
\end{problem}

A compressed associator $\bar\ph(a,b)$ already contains
 a lot of information.
One may try to pass through the LM-BN construction \cite{LM3, BN2} to get
 a well-defined invariant of knots in a quotient of the algebra
 $A$ of chord diagrams.
Let $\De\ph$ be the difference between the exact logarithm of a
 Drinfeld associator $\ph(a,b)\in\hat L_3$ and
 its compressed image $\bar\ph(a,b)\in\wh{\bar L}_3$.
Then
$$\De\Phi=\exp(\bar\ph+\De\ph)-\exp(\bar\ph)
 =\De\ph+\dfrac{1}{2}\Bigl((\De\ph)^2+\De\ph\bar\ph+\bar\ph\De\ph\Bigr)
 +\cdots$$
 is the error of a Drinfeld associator in the completion $\hat A_3$.
By the definition of the quotient $\bar L_3$ the error $\De\ph$
 is expressed via differences $[w_{kl}ab]-[a^kb^lab]$,
 where $w_{kl}$ is a word containing exactly $k$ letters $a$
 and $l$ letters $b$.

One wants to factorize the algebra $A$ in
 such a way that the error $\De\Phi$ dissappears
 and the LM-BN construction leads to a knot invariant in this quotient.
Since all terms of $\De\Phi$ contain $\De\ph$ as a factor, then
 it suffices to kill chord diagrams containing $[w_{kl}ab]-[a^kb^lab]$.
Roughly, the LM-BN construction maps chord diagrams
 on $n$ vertical strands onto chord diagrams on the circle.
Let $\Si$ be a sum of commutators containing the symbols $a,b,c$.
Define \emph{the closure} $\wh{\Si}=0$
 as the relation $\Si=0$ considered in the algebra $A$ of chord diagrams.
Formally, one draws the relation $\Si=0$ on 3 vertical strands
 and assume that these strands are three arcs of a circle.
For instance, the 4T relation from Definition~1.1b is
 the closure of $\Si=[a,b]-[b,c]$.
Since $[abab]=[baab]$ holds in $L_3$, then at the stage of chord diagrams
 the first non-trivial relation is the closure
 of $\bigl[[ab],[aab]\bigr]=[ab]\cdot[aab]-[aab]\cdot[ab]$.

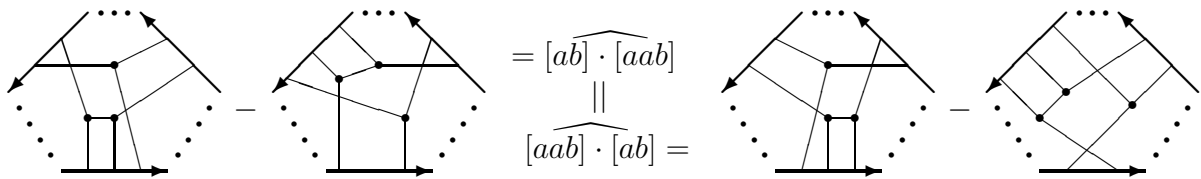
\begin{picture}(450,70)(0,0)

{\thicklines
\put(20,0){\vector(1,0){40}}
\put(30,60){\vector(-1,-1){30}}
\put(80,30){\vector(-1,1){30}}
\put(85,20){$-$}
\put(120,0){\vector(1,0){40}}
\put(130,60){\vector(-1,-1){30}}
\put(180,30){\vector(-1,1){30}}
\put(190,40){$=\wh{[ab]\cdot[aab]}$}
\put(220,25){$||$}
\put(195,5){$\wh{[aab]\cdot[ab]}=$}
\put(290,0){\vector(1,0){40}}
\put(300,60){\vector(-1,-1){30}}
\put(350,30){\vector(-1,1){30}}
\put(355,20){$-$}
\put(390,0){\vector(1,0){40}}
\put(400,60){\vector(-1,-1){30}}
\put(450,30){\vector(-1,1){30}}
}

\multiput(35,60)(5,0){3}{\circle*{2}}
\multiput(16,6)(-4,6){4}{\circle*{2}}
\multiput(64,6)(4,6){4}{\circle*{2}}

\put(30,20){\circle*{3}}
\put(40,20){\circle*{3}}
\put(40,40){\circle*{3}}
\put(30,20){\line(0,-1){20}}
\put(30,20){\line(-1,3){10}}
\put(30,20){\line(1,0){10}}
\put(40,20){\line(0,-1){20}}
\put(40,20){\line(3,2){30}}
\put(40,40){\line(1,-4){10}}
\put(40,40){\line(-1,0){30}}
\put(40,40){\line(2,1){20}}

\multiput(135,60)(5,0){3}{\circle*{2}}
\multiput(116,6)(-4,6){4}{\circle*{2}}
\multiput(164,6)(4,6){4}{\circle*{2}}

\put(125,35){\circle*{3}}
\put(140,40){\circle*{3}}
\put(150,20){\circle*{3}}
\put(125,35){\line(0,-1){35}}
\put(125,35){\line(3,1){15}}
\put(125,35){\line(-1,1){10}}
\put(140,40){\line(-1,1){15}}
\put(140,40){\line(1,0){30}}
\put(150,20){\line(0,-1){20}}
\put(150,20){\line(-3,1){45}}
\put(150,20){\line(1,3){10}}

\multiput(305,60)(5,0){3}{\circle*{2}}
\multiput(286,6)(-4,6){4}{\circle*{2}}
\multiput(334,6)(4,6){4}{\circle*{2}}

\put(310,20){\circle*{3}}
\put(310,40){\circle*{3}}
\put(320,20){\circle*{3}}
\put(310,20){\line(0,-1){20}}
\put(310,20){\line(1,0){10}}
\put(310,20){\line(-3,2){30}}
\put(320,20){\line(0,-1){20}}
\put(320,20){\line(1,3){10}}
\put(310,40){\line(-1,-4){10}}
\put(310,40){\line(-2,1){20}}
\put(310,40){\line(1,0){30}}

\multiput(405,60)(5,0){3}{\circle*{2}}
\multiput(386,6)(-4,6){4}{\circle*{2}}
\multiput(434,6)(4,6){4}{\circle*{2}}

\put(390,20){\circle*{3}}
\put(400,30){\circle*{3}}
\put(425,25){\circle*{3}}
\put(390,20){\line(3,-2){30}}
\put(390,20){\line(-1,1){15}}
\put(390,20){\line(1,1){10}}
\put(400,30){\line(-1,1){15}}
\put(400,30){\line(3,2){30}}
\put(400,0){\line(1,1){40}}
\put(395,55){\line(1,-1){30}}

\end{picture}
\vspace{2mm}

In the above figure the relation was drawn briefly by using
 STU relations from \cite{BN1}.

After Vassilev's paper \cite{Vas} one usually uses another definition
 of Vassiliev invariants via chord diagrams \cite{BN1}.
Roughly, a Vassiliev invariant of framed knots is the composition of
 the Kontsevich integral and \emph{a weight function} on $A$, i.e.
 a linear function on chord diagrams, satisfying the 4T relations,
 see Definition~1.1b.

\begin{definition}[compressed algebra $\bar A$ of chord diagrams,
 compressed Vassiliev invariants]

\noindent
(a) For a word $w_{kl}$ containing exactly $k$ letters $a$ and
 exactly $l$ letters $b$, put
 $\Si(w_{kl}):=[w_{kl}ab]-[a^kb^lab]$, $k,l\geq 1$.
Let $\bar A$ be the quotient of the classical algebra
 $A$ of chord diagrams on the circle by the ideal generated
 by the relations $\wh{\Si(w_{kl})}=0$ for all $k,l\geq 1$.
\smallskip

\noindent
(b) \emph{A compressed weight function} is a linear function
 on the compressed algebra $\bar A$.
In other words, a Vassiliev invariant is \emph{compressed},
 if the corresponding weight function satisfies
 $\wh{\Si(w_{kl})}=0$ for all $k,l\geq 1$.
\emph{The compressed Kontsevich integral} $\bar Z_K$ of a knot $K$
 is the image of the classical Kontsevich integral $Z_K$
 under the natural projection $A\to\bar A$.
\end{definition}

Vassiliev invariants of degrees 2, 3, 4 are compressed
 ones, i.e. the theory is not empty.

\begin{problem}
(a) Check carefully that the LM-BN construction \cite{LM3, BN2},
 for a compressed Drinfeld associator, gives rise to a well-defined
 knot invariant in the compressed algebra $\bar A$.
Does the resulting invariant depend on a particularly chosen
 compressed associator?
\smallskip

\noindent
(b) Which quantum invariants are compressed Vassiliev invariants?
\smallskip

\noindent
(c) Describe all compressed Vassiliev invariants (as linear
 functions on the algebra $\bar A$).
\smallskip

\noindent
(d) Compute the compressed Kontsevich integral for non-trivial
 knots, e.g. torus knots.
\smallskip

\noindent
(e) Which knots can be classified via compressed Vassiliev invariants?
\end{problem}


\section*{Appendix: explicit formulae}

\noindent
{\bf Claim A.1}.
For each $n\geq 1$, the following formulae hold:
$$\begin{array}{ll}
C_{2,2n}=C_{2n,2}=B_{2n},&
C_{4,2n+1}=-C_{2n+1,4}=4B_{2n+4}+4B_{2n+2},\\
C_{2,2n+1}=-C_{2n+1,2}=2B_{2n+2}, &
C_{5,2n}=-C_{2n,5}=5B_{2n+4}+10B_{2n+2}+B_{2n},\\
C_{3,2n}=-C_{2n,3}=3B_{2n+2}+B_{2n},&
C_{5,2n+1}=C_{2n+1,5}=10B_{2n+4}+5B_{2n+2},\\
C_{3,2n+1}=C_{2n+1,3}=3B_{2n+2},&
C_{6,2n}=C_{2n,6}=15B_{2n+4}+15B_{2n+2}+B_{2n},\\
C_{4,2n}=C_{2n,4}=6B_{2n+2}+B_{2n},&
C_{6,2n+1}=-C_{2n+1,6}=6B_{2n+6}+20B_{2n+4}+6B_{2n+2}.
\end{array}$$
\emph{Proof} follows from Lemma~2.11 and
 \emph{the first key point}: $B_{2n+1}=0$ for each $n\geq 1$.
\qed
\medskip

By Claim A.1 one can easily compute the numbers $C_{mn}$ for
 $m+n\leq 12$.
\smallskip

{\bf Table A.2} of the extended Bernoulli numbers $C_{mn}$.
\medskip

{\renewcommand{\arraystretch}{2.2}
\begin{tabular}{c|c|c|c|c|c|c|c|c|c|c|c}
$m \backslash n$ & 1 & 2 & 3 & 4 & 5 & 6 & 7 & 8 & 9 & 10 & 11\\
\hline

1 & $-\dfrac{1}{2}$ & $\dfrac{1}{6}$ & $0$ & $-\dfrac{1}{30}$ & $0$
  & $\dfrac{1}{42}$ & $0$ & $-\dfrac{1}{30}$ & $0$
  & $\dfrac{5}{66}$ & $0$\\
\hline

2 & $-\dfrac{1}{6}$  & $\dfrac{1}{6}$ & $-\dfrac{1}{15}$
  & $-\dfrac{1}{30}$  & $\dfrac{1}{21}$ & $\dfrac{1}{42}$
  & $-\dfrac{1}{15}$ & $-\dfrac{1}{30}$
  & $\dfrac{5}{33}$  & $\dfrac{5}{66}$\\
\hline
3 & $0$ & $\dfrac{1}{15}$ & $-\dfrac{1}{10}$ & $\dfrac{4}{105}$
  & $\dfrac{1}{14}$ & $-\dfrac{8}{105}$ & $-\dfrac{1}{10}$
  & $\dfrac{32}{165}$ & $\dfrac{5}{22}$\\
\hline
4 & $\dfrac{1}{30}$ & $-\dfrac{1}{30}$ & $-\dfrac{4}{105}$
  & $\dfrac{23}{210}$ & $-\dfrac{4}{105}$ & $-\dfrac{37}{210}$
  & $\dfrac{28}{165}$ & $\dfrac{139}{330}$\\
\hline
5 & $0$ & $-\dfrac{1}{21}$ & $\dfrac{1}{14}$ & $\dfrac{4}{105}$
  & $-\dfrac{3}{14}$ & $\dfrac{16}{231}$ & $\dfrac{13}{22}$\\
\hline
6 & $-\dfrac{1}{42}$ & $\dfrac{1}{42}$ & $\dfrac{8}{105}$
  & $-\dfrac{37}{210}$ & $-\dfrac{16}{231}$ & $\dfrac{305}{462}$\\
\hline
7 & $0$ & $\dfrac{1}{15}$ & $-\dfrac{1}{10}$
  & $-\dfrac{28}{165}$ & $\dfrac{13}{22}$\\
\hline
8 & $\dfrac{1}{30}$ & $-\dfrac{1}{30}$
  & $-\dfrac{32}{165}$ & $\dfrac{139}{330}$\\
\hline
9 & 0 & $-\dfrac{5}{33}$ & $\dfrac{5}{22}$\\
\hline
10 & $-\dfrac{5}{66}$ & $\dfrac{5}{66}$\\
\hline
11 & 0
\end{tabular} }
\bigskip

By Table~A.2 one can calculate the function $C(\la,\mu)$ up to
 degree~10, see Definition~2.4b.
\bigskip

\noindent
{\bf Example A.3}.
Up to degree 10 one has
 $C(\la,\mu)=-\dfrac{1}{2}+\dfrac{1}{12}(\la-\mu)+\dfrac{1}{24}\la\mu\;+$
$$+\left(\frac{1}{720}\mu^3+\frac{1}{180}\la\mu^2
         -\frac{1}{180}\la^2\mu-\frac{1}{720}\la^3\right)
  -\left(\frac{1}{1440}\la^3\mu+\frac{1}{360}\la^2\mu^2
        +\frac{1}{1440}\la\mu^3 \right)+$$
$$+\frac{1}{7!}\left( \frac{\la^5}{6}+\la^4\mu
        +\frac{4}{3}(\la^3\mu^2-\la^2\mu^3)
        -\la\mu^4-\frac{\mu^5}{6} \right)
  +\frac{1}{7!}\left( \frac{\la^5\mu}{12}+\frac{\la^4\mu^2}{2}
        +\frac{23}{24}\la^3\mu^3+\frac{\la^2\mu^4}{2}\
        +\frac{\la\mu^5}{12} \right)$$
$$+\frac{1}{7!}\left( \frac{1}{240}\mu^7+\frac{1}{30}\la\mu^6
        +\frac{4}{45}\la^2\mu^5+\frac{1}{15}\la^3\mu^4
        -\frac{1}{15}\la^4\mu^3-\frac{4}{45}\la^5\mu^2
        -\frac{1}{30}\la^6\mu-\frac{1}{240}\la^7 \right)-$$
$$-\frac{1}{8!}\left( \frac{1}{60}\la^7\mu+\frac{2}{15}\la^6\mu^2
        +\frac{37}{90}\la^5\mu^3+\frac{3}{5}\la^4\mu^4
        +\frac{37}{90}\la^3\mu^5+\frac{2}{15}\la^2\mu^6
        +\frac{1}{60}\la\mu^7 \right)+$$
$$+\frac{1}{11!}\left( \frac{\la^9}{6}+\frac{25}{3}\la^8\mu
        +32\la^7\mu^2+56\la^6\mu^3+32(\la^5\mu^4-\la^4\mu^5)
        -56\la^3\mu^6-32\la^2\mu^7-\frac{25}{3}\la\mu^8
        -\frac{\mu^9}{6} \right)$$
$$+\frac{\la\mu}{11!}\left( \frac{5}{12}(\la^8+\mu^8)
        +\frac{\la^7\mu+\la\mu^7}{24}
        +\frac{139}{8}(\la^6\mu^2+\la^2\mu^6)
        +39(\la^5\mu^3+\la^3\mu^5)
        +\frac{305}{6}\la^4\mu^4 \right).\eqno{\bt}$$
\medskip

\noindent
{\bf Proposition A.4}.
Let $L$ be the Lie algebra freely generated by the symbols $P,Q$.
Under $L\to L/\bigl[[L,L],[L,L]\bigr]$
 the Hausdorff series $H=\log\bigl(\exp(P)\cdot\exp(Q)\bigr)$
 maps onto
$$\bar H=(P+Q)
   +\frac{1}{2}[PQ]
   +\left(\frac{1}{12}[P^2Q]-\frac{1}{12}[QPQ]\right)
   -\frac{1}{24}[PQPQ]+$$
$$+\left( \frac{[Q^3PQ]-[P^4Q]}{720}+\frac{[PQ^2PQ]-[P^2QPQ]}{180} \right)
  +\left( \frac{[P^3QPQ]+[PQ^3PQ]}{1440}+\frac{[P^2Q^2PQ]}{360} \right)+$$
$$+\frac{1}{7!}\left( \frac{1}{6}[P^6Q]+[P^4QPQ]
  +\frac{4}{3}[P^3Q^2PQ]-\frac{4}{3}[P^2Q^3PQ]
  -[PQ^4PQ]-\frac{1}{6}[Q^5PQ]\right)$$
$$-\frac{1}{7!}\left( \frac{1}{12}[P^5QPQ]+\frac{1}{2}[P^4Q^2PQ]
  +\frac{23}{24}[P^3Q^3PQ]+\frac{1}{2}[P^2Q^4PQ]
  +\frac{1}{12}[PQ^5PQ]\right)$$
$$+\frac{1}{7!}\left( \frac{[P^8Q]-[Q^7PQ]}{240}
  +\frac{[P^6QPQ]-[PQ^6PQ]}{30}
  +\frac{4}{45}([P^5Q^2PQ]-[P^2Q^5PQ]) \right)$$
$$+\frac{[P^4Q^3PQ]-[P^3Q^4PQ]}{15\cdot 7!}
  -\frac{1}{8!}\left( \frac{[P^7QPQ]+[PQ^7PQ]}{60}
  +\frac{2}{15}([P^6Q^2PQ]+[P^2Q^6PQ]) \right)$$
$$-\frac{1}{8!}\left( \frac{37}{90}([P^5Q^3PQ]+[P^3Q^5PQ])
  -\frac{3}{5}[P^4Q^4PQ] \right)+\mbox{ (higher degree terms).}$$

\noindent
\begin{proof}
The coefficient of the term $[P^{n-1}Q^{m-1}PQ]$ in $\bar H$ coincides
 with the coefficient of the term $\la^{n-1}\mu^{m-1}$ in
 $C(\la,\mu)$.
Hence Proposition~A.4 follows from Example~A.3.
\end{proof}
\medskip

\noindent
{\bf Proposition A.5}.
(a) The even part of any compressed Drinfeld associator
 $\bar\ph(a,b)$ is
$$Even\bigl(\bar\ph(a,b)\bigr)
 =\dfrac{[ab]}{6}-\dfrac{4[a^3b]+[abab]+4[b^2ab]}{360}
  +\dfrac{[a^5b]+[b^4ab]}{945}+$$
$$+\Bigl( \be_{31}+\dfrac{20}{3\cdot 7!}\Bigl)
   \bigl([a^3bab]+[ab^3ab]\bigr)
  +\Bigl( 2\be_{31}+\dfrac{13}{7!}\Bigl)[a^2b^2ab]
  -\dfrac{[a^7b]+[b^6ab]}{9450}+$$
$$+\Bigl( \dfrac{\be_{31}}{6}+\be_{41}-\dfrac{1}{4200}\Bigl)
   \bigl([a^5bab]+[ab^5ab]\bigr)
  +\Bigl( \dfrac{2}{3}\be_{31}+3\be_{41}-\dfrac{113}{45\cdot 7!}\Bigl)
   \bigl([a^4b^2ab]+[a^2b^4ab]\bigr)+$$
$$+\Bigl( \be_{31}+4\be_{41}-\dfrac{947}{5\cdot 9!}\Bigl)
   [a^3b^3ab]
  +\dfrac{[a^9b]+[b^8ab]}{93555}+\cdots$$
\smallskip

\noindent
(b) Up to degree 9 the odd part of any compressed Drinfeld associator
 $\bar\ph(a,b)\in\hat{\bar L}_3$ is
$$Odd\bigl(\bar\ph(a,b)\bigr)
 =\ti\be_{00}\bigl([a^2b]+[bab]\bigr)
  +\Bigl(\ti\be_{10}+\dfrac{\ti\be_{00}}{6}\Bigr)
   \bigl([a^4b]+[b^3ab]\bigr)+$$
$$+\Bigl(2\ti\be_{10}+\dfrac{\ti\be_{00}}{2}\Bigr)
   \bigl([a^2bab]+[ab^2ab]\bigr)
  +\Bigl(\ti\be_{20}+\dfrac{\ti\be_{10}}{6}+\dfrac{\ti\be_{00}}{120}\Bigr)
   \bigl([a^6b]+[b^5ab]\bigr)+$$
$$+\Bigl(3\ti\be_{20}+\dfrac{2}{3}\ti\be_{10}+\dfrac{\ti\be_{00}}{24}\Bigr)
   \bigl([a^4bab]+[ab^4ab]\bigr)
  +\Bigl(5\ti\be_{20}+\dfrac{7}{6}\ti\be_{10}+\dfrac{\ti\be_{00}}{12}\Bigr)
   \bigl([a^3b^2ab]+[a^2b^3ab]\bigr)+$$
$$+\Bigl(\ti\be_{30}+\dfrac{\ti\be_{20}}{6}
    +\dfrac{\ti\be_{10}}{120}+\dfrac{\ti\be_{00}}{7!}\Bigr)
   \bigl([a^8b]+[b^7ab]\bigr)
  +\Bigl(4\ti\be_{30}+\dfrac{5}{6}\ti\be_{20}
    +\dfrac{\ti\be_{10}}{20}+\dfrac{\ti\be_{00}}{720}\Bigr)
   \bigl([a^6bab]+[ab^6ab]\bigr)+$$
$$+\Bigl(\ti\be_{31}+9\ti\be_{30}+2\ti\be_{20}
    +\dfrac{2}{15}\ti\be_{10}+\dfrac{\ti\be_{00}}{240}\Bigr)
   \bigl([a^5b^2ab]+[a^2b^5ab]\bigr)$$
$$+\Bigl(3\ti\be_{31}+13\ti\be_{30}+3\ti\be_{20}
    +\dfrac{5}{24}\ti\be_{10}+\dfrac{\ti\be_{00}}{144}\Bigr)
   \bigl([a^4b^3ab]+[a^3b^4ab]\bigr)
   +\mbox{(higher degree terms)}.$$

\noindent
\emph{Proof} follows from Theorem~1.5c by using routine computations.
The degree 6 part of the series $\bar\ph^B(a,b)$ from Example~1.4
 is obtained from A.5(a) for $\be_{31}=-\dfrac{8}{3\cdot 7!}$.
\qed
\bigskip

\noindent
{\bf Claim A.6}.
(a) The even part of the Drinfeld series $f^D(\la,\mu)$ from
 Theorem~6.2 starts with
{\small
$$Even(f^D(\la,\mu))=\dfrac{1}{6}-\dfrac{\la^2+\mu^2}{90}
 -\dfrac{\la\mu}{360}+\dfrac{\la^4+\mu^4}{945}
 +\Bigl(\dfrac{9}{2}\te_3^2+\dfrac{1}{1260}\Bigr)(\la^3\mu+\la\mu^3)
  +\Bigl(9\te_3^2+\dfrac{23}{3\cdot 7!}\Bigr)\la^2\mu^2
 -\dfrac{\la^6+\mu^6}{9450}$$
$$+\Bigl(15\te_3\te_5-\dfrac{2}{3\cdot 7!}\Bigr)(\la^5\mu+\la\mu^5)
  +\Bigl(45\te_3\te_5+\dfrac{3}{4}\te_3^2-\dfrac{61}{45\cdot 7!}\Bigr)
   (\la^4\mu^2+\la^2\mu^4)
  +\Bigl(60\te_3\te_5+\dfrac{3}{2}\te_3^2-\dfrac{499}{5\cdot 9!}\Bigr)
   \la^3\mu^3+\cdots$$ }
\smallskip

\noindent
(b) Up to degree 7 the odd part of the Drinfeld series $f^D(\la,\mu)$ from
 Theorem~6.2 is
$$Odd(f^D(\la,\mu))=-3\te_3(\la+\mu)-5\te_5(\la^3+\mu^3)
  -\Bigl(10\te_5+\dfrac{\te_3}{2}\Bigr)(\la^2\mu+\la\mu^2)
  -7\te_7(\la^5+\mu^5)$$
$$-\Bigl(21\te_7+\dfrac{5}{6}\te_5-\dfrac{\te_3}{30}\Bigr)
   (\la^4\mu+\la\mu^4)
  -\Bigl(35\te_7+\dfrac{5}{3}\te_5-\dfrac{\te_3}{24}\Bigr)
   (\la^3\mu^2+\la^2\mu^3)
  -9\te_9(\la^7+\mu^7)-$$
$$-\bigl(36\te_9+\dfrac{7}{6}\te_7-\dfrac{\te_5}{18}+\dfrac{\te_3}{315}\bigr)
   (\la^6\mu+\la\mu^6)
  -\Bigl(\dfrac{9}{2}\te_3^3+84\te_9+\dfrac{7}{2}\te_7-\dfrac{\te_5}{8}
   +\dfrac{\te_3}{180}\Bigr)(\la^5\mu^2+\la^2\mu^5)-$$
$$-\Bigl(\dfrac{27}{2}\te_3^3+126\te_9+\dfrac{35}{6}\te_7-\dfrac{7}{36}\te_5
   +\dfrac{\te_3}{144}\Bigr)(\la^4\mu^3+\la^3\mu^4)
   +(\mbox{higher degree terms})$$
\begin{proof}
Rewrite the formula~(6.2) in a more explicit form:
{\small
$$f^D(\la,\mu)=-2\te_2-3\te_3(\la+\mu)-4\te_4(\la^2+\mu^2)
  +\Bigl(2\te_2^2-6\te_4\Bigr)\la\mu-5\te_5(\la^3+\mu^3)
  +\Bigl(6\te_2\te_3-10\te_5\Bigr)(\la^2\mu+\la\mu^2)$$
$$-6\te_6(\la^4+\mu^4)
  +\Bigl(\dfrac{9}{2}\te_3^2+8\te_2\te_4-15\te_6\Bigr)(\la^3\mu+\la\mu^3)
  +\Bigl(-\dfrac{4}{3}\te_2^3+9\te_3^2+12\te_2\te_4-20\te_6\Bigr)\la^2\mu^2
  -7\te_7(\la^5+\mu^5)$$
$$+\Bigl(10\te_2\te_5+12\te_3\te_4-21\te_7\Bigr)(\la^4\mu+\la\mu^4)
  +\Bigl(-6\te_2^2\te_3+20\te_2\te_5+30\te_3\te_4-35\te_7\Bigr)
   (\la^3\mu^2+\la^2\mu^3)
  -8\te_8(\la^6+\mu^6)$$
$$+\Bigl(12\te_2\te_6+15\te_3\te_5+8\te_4^2-28\te_8\Bigr)(\la^5\mu+\la\mu^5)
  +\Bigl(30\te_2\te_6+45\te_3\te_5+24\te_4^2-9\te_2\te_3^2-8\te_2^2\te_4
   -58\te_8\Bigr)(\la^4\mu^2+\la^2\mu^4)$$
$$\Bigl(\dfrac{2}{3}\te_2^4-18\te_2^2\te_3-12\te_2^2\te_4
   +40\te_2\te_6+60\te_3\te_5+34\te_4^2-70\te_8\Bigr)\la^3\mu^3
  +\Bigl( 14\te_2\te_7+18\te_3\te_6+20\te_4\te_5-36\te_9 \Bigr)
   (\la^6\mu+\la\mu^6)$$
$$-9\te_9(\la^7+\mu^7)
  +\Bigl( -24\te_2\te_3\te_4-\dfrac{9}{2}\te_3^3-10\te_2^2\te_5+42\te_2\te_7
   +63\te_3\te_6+70\te_4\te_5-84\te_9 \Bigr)(\la^5\mu^2+\la^2\mu^5)+$$
$$+\Bigl( 4\te_2^3\te_3-60\te_2\te_3\te_4-\dfrac{27}{2}\te_3^3-20\te_2^2\te_5
   +70\te_2\te_7+105\te_3\te_6+120\te_4\te_5-126\te_9 \Bigr)
   (\la^4\mu^3+\la^3\mu^4)+\cdots$$ }
It remains to substitute $\te_{2n}$ from Example~6.4
 and split the even and odd parts.
Observe that the Drinfeld series $f^D(\la,\mu)$ is obtained from
 Proposition~A.5 for the parameters
$$\be_{31}=\dfrac{9}{2}\te_3^2-\dfrac{8}{3\cdot 7!},\quad
  \be_{41}=15\te_3\te_5-\dfrac{3}{4}\te_3^2+\dfrac{44}{45\cdot 7!},\quad
  \ti\be_{00}=-3\te_3,\quad
  \ti\be_{10}=-5\te_5+\dfrac{\te_3}{2},$$
$$\ti\be_{20}=-7\te_7+\dfrac{5}{6}\te_5-\dfrac{7}{120}\te_3,\quad
  \ti\be_{30}=-9\te_9+\dfrac{7}{6}\te_7-\dfrac{7}{72}\te_5
   +\dfrac{31}{7!}\te_3,\quad
  \ti\be_{31}=-\dfrac{9}{2}\te_3^3-3\te_9+\dfrac{\te_3}{630}.$$
The above formulae reprove explicitly Proposition~6.9 up to
 degree 7.
\end{proof}

\noindent
{\bf Example A.7}.
The first and second distiguished even compressed Drinfeld associators
 from Corollary~1.6a--b of Subsection~1.3 are defined by the following series:
$$(a)\quad
  f^I(\la,\mu)=\dfrac{1}{6}-\dfrac{4\la^2+\la\mu+4\mu^2}{360}
  +\dfrac{\la^4+\mu^4}{945}
  +\dfrac{20}{4\cdot 7!}\bigl(\la^3\mu+\la\mu^3\bigr)
  +\dfrac{13}{7!}\la^2\mu^2-$$
$$-\dfrac{\la^6+\mu^6}{9450}
  -\dfrac{\la^5\mu+\la\mu^5}{4200}
  -\dfrac{113}{45\cdot 7!}\bigl(\la^4\mu^2+\la^2\mu^4\bigr)
  -\dfrac{947}{5\cdot 9!}\la^3\mu^3
  +\dfrac{\la^8+\mu^8}{93555}+\cdots$$

$$(b)\quad
  f^{II}(\la,\mu)=\dfrac{1}{6}-\dfrac{4\la^2+\la\mu+4\mu^2}{360}
  +\dfrac{\la^4+\mu^4}{945}
  -\dfrac{53}{6\cdot 7!}\bigl(\la^3\mu+\la\mu^3\bigr)
  -\dfrac{18}{7!}\la^2\mu^2-$$
$$-\dfrac{\la^6+\mu^6}{9450}
  +\dfrac{\la^5\mu+\la\mu^5}{11200}
  -\dfrac{13}{90\cdot 7!}\bigl(\la^4\mu^2+\la^2\mu^4\bigr)
  -\dfrac{431}{5\cdot 9!}\la^3\mu^3
  +\dfrac{\la^8+\mu^8}{93555}+\cdots$$
\smallskip

\noindent
The proof is a computation following from
 Corollary~1.6.
Both series are obtained from
 Proposition~A.5 for $\be_{31}=\be_{41}=0$ and
 $\be_{31}=-\dfrac{31}{2\cdot 7!}$,
 $\be_{41}=\dfrac{127}{30\cdot 7!}$, respectively.
\hfill $\bt$



\begin{thebibliography}{References}

\bibitem{Bak}
H.~Baker,
\emph{On a law of combination of operators (second paper)},
Proceedings London Math. Soc. {\bf 1} (1898), N 29, p.~14--32.

\bibitem{BN1}
D.~Bar-Natan,
\emph{On the Vassiliev knot invariants},
Topology {\bf 34} (1995), p.~423--472.

\bibitem{BN2}
D.~Bar-Natan,
\emph{Non-associative tangles},
in \emph{Geometric Topology}
 (proceedings of the Georgia International Topology Conference),
 (W.~H.~Kazez, ed.), p.~139--183, Providence, 1997.

\bibitem{BN3}
D.~Bar-Natan,
\emph{On associators and the Grothendieck-Teichmuller group I},
Selecta Mathematica, New Series {\bf 4} (1998), p.~183--212.

\bibitem{BLT}
D.~Bar-Natan, T.~Q.~T.~Le, D.~P.~Thurston,
\emph{Two applications of elementary knot theory to Lie algebras
 and Vassiliev invariants},
Geometry and Topology {\bf 7} (2003), p.~1--31.

\bibitem{Cam}
J.~Campbell,
\emph{Alternants and continuous groups},
Proc. London Math. Soc. {\bf 2(3)} (1905), p.~24--47.

\bibitem{Car1}
P.~Cartier,
\emph{Construction combinatoire des invariantes de
 Vassiliev-Kontsevich des noeuds},
C.~R.~Acad.~Sci.~Paris {\bf 316}, S\'erie I (1993), p.~1205--1210.

\bibitem{Car2}
P.~Cartier,
\emph{Fonctions polylogarithmes, nombres polyz\^etas et groupes
 pro-unipotents},
S\'eminaire Bourbaki, v.~2000/2001, Ast\'erisque No.~282 (2002),
 Exp. No.~885, viii, p.~137--173.

\bibitem{Dr1}
V.~Drinfeld,
\emph{On quasi-Hopf algebras},
Leningrad Math. J. {\bf 1} (1990), p.~1419--1457.

\bibitem{Dr2}
V.~Drinfeld,
\emph{On quasi-triangular quasi-Hopf algebras and a group closely
 related with $\mathrm{Gal}(\bar\Q/Q)$},
Leningrad Math. J. {\bf 2} (1990), p.~829--860.

\bibitem{Hau}
F.~Hausdorff,
\emph{Die symbolische exponentialformel in der gruppentheorie},
Leipziger Berichte {\bf 58} (1906), p.~19--48.

\bibitem{Kas}
C.~Kassel,
\emph{Quantum groups},
Springer-Verlag, Graduate Texts in Math. {\bf 155},
New York, 1994.

\bibitem{Kon}
M.~Kontsevich,
Vassiliev's knot invariants,
Adv. in Soviet Math., {\bf 16} (1993), N 2, p.~137--150.

\bibitem{LM1}
T.~Q.~T.~Le, J.~Murakami,
\emph{Representations of the category of tangles by Kontsevich's
 iterated integral},
Commun. Math. Physics {\bf 168} (1995), p.~535--562.

\bibitem{LM2}
T.~Q.~T.~Le, J.~Murakami,
\emph{The universal Vassiliev-Kontsevich invariant for framed oriented links},
Composition Math. {\bf 102} (1996), p.~42--64, see also
 hep-th/9401016.

\bibitem{LM3}
T.~Q.~T.~Le, J.~Murakami,
\emph{Parallel version of the universal Vassiliev-Kontsevich invariant},\\
Journal of Pure and Applied Algebra {\bf 121} (1997), p.~271--291.

\bibitem{Lie}
J.~Lieberum,
\emph{The Drinfeld associator of $\mathrm{gl}(1|1)$},
arXiv:math.QA/0204346.


\bibitem{Oht}
\emph{Problems on invariants of knots and 3-manifolds},
edited by T.~Ohtsuki, T.~Kohno, T.~Le, J.~Murakami, V.~Turaev,
Geometry and Topology Monographs {\bf 4} (2002), p.~377--572.

\bibitem{Piu}
S.~Piunikhin,
\emph{Combinatorial expression for universal Vassiliev link invariant},
Commun. Math. Physics {\bf 168} (1995), p.~1--22.

\bibitem{Reu}
C.~Reutenauer,
\emph{Free Lie algebras},
London Math. Soc. Monographs, New Series, v.~7, 1993.

\bibitem{Vas}
V.~Vassiliev,
\emph{Cohomology of knot spaces},
Theory of singularities and its applications (V.~Arnold, ed.),
 Amer. Math. Soc., Providence, 1990.

\end{thebibliography}
\end{document}